\documentclass[12pt, a4paper, titlepage]{article}

\usepackage{latexsym}
\usepackage{color}

\usepackage{amssymb}
\usepackage[centertags]{amsmath}
\usepackage{verbatim}
\usepackage{epic, eepic}

\usepackage{graphicx}

\title{Transportation to random zeroes\\ by the gradient flow}

\author{Fedor Nazarov\thanks{Partially
supported by the National Science Foundation, DMS grant 0501067},
Mikhail Sodin\thanks{Partially supported by the Israel Science
Foundation of the Israel Academy of Sciences and Humanities, grant
357/04}, \ Alexander Volberg\thanks{Partially supported by the
National Science Foundation, DMS grant 0501067}}

\date{}

\begin{document}

\newtheorem{theorem}[equation]{Theorem}
\newtheorem{lemma}[equation]{Lemma}
\newtheorem{question}[equation]{Question}
\newtheorem{proposition}[equation]{Proposition}
\newtheorem{corollary}[equation]{Corollary}
\newtheorem{conjecture}[equation]{Conjecture}
\newtheorem{definition}[equation]{Definition}
\newtheorem{problem}[equation]{Problem}
\newtheorem{note}[equation]{Note}
\newtheorem{myItem}[equation]{}
\newtheorem{claim}[equation]{Claim}

\numberwithin{equation}{section}

\newcommand{\tendsk}{\xrightarrow[k\to\infty]{}}

\renewcommand{\phi}{\varphi}

\renewcommand{\(}{\bigl(}
\renewcommand{\)}{\bigr)\vphantom{)}}

\newcommand{\N}{\operatorname{N}}
\newcommand{\mes}{\operatorname{mes}}
\newcommand{\const}{\operatorname{const}}
\newcommand{\Prob}{\operatorname{Prob}}
\newcommand{\Var}{\operatorname{Var}}
\newcommand{\Homeo}{\operatorname{Homeo}}
\newcommand{\entier}{\operatorname{entier}}
\newcommand{\Aut}{\operatorname{Aut}}
\newcommand{\Iso}{\operatorname{Iso}}
\newcommand{\MALG}{\operatorname{MALG}}
\newcommand{\dist}{\operatorname{dist}}
\renewcommand{\Re}{\operatorname{Re}}
\renewcommand{\Im}{\operatorname{Im}}
\newcommand{\spec}{\operatorname{spec}}
\newcommand{\sign}{\operatorname{sgn}}
\newcommand{\Tra}{\operatorname{Tra}}
\newcommand{\Ra}{\operatorname{Ra}}
\newcommand{\Di}{\operatorname{Di}}
\newcommand{\Const}{\operatorname{Const}}

\renewcommand{\div}{\operatorname{div}}

\newcommand{\pd}{\partial}

\newcommand{\SF}{{\text{S$\sigma$F}}}
\newcommand{\K}{\mathrm K}
\newcommand{\bF}{\mathrm F}
\newcommand{\bL}{\mathrm L}
\newcommand{\U}{\mathrm U}
\newcommand{\bN}{\mathbf N}
\newcommand{\Om}{\Omega}
\newcommand{\om}{\omega}
\newcommand{\al}{\alpha}
\newcommand{\be}{\beta}
\newcommand{\eps}{\varepsilon}
\newcommand{\ga}{\gamma}
\newcommand{\si}{\sigma}
\newcommand{\la}{\lambda}
\newcommand{\de}{\delta}
\newcommand{\De}{\Delta}
\newcommand{\F}{\mathcal F}
\newcommand{\E}{\mathcal E}
\newcommand{\X}{\mathcal X}
\newcommand{\A}{\mathcal A}
\newcommand{\B}{\mathcal B}
\newcommand{\cC}{\mathcal C}
\newcommand{\M}{\mathcal M}
\newcommand{\cQ}{\mathcal Q}
\newcommand{\fav}{{\text{fav}}}
\newcommand{\Omfav}{\Om^\fav}
\newcommand{\Ffav}{\F^\fav}
\newcommand{\Pfav}{P^\fav}
\newcommand{\Mfav}{\M^\fav}
\newcommand{\Gfav}{G^\fav}
\newcommand{\Xcna}{X^{\text{cna}}}
\newcommand{\Xf}{{X^{\text{f}}}}

\newcommand{\Q}{\partial Q}

\renewcommand{\epsilon}{\varepsilon}

\def\done{{1\hskip-2.5pt{\rm l}}}

\newcommand{\Ex}{\mathbb E}
\renewcommand{\Pr}[1]{\,\mathbb P\,\{\,#1\,\}\,}
\renewcommand{\le}{\leqslant}
\renewcommand{\ge}{\geqslant}

\newcommand{\R}{\mathbb R}
\newcommand{\C}{\mathbb C}
\newcommand{\Z}{\mathbb Z}
\newcommand{\D}{\mathbb D}

\newcommand{\T}{\mathbb T}
\newcommand{\bB}{\mathbb B}
\newcommand{\cE}[2]{\mathbb{E}\,\(\,#1\,\big|\,#2\,\)\,}
\newcommand{\scE}[2]{\,\mathbb{E}(#1|#2)}
\newcommand{\CE}[2]{\mathbb{E}\,\bigg(#1\,\bigg|\,#2\bigg)\,}
\newcommand{\cP}[2]{\mathbb{P}\,\(\,#1\,\big|\,#2\,\)\,}

\newcommand{\ti}{\widetilde}
\newcommand{\crit}{\operatorname{Crit} U}

\newcommand{\sif}{$\sigma$-field}

\newcommand{\Wi}[1]{\rule{0pt}{1pt}\boldsymbol{:}#1\boldsymbol{:}\rule{0pt}{1pt}}
\newcommand{\Zero}{\mathbf0}
\newcommand{\One}{\mathbf1}
\newcommand{\imply}{\;\;\;\Longrightarrow\;\;\;}
\newcommand{\rimply}{\;\;\;\Longleftarrow\;\;\;}
\newcommand{\imp}{$ \Longrightarrow $ }
\newcommand{\impl}{\;\Longrightarrow\;}
\renewcommand{\equiv}{\;\;\;\Longleftrightarrow\;\;\;}
\newcommand{\equi}{\;\Longleftrightarrow\;}
\newcommand{\equ}{$ \;\; \Longleftrightarrow \;\; $}

\def\emailwww#1#2{\par\qquad {\tt #1}\par\qquad {\tt #2}\medskip}
\def\One{{1\hskip-2.5pt{\rm l}}}

\newenvironment{myitemize}{\begin{list}{$\bullet$}
{\setlength{\topsep}{1mm}
\setlength{\partopsep}{0mm}
\setlength{\itemsep}{1mm}
\setlength{\parsep}{0mm}
\setlength{\parskip}{0mm}}}
{\end{list}}

\maketitle

\begin{abstract}
We consider the zeroes of the random Gaussian entire function
\[
f(z) = \sum_{k=0}^\infty \xi_k \frac{z^k}{\sqrt{k!}}
\]
($\xi_0, \xi_1, \dots $ are Gaussian i.i.d. complex random
variables) and show that their basins under the gradient flow of
the random potential $U(z)=\log|f(z)| - \frac12 |z|^2$ partition
the complex plane into domains of equal area.

We find three characteristic exponents $1$, $\frac85$, and $4$ of
this random partition: the probability that the diameter of a
particular basin is greater than $R$ is exponentially small in
$R$; the probability that a given point $z$ lies at a distance
larger than $R$ from the zero it is attracted to decays as
$e^{-R^{8/5}}$; and the probability that, after throwing away
$1\%$ of the area of the basin, its diameter is still larger than
$R$ decays as $e^{-R^4}$.

We also introduce a combinatorial procedure that modifies a small
portion of each basin in such a way that the probability that the
diameter of a particular modified basin is greater than $R$
decays as $e^{-cR^4(\log R)^{-3/2}}$.

\end{abstract}

\tableofcontents

\section{Introduction and main results}

Let $\mathcal Z$ be a random point process in $\R^d$ with the
distribution invariant with respect to the isometries of $\R^d$.
Suppose that $\mathcal Z$ has intensity $1$; that is, the mean
number of points of $\mathcal Z$ per unit volume equals $1$. The
{\em transportation} (a.k.a.  ``matching'', ``allocation'',
``marriage'', etc.) of the Lebesgue measure $m_d$ to $\mathcal Z$
is a (random) measurable map $T\colon \R^d\to \mathcal Z$ that
pushes forward the Lebesgue measure $m_d$ to the counting measure
$\displaystyle n_\mathcal Z = \sum_{a\in \mathcal Z} \delta_a$ of
the set $\mathcal Z$ ( $\delta_a$ is the unit mass at $a$). In
other words, the whole space $\R^d$ is split into disjoint random
sets $B(a)$ of the Lebesgue measure $1$ indexed by $a\in \mathcal
Z$. Because of the invariance of the process $\mathcal Z$, it is
natural to assume that the transportation $T$ has an invariant
distribution; i.e., that the distribution of the vector $T(x)-x$
does not depend on $x$.  The better $T$ is localized, the more
uniformly the process $\mathcal Z$ is spread over $\R^d$. Thus it
is interesting to know the optimal rate of decay of the
probability tails $\Pr {|T(x)-x|>R}$ as $R\to\infty$. A
constructive counterpart is to find an {\em explicit} and
well-localized way to transport the Lebesgue measure $m_d$ to the
point process $\mathcal Z$.

\medskip

The transportation of the Lebesgue measure $m_d$ to the Poisson
process in $\R^d$ was recently developed by Hoffman, Holroyd and
Peres \cite{HHP1, HHP2} (a finite volume version was studied
earlier by Ajtai, Koml\'os and Tusn\'ady~\cite{AKT}, Leighton and
Shor~\cite{LS}, and Talagrand~\cite{T}). In this paper, we
consider the random zero point set $\mathcal Z_f = f^{-1}(0)$ of a
Gaussian entire function $f$ in $\C$ and study the transportation
of the two-dimensional Lebesgue measure $m_2$ to $\mathcal Z_f$.

\medskip Let
\[
f(z) = \sum_{k\ge 0}\xi_k\frac{z^k}{\sqrt{k!}}
\]
where $\xi_k$ are independent standard complex Gaussian random
variables (i.e., the density of $\xi_k$ on the complex plane $\C$
is $\frac1\pi e^{-|z|^2}$). We shall call such a random function a
Gaussian Entire Function (G.E.F.).

The (random) zero set $\mathcal Z_f$ of this function is known as
``flat chaotic analytic zero points" \cite{Han1, Han2, Leb}. It is
distinguished by the invariance of its distribution with respect
to the isometries of $\C$; i.e., rotations and translations, see
\cite[Part~I]{ST} for details and references. Note that the
intensity of the zero process $\mathcal Z_f $ equals
$\frac1{\pi}$. In \cite[part~II]{ST}, the question about the
existence of a well-localized transportation of the area measure
to the zero set of the Gaussian Entire Function in $\mathbb C$ was
studied. Using the Hall matching lemma and some potential theory,
the authors of \cite{ST} proved the existence of a transportation
with sub-Gaussian decay of the tail probability. Unfortunately,
the proof one obtains on this way is a pure existence proof giving
no idea of what the transportation in question looks like.

\medskip

The aim of this paper is to carry out another approach that was
suggested but not followed in \cite[part~II]{ST}, namely, the
transportation by the gradient flow of a random potential. The
main advantage of this approach is that it provides a quite
natural and explicit construction for the desired transportation.

\medskip
Let $U(z)=\log|f(z)|-\frac12{|z|^2}$ be the random potential
corresponding to the G.E.F. $f$. The distribution of $U$ is also
invariant with respect to the isometries of the complex plane, see
\cite[part~I]{ST} or Section~2.2 below. We shall call any integral
curve of the differential equation
\[
\frac{dZ}{dt} = -\nabla U(Z)
\]
a {\em gradient curve} of the random potential $U$.

We orient the gradient curves in the direction of decay of $U$
(this is the reason for our choice of the minus sign in the
differential equation above). If $z\notin \mathcal Z_f$, and
$\nabla U(z)\ne 0$, by $\Gamma_z$ we denote the (unique) gradient
curve that passes through the point $z$.

\begin{definition}[the basin] Let $a$ be a zero of the G.E.F. $f$.
The {\it basin} of $a$ is the set
$$
B(a)=\{z\in \C\colon \,\nabla U(z)\ne 0,\text{ and } \Gamma_z
\text{ terminates at }a\}\,.
$$
\end{definition}

The picture below may help the reader to visualize this
definition. It shows the random zeroes and the trajectories of
various points under the gradient flow.

\begin{figure}[h]
\hskip4cm \scalebox{0.3}{\includegraphics{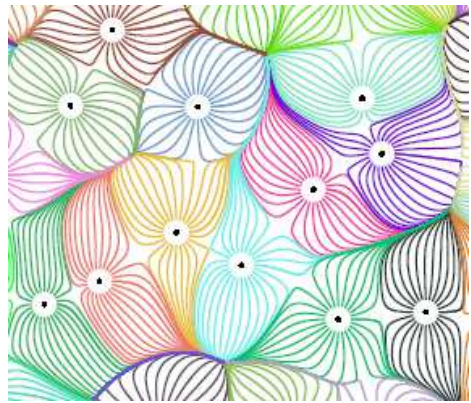}} \caption{The
basins $B(a)$ and trajectories of the gradient field}
\end{figure}

Clearly, each basin $B(a)$ is a connected open set, and $B(a')\cap
B(a'')=\varnothing$ if $a'$ and $a''$ are two different zeroes of
$f$. If the basin $B(a)$ is bounded and the boundary of $B(a)$ is
nice, then $\frac{\partial U}{\partial n} = 0$ on $\partial B(a)$
and therefore, applying the Green formula and observing that the
distributional Laplacian of $U$ equals $\Delta U = 2\pi
\sum_{a\in\mathcal Z_f} \delta_a - 2m_2$, one gets
$$
1-\frac {m_2 B(a)}{\pi}=\frac 1{2\pi}\iint_{B(a)}\Delta
U(z)\,dm_2(z) = \frac1{2\pi}\int_{\partial B(a)}\frac {\partial
U}{\partial n}(z)\,|dz|=0\,;
$$
i.e., $m_2 B(a) = \pi$.

\medskip
Now we are ready to formulate our main results:
\begin{theorem}[partition]\label{th.partition}
Almost surely, each basin is bounded by finitely many smooth
gradient curves (and, thereby, has area $\pi$), and
\[
\C=\bigcup_{a\in\mathcal Z_f} B(a)
\]
up to a set of measure $0$ (more precisely, up to countably many
smooth boundary curves)\,.
\end{theorem}

Consider the random set
\[
S=\bigcup_{a\in\mathcal Z_f} \partial B(a)\,;
\]
that is, the union of all ``singular'' gradient curves; i.e., the
curves that do not terminate at $\mathcal Z_f$. Due to the
translation invariance of the random potential $U$, the
probability $ \Pr{ z_0\in S } $ does not depend on the choice of
the point $z_0\in\C$, hence vanishes:
\[
\Pr{ 0\in S } = \frac1{\pi} \iint_{\mathbb D} \Pr{ z\in S }\,
dm_2(z) = \int_\Omega m_2(S\cap \mathbb D)\, d\mathbb P = 0
\]
(here $\Omega$ is the probability space and $\mathbb D$ is the
unit disk). Thus, almost surely, any given point $z\in\C$ belongs
to some basin.

By $B_z$ we denote the basin that contains the point $z$. By
$\operatorname{diam}(A)$ we denote the diameter of a set
$A\subset\C$. We denote by $C$ and $c$ absolute (numerical)
constants that may change from one line to another.
\begin{theorem}[diameter of the basin]\label{th.diameter}
For any point $z\in \C$ and any $R\ge 1$,
\[
ce^{-CR(\log R)^{3/2}} \le \Pr{\operatorname{diam}(B_z) > R} \le
Ce^{-cR(\log R)^{3/2}}\,.
\]
\end{theorem}

\medskip
The proof of Theorem~\ref{th.diameter} relies on the following
auxiliary theorem. Let $Q(w,s)$ be the square centered at $w$ with
side length $2s$ and let $\Q(w,s)$ be its boundary.

\begin{theorem}[long gradient curve]\label{th.longcurve}
Let $ R\ge 1 $. The probability of the event that there exists a
gradient curve joining $\Q(0,R)$ with $\Q(0,2R)$ does not exceed
$Ce^{-cR(\log R)^{3/2}}$.
\end{theorem}

\medskip
The proof of this theorem is, unfortunately, quite involved. For a
weaker upper bound $Ce^{-cR\sqrt{\log R}}$ that has a simpler
proof, see the first version of this work posted in the {\tt
arxiv} \cite{NSV}. The approach in \cite{NSV} may be more suitable
for extensions to point processes of different nature: recently,
using a similar approach, Chatterjee, Peled, Peres, and Romik
found counterparts of Theorems~\ref{th.diameter} and
\ref{th.longcurve} for the Poisson process in $\R^d$ with $d\ge 3$
\cite{CPPR}. It might be helpful for the reader to look at
\cite{NSV} prior to reading the proof of the long gradient curve
theorem given here.

\medskip
Let $a_z$ be the random zero whose basin contains a given point
$z\in\C$. In other words, the gradient curve $\Gamma_z$ terminates
at $a_z$. It appears that the probability $ \Pr{ |z-a_z| > R }$ is
much smaller than the probability $ \Pr{\operatorname{diam}(B_z) >
R}$:
\begin{theorem}[distance to the sink]\label{th.distance}
For any point $ z\in\C $ and any $ R\ge 1 $,
\[
ce^{-CR^{8/5}} \le \Pr{ |z - a_z| > R} \le Ce^{-cR^{8/5}}\,.
\]
\end{theorem}

\medskip
This is related to long, thin ``tentacles'' seen on the picture
around some basins. They increase the typical diameter of basins
though the probability that a given point $z$ lies in such a
tentacle is very small.

Let $D(w, r)$ be the disk of radius $r$ centered at $w$.
\begin{theorem}[diameter of the core]\label{th.almost}
For any $z\in\C$, any $\epsilon>0$, and any $ R\ge 1$,
\[
c(\epsilon) e^{ -C(\epsilon) R^4} \le \Pr { m_2\left( B_z
\setminus D(a_z, R) \right) \ge \epsilon } \le C(\epsilon) e^{
-c(\epsilon) R^4}\,.
\]
Here $c(\epsilon)$, $C(\epsilon)$ are positive constants that
depend only on $\epsilon$.
\end{theorem}

\medskip
There exists a combinatorial procedure that allows one to cut the
tentacles off and to get an almost optimal estimate for the
diameters of the modified basins.
\begin{theorem}[modified basins]\label{th.cutoff}
Given $\epsilon>0$, there exist open pairwise disjoint sets
$B'(a)$ with the following properties:

\smallskip\par\noindent (i) $m_2 B'(a) = \pi$;

\smallskip\par\noindent (ii) $\displaystyle \C =
\bigcup_{a\in \mathcal Z_f} B'(a)$ (up to a set of measure $0$);

\smallskip\par\noindent (iii) $m_2\left(B(a) \bigcap B'(a) \right)
\ge \pi-\epsilon$;

\smallskip\par\noindent (iv) for any $z\in\C$, and any $R \ge 2$,
\[
\Pr {\operatorname{diam}(B')_z > R } \le C(\epsilon) e^{-cR^4(\log
R)^{-3/2}}.
\]
Here $(B')_z$ is the modified basin that contains the point $z$.
\end{theorem}

\medskip
The estimate in item (iv) is not as good as the tail estimate
$e^{-cR^4 (\log R)^{-1}}$ that can be obtained by modification of
the proof in \cite[part~II]{ST}, but it comes fairly close.

\medskip Now, a few words about the tools we use in the proofs.
First of all, it is the  ``almost independence" of the
localizations of a G.E.F. to distant disks (Theorem~\ref{thm3.1}),
which may be useful in other problems as well. In the proof of the
long gradient curve theorem, we use lower bounds for the
determinants of large covariance matrices of some Gaussian complex
random variables. These bounds are proved in Section~5. The proofs
of the distance to the sink theorem~\ref{th.distance}, the
diameter of the core theorem~\ref{th.almost}, and the modified
basins theorem~\ref{th.cutoff} are based on a version of the
length and area principle (Proposition~\ref{lemma8.a}).

\subsubsection*{Acknowledgment} Boris Tsirelson suggested
the idea of the transportation by the gradient flow. Manjunath
Krishnapur kindly provided us with inspiring computer generated
pictures of this transportation (including the one put in the
introduction). Yuval Peres helped us with the proof of
Lemma~\ref{lemma.16}. Leonid Polterovich helped us with
presentation in Section~8.1.  Bernie Shiffman explained to us a
connection with the statistics of critical points computed in
\cite{DSZ}. We thank all of them for numerous helpful discussions.
We thank Ron Peled and the referee for reading the paper carefully
and suggesting a number of corrections.

The first version of this paper was written while the second named
author was visiting Michigan State University and University of
California Berkeley in Fall 2005. He thanks both these
institutions for their generous hospitality.

\section{Preliminaries}

\subsection{Basic facts about complex Gaussian random variables}

We fix some probability space $(\Omega,\mathbb P)$ and some (very
big) family $\{\Xi_j\}_{j\in J}$ of  independent standard complex
Gaussian random variables on that probability space (i.e., the
density of $\Xi_j$ on the complex plane $\C$ is $\frac1\pi
e^{-|z|^2}$). Every complex Gaussian random variable in this paper
will be just a (possibly infinite) linear combination of $\Xi_j$
with square summable coefficients. Such a complex Gaussian random
variable $\xi$ is standard if $\Ex|\xi|^2=1$.

A useful remark is that if $\eta_k$ are standard complex Gaussian
random variables and $a_k\in \C$ satisfy $\sum_k|a_k|<+\infty$,
then $\sum_k a_k\eta_k$ can be represented as $a\eta$ where $0\le
a\le\sum_k|a_k|$ and $\eta$ is some standard Gaussian random
variable.

We shall start with simple probabilistic estimates.

\begin{lemma}\label{lemma2.1} Let $\eta_k$ be standard
complex Gaussian random variables (not necessarily independent).
Let $a_k>0$, $a=\sum_k{a_k}$. Then, for every $t>0$,
\[
\Pr {\sum_{k}a_k|\eta_k|>t}  \le 2e^{-\frac 12a^{-2}t^2}\,.
\]
\end{lemma}
\par\noindent{\em Proof:} Without loss of generality, $a=1$. We have
\begin{multline*}
\Pr { \sum_{k}a_k|\eta_k|>t } \le e^{-\frac12
{t^2}}\Ex\exp\left\{\frac12\Bigl(\sum_{k}a_k|\eta_k|\Bigr)^2\right\}
\\
\le e^{-\frac12{t^2}}\sum_k a_k
\Ex\exp\left\{\tfrac12|\eta_k|^2\right\} =
e^{-\frac12{t^2}}\Ex\exp\left\{\tfrac12|\eta|^2\right\}
\end{multline*}
where $\eta$ is a standard complex Gaussian random variable. But
\[
\Ex\exp\left\{\tfrac12|\eta|^2\right\}=\frac
1\pi\iint_{\C}e^{-\frac12
|z|^2}\,dm_2(z)=2\int_0^{+\infty}re^{-\frac 12 r^2}\,dr=2\,.
\]
\nopagebreak\hfill $\Box$

\begin{lemma}\label{lemma2.2} Let $ \{\xi_i\}_{1\le i \le n} $ be
complex Gaussian random variables, and let $\Gamma = \left(
\gamma_{ij}\right)$ be their covariance matrix; i.e.,
$\gamma_{ij}=\Ex \xi_i \bar \xi_j $. Suppose $\operatorname{\det}
\Gamma \ge 1$. Then
\[
\Pr {|\xi_i|\le \epsilon, 1\le i \le n} \le  \epsilon^{2n}\,.
\]
\end{lemma}
\noindent{\em Proof:} The joint density function of the variables
$\xi_i$ is
\[
\frac1{\pi^n\operatorname{det} \Gamma} e^{-\langle \Gamma^{-1}\xi,
\xi \rangle} \le \pi^{-n}\,.
\]
Thus
\[
\Pr { |\xi_i|\le \epsilon, 1\le i \le n } \le \frac1{\pi^n}
\idotsint\limits_{|\xi_1|<\epsilon,\, ...,\, |\xi_n|<\epsilon}
dm_2(\xi_1)\, ...\, dm_2(\xi_n) = \epsilon^{2n}\,.
\]
\nopagebreak\hfill $\Box$

\medskip
Now we want to elaborate on the well-known fact that a family
$\{\xi_i\}_{i\in I}$ of complex Gaussian random variables is
independent if and only if the covariances $\Ex\xi_i\bar\xi_j$
vanish for $i\ne j$.

\begin{lemma}\label{lemma2.3} Let $\xi_k$ be standard
complex Gaussian random variables whose covariances
$\gamma_{ij}=\Ex \xi_i\bar \xi_j$ satisfy
$$
\sum_{j\,:\,j\ne i}|\gamma_{ij}| \le \sigma  \le \frac13
\quad\text{ for all }i\,.
$$
Then $\xi_k=\zeta_k+ b_k\eta_k$ where $\zeta_k$ are independent
standard complex Gaussian random variables, $\eta_k$ are standard
complex Gaussian random variables, and $b_k\in [0,\sigma)$.
\end{lemma}
\noindent{\em Proof:} Note that $ \|M\| = {\displaystyle
\sup_i}\sum_{j}|m_{ij}|$ defines a norm on matrices $M=(m_{ij})$
(more precisely, it is the norm of $M$ as an operator in
$\ell^\infty$). Now let $ \Gamma = (\gamma_{ij}) $ be the
covariance matrix of the family $\xi_k$. Note that
$\Gamma=I-\Delta$ where $I$ is the identity matrix and
$\|\Delta\|\le \sigma$. Then $\Gamma^{-\frac12}=I-\Delta'$ with
$\|\Delta'\| \le \sigma$. Indeed, using the Taylor series
$\displaystyle (1-z)^{-\frac 12}=1+\frac 12z+\sum_{\ell\ge
2}\alpha_\ell z^\ell$ and observing that $|\alpha_\ell|<1$ for all
$\ell\ge 2$, we get
$$
\|\Delta'\|=\Bigl\|\frac 12\Delta+\sum_{\ell\ge 2}\alpha_\ell
\Delta^\ell \Bigr\|\le \frac \sigma 2+\sum_{\ell\ge 2}\sigma^\ell
=\frac\sigma 2+\frac{\sigma^2}{1-\sigma}\le \sigma\,.
$$
It remains to put $\displaystyle \zeta_k=\sum_j (\Gamma^{-\frac
12})_{kj}\xi_j$, $\displaystyle b_k\eta_k=\sum_j
\Delta'_{kj}\xi_j$. \nopagebreak\hfill $\Box$

\subsection{Operators $T_w$ and shift invariance}

The main thing we need from $f$ and $U$ is their shift invariance.
It is literally true that $U$ is shift invariant (as a random
process) but $U$ is a little bit less convenient than $f$ to work
with because, firstly, it is not a Gaussian process and, secondly,
it has singularities. The random function $f$ itself is not shift
invariant, but there is a simple transformation that makes a shift
of $f$ a G.E.F. again.

For a function $f:\C\to\C$ and a complex number $w\in \C$, define
$$
T_wf(z)=f(w+z)e^{-z\overline w}e^{-\frac 12|w|^2}\,.
$$
\begin{lemma}\label{lemma2.4} Let $f\colon \C\to\C$ be an
arbitrary function and let $w\in\C$. Let
$U(z)=\log|f(z)|-\frac12{|z|^2}$ and let $U_w(z)=\log|T_w
f(z)|-\frac12{|z|^2}$. Then $U(w+z)=U_w(z)$.
\end{lemma}
\noindent{\em Proof:}
\begin{multline*}
U_w(z)=\log|T_w f(z)|-\tfrac12{|z|^2}=\log|f(w+z)|-\Re z\overline
w-\tfrac 12|w|^2-\tfrac12{|z|^2}
\\
=\log|f(w+z)|-\tfrac12{|w+z|^2}=U(w+z)\,.
\\ \hfill\Box
\end{multline*}

\begin{lemma}\label{lemma2.5} For any $w',w''\in \C$, we
have
$$
T_{w'+w''}f=e^{i\Im w'\overline{w''}}T_{w'}T_{w''}f\,.
$$
\end{lemma}
\noindent{\em Proof:}
\begin{multline*}
(T_{w'}T_{w''}f)(z)=(T_{w''}f)(w'+z)e^{-z\overline {w'}}e^{-\frac
12|w'|^2}
\\
=f(w''+w'+z)e^{-(z+w')\overline {w''}}e^{-\frac 12|w''|^2}
e^{-z\overline {w'}}e^{-\frac 12|w'|^2}
\\
=f(w'+w''+z)e^{-i\Im w'\overline {w''}} e^{-z\overline
{(w'+w'')}}e^{-\frac 12|w'+w''|^2}
\\
=e^{-i\Im w'\overline {w''}}(T_{w'+w''}f)(z) \,.
\\ \hfill\Box
\end{multline*}

\begin{lemma}\label{lemma2.6} Let $f$ be a G.E.F.
Then $T_w f$ is also a G.E.F.
\end{lemma}
\noindent{\em Proof:} It suffices to check that the covariances of
these two complex Gaussian processes are the same. Recalling that,
for a G.E.F. $f$, we have $\Ex f(z')\overline
{f(z'')}=e^{z'\overline{z''}}$, we get
\begin{multline*}
\Ex (T_w f)(z')\overline {(T_w f)(z'')}=e^{-z'\overline w -
\overline{z''}w}e^{-|w|^2}\Ex f(w+z')\overline {f(w+z'')}
\\
= e^{-z'\overline
w-\overline{z''}w}e^{-|w|^2}e^{(w+z')\overline{(w+z'')}}=e^{z'\overline{z''}}
=\Ex f(z')\overline {f(z'')}\,.
\\ \hfill\Box
\end{multline*}

\medskip
Note that one can give another proof of this lemma using the fact
that the functions $\displaystyle \frac{z^n}{\sqrt{n!}}$ form an
orthonormal basis in the Fock-Bargmann space $\mathcal H$ (that
is, in the closure of the analytic polynomials in the weighted
space $ L^2_{\mathbb C}(\frac1{\pi}e^{-|z|^2}\, dm_2) $\,) and
that $T_w$ is a unitary operator on $\mathcal H$.

\medskip
Lemma~\ref{lemma2.6} together with Lemma~\ref{lemma2.4}
immediately imply that the random potential $U$ corresponding to a
G.E.F. $f$ is shift invariant (as a random process).

\section{Almost independence}

Let now $\displaystyle T_w f(z)=\sum_{k\ge
0}\frac{\xi_k(w)}{\sqrt{k!}}z^k$ be the (random) Taylor series of
$T_w f$ at $0$. Lemma~\ref{lemma2.6} implies that, for a fixed
$w\in\C$, $\xi_k(w)$ are independent standard Gaussian random
variables, but, of course, the covariances between $\xi_k(w')$ and
$\xi_j(w'')$ may be nontrivial for $w'\neq w''$.

\begin{lemma}\label{lemma3.1}
$$
|\Ex \xi_j(w')\overline{\xi_k(w'')}|\le
5^{\frac12(j+k)}e^{-\frac14|w'-w''|^2}\,.
$$
\end{lemma}
\noindent{\em Proof:} Let $w=w''-w'$. Since, according to
Lemma~\ref{lemma2.5}, $T_{w''}f=e^{i\Im w\overline{w'}}T_w
T_{w'}f$, the random variable $\xi_k(w'')$ equals $e^{i\Im
w\overline{w'}}\sqrt{k!}$ times the $k$-th Taylor coefficient of
the function
$\sum_\ell\xi_\ell(w')\frac{(w+z)^\ell}{\sqrt{\ell!}}e^{-z\overline
w}e^{-\frac12|w|^2}$. Hence the absolute value of the covariance
in question is just $\sqrt{k!}$ times the absolute value of the
$k$-th Taylor coefficient of the function
$\frac1{\sqrt{j!}}(w+z)^je^{-z\overline w}e^{-\frac12|w|^2}$.
According to the Cauchy inequality, this coefficient does not
exceed
$$
\frac1{\sqrt{j!}}\rho^{-k}\max_{|z|=\rho
}\left|(w+z)^je^{-z\overline w}e^{-\frac12|w|^2}\right|\le
\frac1{\sqrt{j!}}\rho^{-k}(|w|+\rho)^je^{\rho|w|}e^{-\frac12|w|^2}
$$
for any $\rho>0$. Choosing $\rho=\frac14|w|$, we get the estimate
$$
|\Ex \xi_j(w')\overline{\xi_k(w'')}|\le
\frac{\sqrt{k!}}{\sqrt{j!}}4^k\left(\frac54\right)^j|w|^{j-k}
e^{-\frac14|w|^2}\,.
$$
Exchanging the roles of $w'$ and $w''$, we get the symmetric
inequality
$$
|\Ex \xi_j(w')\overline{\xi_k(w'')}|\le
\frac{\sqrt{j!}}{\sqrt{k!}}4^j\left(\frac54\right)^k
|w|^{k-j}e^{-\frac14|w|^2}\,.
$$
Taking the geometric mean of these two estimates, we get the
statement of the lemma. \nopagebreak\hfill $\Box$
\medskip

Our next aim is to show that the G.E.F. $T_{w_j}f$ can be
simultaneously approximated by independent G.E.F. in the disk
$|z|\le r$ if all distances between the points $w_j$ are much
greater than $r$. More precisely, the following statement holds.

\begin{theorem}[almost independence]\label{thm3.1} For every
$N>0$, there exists $A=A(N)>0$ such that, for all $r>1$, and for
all families of points $w_j\in\C$ satisfying $|w_i-w_j|\ge Ar$,
$i\ne j$, we can write
$$
T_{w_j}f=f_j+h_j
$$
where $f_j$ are independent G.E.F. and $h_j$ are random analytic
functions satisfying
$$
\Pr {\max_{|z|\le r}|h_j(z)|>e^{-N r^2}} \le
2\exp\left\{-\tfrac12\exp\{Nr^2\}\right\}\,.
$$
\end{theorem}
\noindent{\em Proof:} Fix two constants $A \gg B \gg 1$ to be
chosen later. Consider the standard complex Gaussian random
variables $\xi_k(w_j)$ with $k\le B^2r^2$. We want to apply
Lemma~\ref{lemma2.3}. To this end, we need to estimate the sum of
covariances $\displaystyle \sum_{(k,j)\colon (k,j)\ne(\ell,i)}|\Ex
\xi_k(w_j)\overline{\xi_\ell(w_i)}|$. Recall that $\Ex
\xi_k(w_j)\overline{\xi_\ell(w_i)}=0$ if $i=j$. For $j\ne i$, we
can use Lemma~\ref{lemma3.1}, which yields
$$
\sum_{(k,j)\colon j\ne i}|\Ex
\xi_k(w_j)\overline{\xi_\ell(w_i)}|\le
(B^2r^2+1)5^{B^2r^2}\sum_{j\,:\,j\ne i}e^{-\frac14|w_j-w_i|^2}\,.
$$
It remains to estimate $\sum_{j\,:\,j\ne
i}e^{-\frac14|w_j-w_i|^2}$. Let $\mu$ be the counting measure of
the set $\{w_j\}_{j\ne i}$. We have
$$
\mu(D(w_i,s))\le\left\{\aligned 0,&\quad\text{ if }s<Ar\,,
\\
9A^{-2}r^{-2}s^2&\quad\text{ if }s\ge Ar\,.
\endaligned\right.
$$
(The second estimate follows from the observation that the disks
$D(w_j,\frac12Ar)$ are pairwise disjoint and contained in the disk
$D(w_i,s+\frac12Ar)\subset D(w_i,\frac32s)$ if $|w_j-w_i|\le s$
and $s\ge Ar$). Now, write
\begin{multline*}
\sum_{j\colon j\ne i}e^{-\frac14|w_j-w_i|^2} = \iint_\C
e^{-\frac14 |z-w_i|^2}\,d\mu(z)=\int_{0}^\infty\frac s2e^{-\frac14
s^2}\mu(D(w_i,s))\,ds
\\
\le 9 \int_{Ar}^\infty A^{-2}r^{-2}\frac{s^3}2e^{-\frac14 s^2}\,ds
= 9\left(1+4A^{-2}r^{-2}\right)e^{-\frac14 A^2r^2} \le 10
e^{-\frac14 A^2r^2}\,,
\end{multline*}
provided that $A \gg 1$ and $r>1$. Using this estimate, we finally
get
$$
\sum_{(k,j)\colon j\ne i}|\Ex \xi_k(w_j)\overline{\xi_\ell(w_i)}|
\le 10(B^2r^2+1)5^{B^2r^2}e^{-\frac14 A^2 r^2} \le e^{-\frac 15
A^2 r^2}\,,
$$
provided that $A \gg B \gg 1$ and $r>1$. Applying
Lemma~\ref{lemma2.3}, we conclude that
$\xi_k(w_j)=\zeta_k(w_j)+b_{kj}\eta_k(w_j)$ where $\zeta_k(w_j)$
are independent standard Gaussian random variables, $\eta_k(w_j)$
are standard Gaussian random variables, and $0\le b_{kj}\le
e^{-\frac 15 A^2 r^2}$.

For $k>B^2r^2$, let $\zeta_k(w_j)$ be independent standard complex
Gaussian random variables that are also independent with
$\zeta_{\ell}(w_i)$ for all $\ell\le B^2r^2$ and for all $i$. Put
$$
f_j(z)= \sum_{k\ge 0} \zeta_k(w_j) \frac{z^k}{\sqrt{k!}}
$$
and
$$
h_j(z) = T_{w_j}f(z)-f_j(z) = -\sum_{k>B^2r^2} \zeta_k(w_j)
\frac{z^k}{\sqrt{k!}} + \sum_{k>B^2r^2} \xi_k(w_j)
\frac{z^k}{\sqrt{k!}} + \sum_{k\le B^2r^2} b_{kj}\eta_k(w_j)
\frac{z^k}{\sqrt{k!}}\,.
$$
The G.E.F. $f_j$ are, clearly, independent and all we need to do
now is to show that $h_j$ are small in the disk $|z|\le r$. We
shall use Lemma~\ref{lemma2.1}. It reduces our task to that of
estimating the sum
$$
2\sum_{k>B^2r^2}\frac{r^k}{\sqrt{k!}} +\sum_{k \le B^2r^2} b_{kj}
\frac{r^k}{\sqrt{k!}}\le 2\sum_{k>B^2r^2}\frac{r^k}{\sqrt{k!}}
+e^{-\frac 15 A^2r^2}\sum_{k\le B^2r^2} \frac{r^k}{\sqrt{k!}}\,.
$$
Note that in the series $ \displaystyle
\sum_{k>B^2r^2}\frac{r^k}{\sqrt{k!}} $ the ratio of each term to
the previous one equals $\frac r{\sqrt {k+1}}\le \frac r{
Br}=\frac 1{ B}<\frac 12$ if $B>2$. Hence the sum does not exceed
twice the first term of the series, which is
$$
\frac{1}{\sqrt{k_0!}}r^{k_0}<\left(\frac {\sqrt er}{\sqrt
k_0}\right)^{k_0}\le \left(\frac {\sqrt er}{B r}\right)^{k_0}\le
\left(\frac {\sqrt e}{B}\right)^{B^2r^2}\le e^{-B^2r^2}\,,
$$
provided that $B>e\sqrt e$ (here $k_0$ is the smallest integer
bigger than $B^2r^2$). On the other hand, the Cauchy - Schwarz
inequality yields
$$
\sum_{k \le B^2r^2}
\frac{r^k}{\sqrt{k!}}\le\sqrt{B^2r^2+1}\sqrt{\sum_{k\ge 0}
\frac{r^{2k}}{k!}}=\sqrt{B^2r^2+1}\,e^{\frac 12r^2}\,.
$$
Thus, the sum we need to estimate does not exceed
$$
4e^{-B^2r^2} + e^{-\frac 15 A^2r^2}\sqrt{B^2r^2+1}\,e^{\frac
12r^2} \le e^{-\frac32 Nr^2}\,,
$$
provided that $A \gg B \gg \sqrt N$.

It remains to apply Lemma~\ref{lemma2.1} with $a\le e^{-\frac32
Nr^2}$, $t = e^{-Nr^2}$. \nopagebreak\hfill $\Box$

\section{Size of the potential $U$}

First, we estimate the probability that the maximum of the random
potential $U$ over the disk of radius $\rho$ is large positive.
\begin{lemma}\label{lemma4.1} For $\rho\ge 1$ and $M>0$,
\[
\Pr { \max_{|z|\le\rho} U(z) > M} \le C\rho^2 e^{-ce^{2M}}\,.
\]
\end{lemma}
\medskip\noindent{\em Proof:} Since $U$ is a stationary process,
and since the disk $\{|z|\le\rho\}$ can be covered by $C\rho^2$
copies of the unit disk, it suffices to show that
\[
\Pr { \max_{|z|\le 1} U(z) > M } \le Ce^{-ce^{2M}}\,.
\]
But this probability does not exceed $\displaystyle \Pr
{\max_{|z|\le 1}|f(z)|>e^M }$, which, in its turn, does not exceed
$\displaystyle \Pr {\sum_{k}\frac1{\sqrt{k!}}|\xi_k| > e^M }$.
Estimating the latter probability by Lemma~\ref{lemma2.1}, we get
the desired result. \nopagebreak\hfill $\Box$

\begin{lemma}\label{lemma4.2}
Suppose that $\rho\ge 1$. Then
\[
\Pr { \max_{|z|\le \rho } |f(z)|  \le e^{-3\rho^2}} <
e^{-8\rho^4}\,.
\]
\end{lemma}

\medskip\noindent{\em Proof:} Assume that $\displaystyle \max_{|z|\le \rho } |f(z)|
\le e^{-3\rho^2}$. Then by Cauchy's inequalities for the Taylor
coefficients of analytic functions, we have
\[
|\xi_n| \le \frac{\sqrt{n!}}{\rho^n}\, \max_{|z|\le\rho} |f(z)|
\le \frac{n^{n/2}}{\rho^n}  e^{-3\rho^2}\,, \qquad n=0, 1, 2, \,
...\, .
\]
The probabilities of these independent events do not exceed $ (n
\rho^{-2})^n e^{-6\rho^2}$. Thus
\[
\Pr { \max_{|z|\le \rho} |f(z)| \le e^{-3\rho^2}} \le \prod_{0 \le
n \le 2\rho^2} \left[ ( n \rho^{-2} )^n e^{-6 \rho^2} \right] \le
\left( 2^{2\rho^2} e^{-6 \rho^2}\right)^{2\rho^2} <
e^{-8\rho^4}\,.
\]
\nopagebreak\hfill$\Box$

\begin{theorem}\label{thm4.3} Given $\beta>0$, suppose that $\rho$
is sufficiently large, and that $\log^2\rho \le M \le \rho^2$.
Then the probability of the event
\[
\Big\{ {\rm there\ exists\ a\ curve\ } \gamma\subset\rho^4\D {\rm\
with\ } \operatorname{diam}(\gamma)\ge \beta \rho {\rm\ such\
that\ } \max_\gamma U < -M \Big\}
\]
does not exceed $e^{-c\rho M^{3/2}}$ with the constant $c$
depending on $\beta$.
\end{theorem}

\medskip\par\noindent Recall that by $\D$ we denote the unit disk
in the complex plane centered at the origin, $t\D$ is the disk of
radius $t$ concentric with $\D$.

\bigskip\noindent{\em Proof:} We fix a sufficiently small constant
$a < \min\left( \frac12, \frac{\beta}4 \right)$ and cover the disk
$\rho^4\D$ by the disks $D_j = D(w_j, a\sqrt{M})$, $j\in\mathcal
J$, with bounded multiplicity of covering. Clearly, $\# \mathcal J
\le CM^{-1}\rho^8 $.

Suppose that there exists a curve $\gamma\subset\rho^4\D$ with
diameter at least $\beta \rho$ and  such that $\displaystyle
\max_\gamma U < -M $. Note that if $\gamma$ enters the disk $D_j$,
then it must exit the disk $2D_j = D(w_j, 2a\sqrt{M})$; otherwise,
$4a\sqrt{M}$ (the diameter of $2D_j$) is larger than $ \beta \rho
$ (the diameter of the curve $\gamma$), which is impossible due to
our choice of $a$.

Let $A$ be the constant corresponding to the value $N=a^{-2}$ in
the almost independence theorem~\ref{thm3.1}. Having the curve
$\gamma$ and the constants $a$ and $A$, we choose a sub-collection
of well-separated disks $D_j$, $j\in \mathcal J^*$, with the
following properties:
\begin{itemize}
\item[$\bullet$] $|w_i-w_j|\ge 2Aa\sqrt{M}$ for $j\ne i$;
\item[$\bullet$] the curve $\gamma$ enters each of the disks $D_j$;
\item[$\bullet$] $\displaystyle \# \mathcal J^*= \left\lceil \frac{\beta
\rho}{2Aa\sqrt{M} + 2a\sqrt{M}} \right\rceil = \left\lceil
\frac{\beta \rho}{2(A+1)a\sqrt{M}} \right\rceil$.
\end{itemize}
By $\lceil x \rceil $ we denote the least integer $n\ge x$.
\begin{figure}[h]
\begin{center}
\setlength{\unitlength}{0.00066667in}
\begingroup\makeatletter\ifx\SetFigFont\undefined%
\gdef\SetFigFont#1#2#3#4#5{%
  \reset@font\fontsize{#1}{#2pt}%
  \fontfamily{#3}\fontseries{#4}\fontshape{#5}%
  \selectfont}%
\fi\endgroup%
{\renewcommand{\dashlinestretch}{30}
\begin{picture}(3128,2861)(0,-10)
\allinethickness{1.000pt}%
\put(1423,1423){\ellipse{2830}{2830}}
\put(1411,948){\ellipse{456}{456}}
\put(1886,2123){\ellipse{444}{444}}
\put(1186,1898){\ellipse{468}{468}}
\dashline{60.000}(1436,1411)(2811,1411)
\path(2711.000,1386.000)(2811.000,1411.000)(2711.000,1436.000)
\dashline{60.000}(1886,2111)(1986,1948)
\path(1912.398,2020.164)(1986.000,1948.000)(1955.016,2046.311)
\path(1243.449,1768.483)(1198.000,1861.000)(1194.564,1757.980)
\dashline{60.000}(1198,1861)(1386,986)
\path(1340.551,1078.517)(1386.000,986.000)(1389.436,1089.020)
\allinethickness{3.000pt}%
\path(1986,2161)(1983,2161)(1977,2161)
    (1966,2160)(1950,2159)(1929,2158)
    (1902,2157)(1872,2156)(1840,2154)
    (1807,2152)(1773,2150)(1740,2148)
    (1708,2145)(1678,2143)(1649,2141)
    (1622,2138)(1596,2135)(1571,2132)
    (1547,2129)(1523,2125)(1500,2121)
    (1476,2117)(1451,2112)(1426,2107)
    (1401,2102)(1374,2096)(1348,2089)
    (1321,2082)(1293,2074)(1266,2065)
    (1238,2056)(1211,2045)(1184,2035)
    (1158,2023)(1133,2011)(1110,1999)
    (1088,1986)(1067,1973)(1048,1959)
    (1031,1945)(1016,1930)(1003,1915)
    (991,1900)(982,1883)(974,1867)
    (967,1849)(963,1831)(959,1811)
    (958,1791)(958,1770)(960,1748)
    (964,1726)(969,1703)(975,1680)
    (983,1657)(992,1634)(1003,1611)
    (1015,1589)(1027,1568)(1041,1547)
    (1055,1527)(1070,1508)(1085,1491)
    (1101,1474)(1117,1459)(1134,1444)
    (1153,1429)(1173,1415)(1193,1402)
    (1214,1389)(1236,1377)(1259,1365)
    (1281,1353)(1304,1342)(1327,1331)
    (1350,1321)(1372,1310)(1393,1300)
    (1412,1290)(1431,1280)(1448,1270)
    (1463,1260)(1477,1250)(1488,1239)
    (1499,1228)(1507,1217)(1515,1202)
    (1520,1185)(1524,1167)(1525,1148)
    (1524,1128)(1521,1107)(1517,1085)
    (1511,1063)(1504,1042)(1496,1021)
    (1488,1001)(1479,982)(1471,965)
    (1463,949)(1455,935)(1449,923)
    (1439,907)(1432,894)(1426,885)
    (1421,879)(1417,875)(1415,873)
    (1413,872)(1411,873)
\put(2161,1873){\makebox(0,0)[lb]{{\SetFigFont{10}{12.0}{\rmdefault}{\mddefault}{\updefault}$a\sqrt
M $}}}
\put(1323,1473){\makebox(0,0)[lb]{{\SetFigFont{10}{12.0}{\rmdefault}{\mddefault}{\updefault}$\ge
2Aa\sqrt M $}}}
\put(2898,1373){\makebox(0,0)[lb]{{\SetFigFont{10}{12.0}{\rmdefault}{\mddefault}{\updefault}$\rho^4$}}}
\end{picture}
}
\end{center}
\caption{The curve $\gamma$ and the disks $D_j=D(w_j, a\sqrt{M})$}
\end{figure}
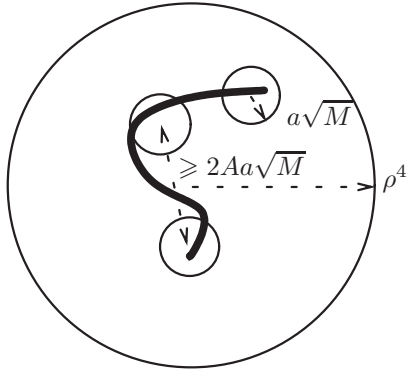

\medskip
Applying Theorem~\ref{thm3.1} with $r=2a\sqrt{M}$, we get
$T_{w_j}f = f_j + h_j$, $j\in\mathcal J^*$, where $f_j$ are
independent G.E.F. and
\[
\Pr {\max_{2a\sqrt{M}\,\D} |h_j| > e^{-4M} } \le 2 \exp \left[
-\frac12 \exp (4M) \right]\,.
\]
If $ \displaystyle \max_{j\in \mathcal J^*} \max_{2a\sqrt{M}\D}
|h_j| \le e^{-4M} $, then, for $z\in (\gamma - w_j) \cap
2a\sqrt{M}\D$, and for big enough $M$,
\[
|f_j(z)| \le e^{-M} e^{\frac12 |z|^2} + e^{-4M} \le e^{-M + 2a^2M}
+ e^{-4M} < e^{-\frac12 M}\,.
\]
Now, we introduce the independent events $(\star_j)$. We say that
the event $(\star_j)$ occurs if there exists a curve $\gamma_j$
that connects the circumferences $\left\{ |z| = a
\sqrt{M}\right\}$ and $\left\{ |z| = 2a \sqrt{M}\right\}$ such
that $|f_j(z)|< e^{-\frac12 M}$ everywhere on $\gamma_j$.

\begin{claim}\label{claim4.0} If the constant $a$ is small enough,
then $\Pr {(\star_j)} \le e^{-cM^2}$.
\end{claim}

\par\noindent{\em Proof of Claim~\ref{claim4.0}:} Consider the
function $\log |f_j|$ subharmonic in the disk $2a\sqrt{M}\D$. By
Lemma~\ref{lemma4.1}, throwing away an event of probability less
than \[ \displaystyle  Ca^2 M e^{-ce^{4a^2M}} <  e^{-cM^2}, \] we
have \[ \displaystyle \max_{z\in 2a\sqrt{M}\,\D} \left[
\log|f_j(z)| - \frac12 |z|^2 \right] \le 2a^2 M \] and hence
\[ \displaystyle \max_{2a\sqrt{M}\,\D} \log|f_j| \le 4a^2 M. \] The
curve $\gamma_j$ connects the circumferences $\{|z|=a\sqrt{M}\}$
and $\{|z|=2a\sqrt{M}\}$. Hence its harmonic measure with respect
to $\left( 2a\sqrt{M}\D \right) \setminus \gamma_j$ is bounded
from below by a positive numerical constant $c_0$ uniformly in the
disk $a\sqrt{M}\D$ (this well-known fact follows, for instance,
from \cite[Theorem~3-6]{Ahlfors}). Thus
\[
\max_{a\sqrt{M}\D} \log |f_j| \le 4a^2 M - \frac{c_0}2 M <
-\frac{c_0}4 M\,,
\]
if the constant $a$ was chosen so small that $a^2<\frac1{16} c_0$.
Then $\frac14 c_0 M > 3 (a\sqrt{M})^2$ and we can apply
Lemma~\ref{lemma4.2} to the function $f_j$ in the disk
$a\sqrt{M}\D$. The lemma yields that the probability that
($\star_j$) happens does not exceed $e^{-cM^2}$.
\nopagebreak\hfill$\Box$

\medskip
We conclude that the existence of a curve $\gamma$ satisfying the
assumptions of the theorem implies existence of a subset $\mathcal
J^*\subset \mathcal J$ with $\# \mathcal J^* = \left\lceil
\frac{\beta \rho}{2(A+1)a\sqrt{M}} \right\rceil$, such that at
least one of the following happens:
\begin{itemize}
\item[(i)] for $j\in \mathcal J^*$, the {\em independent}
events ($\star_j$) occur with $\gamma_j = (\gamma-w_j) \cap
2a\sqrt{M}\D$;
\item[(ii)]
$\displaystyle \max_{j\in\mathcal J^*} \max_{2a\sqrt{M}\,\D} |h_j|
> e^{-4M}$.
\end{itemize}
In the case (i), the probability is bounded by $\displaystyle
\binom{\# \mathcal J}{\# \mathcal J^*} \cdot \left[ e^{-cM^2}
\right]^{\# \mathcal J^*}$. Since $\displaystyle \binom{n}{k} \le
n^k$, $\#J \le C\rho^8$, and $\#J^* \le C\rho$, the first factor
is bounded by $e^{C\rho \log\rho}$. The second factor does not
exceed $e^{-cM^{3/2}\rho}$. Since $\log\rho \le M^{1/2}$, the
whole product is bounded by $e^{-cM^{3/2}\rho}$. In the case (ii),
the probability does not exceed
\[
\# \mathcal J^* \cdot \binom{\# \mathcal J}{\# \mathcal J^*} \cdot
2 e^{-ce^{cM}} < e^{C\rho \log\rho - ce^{cM}} < e^{-ce^{cM}}\,,
\]
which is much less than $e^{-c\rho M^{3/2}}$. This completes the
proof. \nopagebreak\hfill $\Box$

\section{Determinants of covariance matrices}

In this section, we estimate from below the determinant of the
covariance matrix of the complex Gaussian random variables
$\{f'(z_i) - \bar z_i f(z_i)\}_{1\le i \le n}$. This estimate will
be used in the next section when we apply Lemma~\ref{lemma2.2} to
the proof of the long gradient curve theorem~\ref{th.longcurve}.

To warm up, first, we estimate the determinant of the covariance
matrix of random variables $\{f(z_i)\}_{1\le i \le n}$, which has
a simpler structure.
\begin{lemma}[the $1$st determinant estimate]\label{lemma5.1}
Let $\{z_i\}_{1\le i \le n }\subset\C$ be a well-separated
sequence; i.e., for some $\lambda > 0$,
\[
|z_i - z_j| \ge \lambda |i-j|, \qquad 1\le i,j \le n\,,
\]
and let $\Gamma = (\gamma_{ij})$ where $\gamma_{ij} = \Ex
f(z_i)\overline{f(z_j)} = e^{z_i \bar z_j}$, $1\le i,j \le n$.
Then
\[
\operatorname{det} \Gamma \ge \left( c\lambda
\sqrt{n}\right)^{n(n-1)}\,.
\]
\end{lemma}

\medskip\par\noindent{\em Proof: } Without loss of generality,
we suppose that $n\ge 2$ (if $n=1$, the statement is obvious).
Since
\[
e^{z_i \bar z_j} = \sum_{k=0}^\infty \frac{z_i^k}{\sqrt{k!}} \cdot
\frac{\bar z_j^k}{\sqrt{k!}}\,,
\]
we have $\Gamma =  A  A^*$ with the matrix
\[
A = \left(
\begin{matrix}
1 & \frac{z_1}{\sqrt{1!}} & \frac{z_1^2}{\sqrt{2!}} & \dots &
\frac{z_1^k}{\sqrt{k!}} & \dots \\
\vdots & \vdots & \vdots & \vdots & \vdots & \vdots \\
1 & \frac{z_n}{\sqrt{1!}} & \frac{z_n^2}{\sqrt{2!}} & \dots &
\frac{z_n^k}{\sqrt{k!}} & \dots
\end{matrix}
\right)\,.
\]
Hence, by the Cauchy-Binet formula,
\[
\operatorname{det} \Gamma = \sum_t |m_t(A)|^2\,,
\]
where the sum is taken over all principal minors  $m_t( A)$ of the
matrix $ A$. We use only one principal minor
\[
m_0(A) = \left|
\begin{matrix}
1 & \frac{z_1}{\sqrt{1!}} & \dots &
\frac{z_1^{n-1}}{\sqrt{(n-1)!}} \\
\vdots & \vdots & \vdots & \vdots \\
1 & \frac{z_n}{\sqrt{1!}} & \dots &
\frac{z_n^{n-1}}{\sqrt{(n-1)!}}
\end{matrix}
\right| = \frac1{\sqrt{1!2!\,...\, (n-1)!}} \prod_{i<j}
(z_j-z_i)\,.
\]
Since the points $z_1$, ..., $z_n$ are well-separated, we get
\[
|m_0( A)| \ge \lambda^{n(n-1)/2} \sqrt{1!2!\,...\, (n-1)!}\,.
\]
Since $k! \ge k^ke^{-k}$, $k\ge 1$, we have
\begin{multline*}
1!2!\,...\, (n-1)! \ge \exp \left( \sum_{k=1}^{n-1} \left( k\log k
- k \right) \right) \\
\ge \exp \left( \int_0^n x\log x\, dx - n\log n - \frac{n(n-1)}2
\right)
\\ \ge \exp \left( \frac12 n^2\log n - n^2 - n\log n\right)\,,
\end{multline*}
and
\[
|m_0( A)|^2 \ge \left( \lambda \sqrt{n}\right)^{n(n-1)} \cdot
e^{-n^2 - n\log n} \ge \left( c\lambda \sqrt{n}\right)^{n(n-1)}\,,
\]
completing the proof of the lemma. \hfill $\Box$

\medskip In the second estimate, we fix the parameters $n\in\mathbb N$
and $r = B\sqrt{n}$ where $B\gg 1$.

\begin{lemma}[the $2$nd determinant estimate]\label{lemma5.2}
Let $\{z_i\}_{1\le i \le n}\subset \C$ be a collection of points
such that $\displaystyle |z_i - z_j| \ge \frac{r}{n} |i-j|$, and
let $\displaystyle \min_i |z_i|\ge r$. Let $\Gamma$ be the
covariance matrix of the complex Gaussian random variables
$\xi_{i} = f'(z_{i}) - \bar z_{i} f(z_{i})$, $1\le i \le n$. If
$B$ is sufficiently big, then $ \operatorname{det} \Gamma \ge 1 $.
\end{lemma}

\medskip\par\noindent The idea of the proof of this lemma is
similar to that of Lemma~\ref{lemma5.1}, though the proof is more
involved due to a more complicated structure of the covariance
matrix.

\medskip\par\noindent{\em Proof:}
First, we compute the values $\gamma_{ij} = \Ex \xi_i \bar\xi_j$:
\begin{claim}\label{claim5.3} $\gamma_{ij} = \left( 1 - |z_i - z_j|^2
\right) e^{z_i \bar z_j} $.
\end{claim}

\par\noindent{\em Proof of Claim~\ref{claim5.3} :}
\begin{multline*}
\Ex ( f'(z) - \bar z f(z) ) ( \overline{f'(w)} - w\,
\overline{f(w)} ) = \left( \partial_z - \bar z \right) \left(
\partial_{\overline w} - w \right) \Ex f(z)\overline{f(w)}
\\ = \left( \partial_z - \bar z \right)
\left( \partial_{\bar w} - w \right) e^{z \overline w} = \left( 1
+ \overline{w} z - \bar z z - w \overline w + w \bar z\right) e^{z
\overline w}
\\ = \left( 1 - |z-w|^2 \right) e^{z \overline w}\,.
\end{multline*}
\nopagebreak\hfill $\Box$

\medskip Now, we suppose that $n\ge 2$ (if $n=1$, the statement is obvious) and
factor the matrix $ \Gamma $. We have
\begin{multline*}
\gamma_{ij} = \left( \partial_{z_i} - \bar z_i \right) \left(
\partial_{\bar z_j} - z_j \right) \left( \sum_{k=0}^\infty \frac{z_i^k}{\sqrt{k!}}
\cdot \frac{\bar z_j^k}{\sqrt{k!}}  \right)
\\ = \sum_{k=0}^\infty \frac{(k-|z_i|^2)z_i^{k-1}}{\sqrt{k!}}
\cdot \frac{(k-|z_j|^2)\bar z_j^{k-1}}{\sqrt{k!}}\,.
\end{multline*}
Put
\[
A = \left(
\begin{matrix}
- \bar z_1 & \frac{1-|z_1|^2}{\sqrt{1!}} &
\frac{(2-|z_1|^2)z_1}{\sqrt{2!}} & \dots &
\frac{(k-|z_1|^2)z_1^{k-1}}{\sqrt{k!}} & \dots \\
\vdots & \vdots & \vdots & \dots & \vdots & \dots \\
- \bar z_n & \frac{1-|z_n|^2}{\sqrt{1!}} &
\frac{(2-|z_n|^2)z_n}{\sqrt{2!}} & \dots &
\frac{(k-|z_n|^2)z_n^{k-1}}{\sqrt{k!}} & \dots
\end{matrix}
\right)\,.
\]
Then $ \Gamma = A A^* $, and by the Cauchy-Binet formula,
\[
\operatorname{det} \Gamma \ge \sum_{t=1}^{n+1} |M_t(A)|^2
\]
where the sum is taken over $n+1$ principal minors $ M_t(A) $ of
the matrix $A$:
\[
M_t = \operatorname{det} \left(
\frac{(k+t-|z_i|^2)z_i^{k+t-1}}{\sqrt{(k+t)!}}\right)_{1\le i, k
\le n }\,, \qquad t=1, 2, \,...\, n+1\,.
\]
To estimate the sum of the squares of these determinants, we
introduce the determinants of  simpler structure:
\[
\mu_t = \operatorname{det} \left(
\frac{(k+t-|z_i|^2)z_i^{k-1}}{\sqrt{k!}}\right)_{1\le i, k \le n
}\,, \qquad t=1, 2, \,...\, n+1\,.
\]

\begin{claim}\label{claim5.4}
For $ 1\le t \le n+1 $, $|M_t| \ge |\mu_t|$.
\end{claim}
\noindent{\em Proof of Claim~\ref{claim5.4}:} follows by a
straightforward estimate of the ratio
\[
\left| \frac{M_t}{\mu_t} \right|  = \left( |z_1| \cdot |z_2| \cdot
\, ...\, \cdot |z_n| \right)^t \cdot \sqrt{\frac{1!}{(t+1)!} \cdot
\frac{2!}{(t+2)!} \cdot \, ...\, \cdot \frac{n!}{(t+n)!}}\,.
\]
Since $|z_i|\ge r$, we have $ \left( |z_1| \cdot |z_2| \cdot \,
...\, \cdot |z_n| \right)^t \ge r^{tn}$. For each integer $t$
between $1$ and $n+1$,
\begin{multline*}
\sqrt{\frac{1!}{(t+1)!} \cdot \frac{2!}{(t+2)!} \cdot \, ...\,
\cdot \frac{n!}{(t+n)!}} \ge \left( \frac{n!}{(t+n)!}
\right)^{n/2} \\
= \left( \frac1{(n+1)\,...\, (n+t)}  \right)^{n/2} \ge \left(
\frac1{n+t} \right)^{nt/2} \ge \left( \frac1{3n} \right)^{nt/2}\,.
\end{multline*}
Thus
\[
\left| \frac{M_t}{\mu_t} \right| \ge \left( \frac{r^2}{3n}
\right)^{nt/2} = \left( \frac13 B^2 \right)^{nt/2} \ge 1\,,
\]
provided that $\displaystyle B \ge \sqrt{3} $. This proves the
claim. \nopagebreak\hfill $\Box$

\medskip Now, we complete the proof of
Lemma~\ref{lemma5.2}. Observe that $\mu_t$ is a polynomial of
degree $ n $ in $ t $. We use a version of the pigeonhole
principle:

\begin{claim}\label{claim5.5}
Let $P(t)$ be a polynomial of degree $n$ with the leading
coefficient $a$. Then $\displaystyle \max \big\{ |P(t)|\colon
t\in\{1, 2, \,...\,, n+1\} \big\}\ge |a|2^{-n}$.
\end{claim}

\par\noindent{\em Proof of Claim~\ref{claim5.5}: } We have
\[
\max_{t\in \{1, 2, \,...\,, n+1\} } |P(t)| = |a| \cdot \max_{t\in
\{1, 2, \,...\,, n+1\}} |t - \tau_1|\, ...\, |t-\tau_n|\,,
\]
where $\tau_1$, ..., $\tau_n$ are the zeroes of $P(t)$. We have
$n+1$ disjoint $\frac12$-neighbourhoods of the points $1$, $2$,
..., $n+1$ in $\C$. At least one of them is free of the zeroes of
$P(t)$. Hence at the center of this neighbourhood, the absolute
value of $P$ cannot be smaller than $|a|2^{-n}$; whence the claim.
\nopagebreak\hfill $\Box$

\medskip We apply this claim to the polynomial $\mu_t$. Its leading
coefficient equals
\[
a_n = \lim_{t\to\infty} t^{-n}\mu_t = \det \left(
\frac{z_i^{k-1}}{\sqrt{k!}} \right)_{1\le i,k \le n}\,.
\]
We've already estimated this determinant in the proof of the model
lemma~\ref{lemma5.1}. We get $|a_n| \ge \left( c\lambda \sqrt{n}
\right)^{n(n-1)/2}$ with $ \lambda = \frac{r}{n}$ and $ n =
(r/B)^2 $. Then $|a_n| \ge (cB)^{n(n-1)/2}$, and
\[
\max_{t\in \{1, 2, \,...\,, n+1\} } |\mu_t| \ge |a_n|2^{-n} \ge
\left( \frac12 cB \right)^{n(n-1)/2}\,,
\]
whence
\[
\operatorname{det} \Gamma \ge \sum_{t=1}^{n+1} |M_t|^2 \ge
\sum_{t=1}^{n+1} |\mu_t|^2 \ge \left( \frac12 cB \right)^{n(n-1)}
\ge 1\,,
\]
provided that the constant $ B $ is chosen sufficiently large.
This completes the proof of Lemma~\ref{lemma5.2}.
\nopagebreak\hfill $\Box$

\section{The long gradient curve theorem}

Till the end of the proof, we fix $\delta=\frac1{20}$. Everywhere
below we shall assume that $ R\gg 1 $. In the proof, we work with
three scales: starting with the macroscopic $R$-scale, we move to
the intermediate $\sqrt{\log R}$-scale, and then to the
microscopic $\frac1R$-scale.

\subsection{Bad squares}

Suppose that there exists a gradient curve $\Gamma $ connecting
$\partial Q(0,R)$ and $\partial Q(0,2R)$. Due to
Lemma~\ref{lemma4.1}, we can assume that $U\le R^\delta $
everywhere on  $Q(0,2R)$, hence on $\Gamma$: the probability of
the opposite event does not exceed $CR^2 e^{-ce^{2R^\delta}}$.
Suppose that there is a point on $\Gamma$ where $U=-R^\delta$.
Since $\Gamma$ is a gradient curve, if such a point exists, then
it is unique. This point splits $\Gamma$ into two parts:
$\Gamma_1$ where $U<-R^\delta$, and $\Gamma_2$ where $U\ge
-R^\delta$. If $U\ne -R^\delta$ on $\Gamma$, then one of these
parts is empty. One of the curves $\Gamma_1$, $\Gamma_2$ must
connect either $\partial Q(0, R)$ with $\partial Q(0, \sqrt{2}
R)$, or $\partial Q(0, \sqrt{2} R)$ with $\partial Q(0, 2R)$. If
this is the curve $\Gamma_1$, then its diameter is larger than
$cR$. By Theorem~\ref{thm4.3}, the probability of this event does
not exceed $e^{-cR^{1+\frac32 \delta}}$, and we are done.

Thus the proof boils down to the case when the gradient curve
$\Gamma$ connects $\partial Q(0, R)$ with $\partial Q(0, \sqrt{2}
R)$ and $-R^\delta \le U \le R^\delta $ everywhere on $\Gamma$. In
this case,
\begin{equation}\label{eq6.1}
\int_\Gamma |\nabla U(z)|\, |dz| \le 2R^\delta\,;
\end{equation}
that is, the gradient  $\nabla U$ is small in the mean on
$\Gamma$. We will not use anymore that $\Gamma$ is a gradient
curve; starting this moment, it is an arbitrary curve connecting
$\partial Q(0, R)$ with $\partial Q(0, \sqrt{2} R)$ such that
\eqref{eq6.1} happens.

\medskip
We take $\displaystyle  r = \frac14 \sqrt{\delta\log R} $ and fix
the standard partition of the complex plane $\C$ into squares
$Q(w_j,r)$ with side length $2r$. Let $\mathcal J$ be the set of
indices $j$ for which the square $Q(w_j,2r)$ is entirely contained
in the ``square annulus'' $Q(0, \sqrt{2}R)\setminus Q(0, R)$. Note
that $\#\mathcal J\le (R/r)^2$.

\begin{definition}[bad squares] Let $j\in\mathcal J$. We shall call the
standard square $Q(w_j,r)$ {\rm bad} if there exists a curve $
\gamma_j $ joining $\partial Q(w_j, r)$ with $\partial Q(w_j, 2r)$
such that
\begin{equation}\label{eq.bad1}
\int_{\gamma_j}|\nabla U(z)|\,|dz|< r R^{2\delta-1}\,.
\end{equation}
We shall call the square $Q(w_j,r)$ {\rm good} if it is not bad.
\end{definition}

By $\mathcal F\subset \mathcal J$ we denote the family of all
indices $j$ such that the square $Q(w_j, r)$ intersects the curve
$\Gamma$.

\begin{lemma}\label{lemma.14}
At most $ 8 R^{1-\delta}/r$ of the squares $\big\{ Q(w_j, r)
\big\}_{j\in\mathcal F}$ are good.
\end{lemma}
\noindent{\em Proof: } Let $N$ denote the number of good squares
$Q(w_j,r)$. By $\gamma_j$ we denote a connected part of
$\Gamma\cap (Q(w_j,2r)\setminus Q(w_j,r))$ that joins $\Q(w_j,r)$
with $\Q(w_j,2r)$. Since almost every point of the curve $\Gamma$
belongs to at most $4$ squares $Q(w_j,2r)$, we can write
$$
NrR^{2\delta-1}\le \sum_{j\colon Q(w_j,r)\text{ is
good}}\int_{\gamma_j}|\nabla U(z)|\,|dz| \le 4\int_\Gamma |\nabla
U(z)|\,|dz|\le 8 R^{\delta}\,,
$$
whence the estimate. \nopagebreak\hfill $\Box$

\medskip The immediate consequence of Lemma~\ref{lemma.14} is that
the existence of a curve $\Gamma $ connecting $\partial Q(0,R)$
with $\partial Q(0,\sqrt{2}R)$ such that \eqref{eq6.1} happens
implies the existence of a family $\mathcal F$ of squares
$Q(w_j,r)$ of cardinality $\displaystyle L\ge \frac{cR}{r}$ and a
subfamily $\mathcal F'\subset\mathcal F$ of bad squares of
cardinality at least $\displaystyle L - \frac{8R^{1-\delta}}{r}
\ge \frac12 L$ in that family.
\medskip
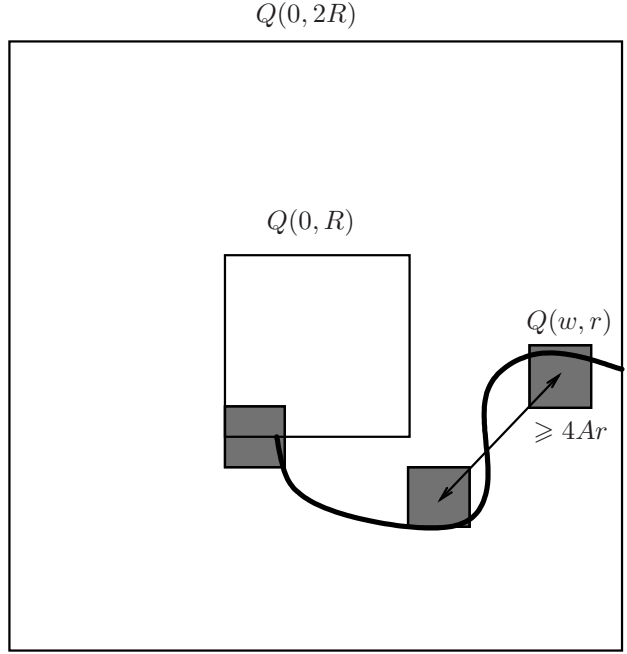
\begin{figure}[h]
\begin{center}
\setlength{\unitlength}{0.00066667in}
\begingroup\makeatletter\ifx\SetFigFont\undefined%
\gdef\SetFigFont#1#2#3#4#5{%
  \reset@font\fontsize{#1}{#2pt}%
  \fontfamily{#3}\fontseries{#4}\fontshape{#5}%
  \selectfont}%
\fi\endgroup%
{\renewcommand{\dashlinestretch}{30}
\begin{picture}(4759,5025)(0,-10)
\put(3987,1625){\makebox(0,0)[lb]{{\SetFigFont{10}{12.0}{\rmdefault}{\mddefault}{\updefault}
$\ge 4Ar$}}}
\put(1912,4850){\makebox(0,0)[lb]{{\SetFigFont{10}{12.0}{\rmdefault}{\mddefault}{\updefault}$Q(0,
2R)$}}}
\put(3999,2462){\makebox(0,0)[lb]{{\SetFigFont{10}{12.0}{\rmdefault}{\mddefault}{\updefault}$Q(w,
r)$}}}
\put(1999,3237){\makebox(0,0)[lb]{{\SetFigFont{10}{12.0}{\rmdefault}{\mddefault}{\updefault}$Q(0,
R)$}}}
\allinethickness{1.000pt}%
\texture{aadddddd ddaaaaaa aa555555 55aaaaaa aad5d5d5 d5aaaaaa aa555555 55aaaaaa
    aadddddd ddaaaaaa aa555555 55aaaaaa aaddd5dd d5aaaaaa aa555555 55aaaaaa
    aadddddd ddaaaaaa aa555555 55aaaaaa aad5d5d5 d5aaaaaa aa555555 55aaaaaa
    aadddddd ddaaaaaa aa555555 55aaaaaa aaddd5dd d5aaaaaa aa555555 55aaaaaa }
\shade\path(4024,2387)(4499,2387)(4499,1900)
    (4024,1900)(4024,2387)
\path(4024,2387)(4499,2387)(4499,1900)
    (4024,1900)(4024,2387)
\shade\path(3087,1437)(3562,1437)(3562,975)
    (3087,975)(3087,1437)
\path(3087,1437)(3562,1437)(3562,975)
    (3087,975)(3087,1437)
\shade\path(1674,1912)(2137,1912)(2137,1437)
    (1674,1437)(1674,1912)
\path(1674,1912)(2137,1912)(2137,1437)
    (1674,1437)(1674,1912)
\path(1674,3087)(3099,3087)(3099,1675)
    (1674,1675)(1674,3087)
\path(3387.610,1276.798)(3337.000,1187.000)(3423.914,1242.417)
\path(3337,1187)(4249,2150)
\path(4198.390,2060.202)(4249.000,2150.000)(4162.086,2094.583)
\path(12,4750)(4737,4750)(4737,12)
    (12,12)(12,4750)
\allinethickness{2.000pt}%
\path(2074,1675)(2074,1674)(2075,1671)
    (2076,1662)(2079,1647)(2082,1625)
    (2087,1597)(2093,1566)(2100,1533)
    (2107,1500)(2115,1469)(2123,1439)
    (2131,1412)(2139,1388)(2148,1366)
    (2157,1346)(2167,1327)(2178,1310)
    (2191,1293)(2201,1280)(2213,1267)
    (2226,1255)(2240,1242)(2255,1229)
    (2271,1217)(2289,1204)(2308,1192)
    (2328,1179)(2350,1167)(2372,1155)
    (2396,1144)(2421,1132)(2447,1121)
    (2473,1111)(2501,1101)(2528,1091)
    (2557,1082)(2586,1073)(2615,1064)
    (2645,1056)(2676,1047)(2701,1041)
    (2727,1035)(2754,1029)(2782,1023)
    (2810,1017)(2840,1011)(2870,1005)
    (2901,1000)(2933,995)(2965,990)
    (2998,985)(3031,981)(3064,977)
    (3098,974)(3131,971)(3163,969)
    (3195,967)(3226,966)(3257,965)
    (3286,965)(3315,966)(3342,968)
    (3368,970)(3392,973)(3416,976)
    (3438,981)(3459,986)(3479,991)
    (3499,999)(3518,1007)(3536,1016)
    (3553,1026)(3569,1037)(3584,1049)
    (3598,1062)(3612,1077)(3624,1092)
    (3636,1109)(3646,1126)(3656,1144)
    (3665,1163)(3672,1183)(3679,1203)
    (3685,1224)(3690,1245)(3694,1266)
    (3698,1287)(3700,1309)(3702,1330)
    (3703,1351)(3704,1373)(3704,1394)
    (3704,1416)(3703,1437)(3702,1461)
    (3701,1485)(3699,1509)(3697,1534)
    (3694,1560)(3692,1586)(3689,1612)
    (3687,1640)(3684,1667)(3682,1695)
    (3679,1723)(3678,1751)(3676,1778)
    (3675,1806)(3674,1833)(3674,1859)
    (3675,1884)(3676,1909)(3679,1932)
    (3681,1955)(3685,1977)(3690,1998)
    (3695,2018)(3702,2037)(3709,2056)
    (3718,2075)(3727,2092)(3738,2110)
    (3750,2127)(3764,2144)(3778,2161)
    (3794,2177)(3811,2192)(3829,2207)
    (3848,2222)(3867,2235)(3888,2248)
    (3909,2260)(3930,2271)(3952,2281)
    (3974,2289)(3996,2297)(4018,2304)
    (4040,2310)(4061,2315)(4083,2319)
    (4105,2322)(4127,2325)(4150,2326)
    (4174,2327)(4199,2327)(4224,2326)
    (4249,2324)(4275,2321)(4301,2318)
    (4328,2314)(4354,2309)(4380,2304)
    (4406,2299)(4431,2293)(4456,2287)
    (4480,2281)(4502,2274)(4524,2268)
    (4544,2262)(4562,2256)(4580,2250)
    (4596,2245)(4611,2240)(4624,2235)
    (4650,2227)(4671,2220)(4688,2214)
    (4702,2210)(4714,2206)(4724,2203)
    (4731,2201)(4735,2200)(4737,2200)
\end{picture}
}
\end{center}
\caption{The gradient curve $\Gamma$ generates separated bad
squares}
\end{figure}

\medskip
Let $A$ be the constant corresponding to $N = 4\delta^{-1}  $ in
the almost independence theorem~\ref{thm3.1}. Let $\mathcal
J'\subset\mathcal J$ satisfy $|w_i-w_j|\ge 4Ar $ for
$i,j\in\mathcal J'$, $i\ne j$. According to Theorem~\ref{thm3.1}
applied to $4r$ instead of $r$, we can represent $T_{w_j} f$ as
$f_j + h_j$ where $f_j$ are independent G.E.F. and all the
functions $h_ j$ are small in the disk $|z|\le 4r$. We set
\[
\Omega_*=\{\max_{j\in\mathcal J'}\max_{|z|\le
4r}|h_j(z)|>R^{-4}\}\,.
\]
Then for any $j\in\mathcal J'$, $\Pr {\displaystyle \max_{|z|\le
4r}|h_j(z)|>R^{-4} } \le 2e^{-\frac12 R^4} $ (recall that $e^{-16
Nr^2}=R^{-4}$ for our choice of $r$ and $N$). Therefore, $\Pr
{\Omega_*} \le 2R^2 e^{-\frac12 R^4} < e^{-c R^4}$.

\medskip
The next proposition is the central part in the proof of the long
gradient curve theorem:
\begin{proposition}\label{subthm6.4}
There exist events $\Omega_j$ with $\Pr {\Omega_j}\le e^{-cr^4}$
depending only on $f_j$ (and, thereby, independent), and such
that, for any $j\in\mathcal J'$,
\[
\{\, Q(w_j, r) {\rm \ is\ bad\ } \} \subset \Omega_j\cup
\Omega_*\,.
\]
\end{proposition}

\medskip Now, using this proposition, we complete the proof of
the long gradient curve theorem.  We choose a family $\mathcal
F''\subset\mathcal F'$ of $cA^{-2}r^{-1}R$ ``$4Ar$-separated
squares'' (that is, all the distances between the centers of these
squares are not less than $4Ar$), and discard the rest of
$\mathcal F'$. From Proposition~\ref{subthm6.4} we see that the
probability that a {\em given} subfamily $\mathcal F''\subset
\mathcal J$ of $cA^{-2}r^{-1}R$ squares is bad does not exceed
\[
(e^{-cr^4})^{cA^{-2}r^{-1}R} + e^{-cR^4} = e^{-cR(\log R)^{3/2}} +
e^{-cR^4} \le 2e^{-cR(\log R)^{3/2}}\,,
\]
provided that $R\gg 1$. At last, we have at most
\[
\left( \# \mathcal J\right)^{\# \mathcal F''} \le \left( CR^2
\right)^{CR} \le e^{CR\log R }
\]
ways to choose $\mathcal F''$ in $\mathcal J$. This does not harm
the previous upper bound. Hence the long gradient curve theorem is
proved  (modulo the proposition). \nopagebreak\hfill $\Box$

\subsection{Proof of the proposition}

Assume that the event $\Omega_*$ does not occur. Then $ T_{w_j} f
= f_j + h_j $ where $\displaystyle \max_{4r\D} |h_j| \le R^{-4}$.
We fix $j\in\mathcal J'$ and aim at building an event $\Omega_j$
depending only on $f_j$ of probability $ \Pr {\Omega_j} \le
e^{-cr^4}  $ and such that, if the square $ Q(w_j, r) $ is bad and
$h_j$ is small as above, then $ \Omega_j $ must occur. To simplify
the notation, we set $w=w_j$.

\medskip
Fix the partition of the complex plane $\C$ into standard squares
with side length $\frac 2R$.
\begin{definition}[black squares]
We shall call a standard square $Q(\zeta,\frac 1R) \subset Q(w, r)
$ {\rm black} if $ \displaystyle \inf_{Q(\zeta, \frac1{R})}
|\nabla U| \le R^{3\delta-1} $. Otherwise, the  square $Q(\zeta,
\frac1{R})$ is called {\rm white}.
\end{definition}

First, we check that if the square $Q(\zeta, \frac1{R}) $ is black
(i.e., the gradient $\nabla U$ is small somewhere in this square),
and the functions $f_j$ and $h_j$ are not too large, then the
function $f'_j (z) - \bar z f_j(z)$ must be small at the center
$\zeta - w$ of the shifted square.
\begin{lemma}\label{sublemma6.12} Suppose that

\smallskip\par\noindent {\rm (i)} the square $Q(\zeta, \frac1R)$ is black;

\smallskip\par\noindent {\rm (ii)} $\displaystyle \max_{4r\D} |h_j| \le R^{-4}$;

\smallskip\par\noindent {\rm (iii)} $\displaystyle \max_{4r\D} |f_j| \le R^\delta $.

\smallskip\par\noindent Then
\begin{equation}\label{eq6.c}
|f_j'(\zeta - w) - ( \overline{\zeta-w} ) f_j(\zeta-w)| <
R^{6\delta -1}\,.
\end{equation}
\end{lemma}

\par\noindent{\em Proof:} We have $U(w + z) = \log |T_w f(z)| - \frac12
|z|^2$, whence $\displaystyle |\nabla U(w+z)| = \left| \frac{(T_w
f)'(z)}{(T_w f) (z)} - \bar z\right|$. Thereby,
\begin{equation}\label{eq6.d}
|f_j' (z) - \bar z f_j (z)| \le |h_j' (z) - \bar z h_j (z)|  +
|\nabla U(w + z)| \left( |f_j(z)| + |h_j(z)|\right)\,.
\end{equation}
Since the square $Q(\zeta, \frac1R)$ is black, there exists a
point $z$ such that $w+z \in Q(\zeta, \frac1{R})$ and $|\nabla U
(w+z)| \le R^{3\delta-1} $. The other terms on the RHS of
\eqref{eq6.d} are readily estimated using assumptions (ii) and
(iii) and Cauchy's inequality for the derivative of an analytic
function. We get
\[
|f_j' (z) - \bar z f_j (z)| \le 2rR^{-4} + R^{3\delta-1} \left(
R^\delta + R^{-4}\right) < R^{5\delta -1}\,.
\]
It remains to replace $z$ by $\zeta-w$ on the LHS.

By Cauchy's inequalities,
\[
\max_{Q(0, 2r)} |f_j'| \le r^{-1} R^\delta\,, \quad \text{\rm and
}\ \max_{Q(0, 2r)} |f_j''| \le 2r^{-2} R^\delta\,.
\]
Hence the operator norm of the differential of $f_j'(z) - \bar z
f_j(z)$ does not exceed \[ 2r^{-2} r R^\delta + 2r r^{-1} R^\delta
+ 2 R^\delta < 5R^\delta \] everywhere in $Q(0, 2r)$. Since
$|z-(\zeta-w)|\le \frac{\sqrt{2}}{R}$ and $R \gg 1$, we are done.
\hfill $\Box$

\medskip

\begin{figure}[h]
\begin{center}
\setlength{\unitlength}{0.00065000in}
\begingroup\makeatletter\ifx\SetFigFont\undefined%
\gdef\SetFigFont#1#2#3#4#5{%
  \reset@font\fontsize{#1}{#2pt}%
  \fontfamily{#3}\fontseries{#4}\fontshape{#5}%
  \selectfont}%
\fi\endgroup%
{\renewcommand{\dashlinestretch}{30}
\begin{picture}(3837,3839)(0,-10)
\allinethickness{1.000pt}%
\dashline{60.000}(725,3087)(3087,3087)(3087,737)
    (725,737)(725,3087)
\path(1212,2624)(2587,2624)(2587,1212)
    (1212,1212)(1212,2624)
\texture{bfffffff ffffffff ffffffff ffbbbbbb bbffffff ffffffff
ffffffff fffbfbfb
    fbffffff ffffffff ffffffff ffbbbbbb bbffffff ffffffff ffffffff fffbfbfb
    fbffffff ffffffff ffffffff ffbbbbbb bbffffff ffffffff ffffffff fffbfbfb
    fbffffff ffffffff ffffffff ffbbbbbb bbffffff ffffffff ffffffff fffbfbfb }
\shade\path(2587,2224)(2687,2224)(2687,2099)
    (2587,2099)(2587,2224)
\path(2587,2224)(2687,2224)(2687,2099)
    (2587,2099)(2587,2224)
\dashline{60.000}(975,2837)(2850,2837)(2850,962)
    (975,962)(975,2837)
\shade\path(2775,2337)(2900,2337)(2900,2212)
    (2775,2212)(2775,2337)
\path(2775,2337)(2900,2337)(2900,2212)
    (2775,2212)(2775,2337)
\shade\path(3050,2162)(3162,2162)(3162,2062)
    (3050,2062)(3050,2162)
\path(3050,2162)(3162,2162)(3162,2062)
    (3050,2062)(3050,2162)
\dashline{60.000}(487,3324)(3337,3324)(3337,487)
    (487,487)(487,3324)
\shade\path(3300,1712)(3412,1712)(3412,1599)
    (3300,1599)(3300,1712)
\path(3300,1712)(3412,1712)(3412,1599)
    (3300,1599)(3300,1712)
\dashline{60.000}(250,3562)(3575,3562)(3575,249)
    (250,249)(250,3562)
\shade\path(3525,1499)(3637,1499)(3637,1387)
    (3525,1387)(3525,1499)
\path(3525,1499)(3637,1499)(3637,1387)
    (3525,1387)(3525,1499)
\path(12,3812)(3800,3812)(3800,12)
    (12,12)(12,3812)
\allinethickness{3.000pt}%
\path(2587,2162)(2589,2164)(2593,2168)
    (2601,2174)(2612,2184)(2626,2196)
    (2642,2209)(2660,2223)(2680,2236)
    (2701,2250)(2724,2262)(2747,2272)
    (2773,2280)(2801,2286)(2831,2289)
    (2862,2287)(2889,2282)(2913,2275)
    (2933,2268)(2949,2261)(2960,2255)
    (2968,2251)(2973,2247)(2976,2244)
    (2978,2241)(2980,2237)(2982,2233)
    (2987,2226)(2994,2216)(3005,2203)
    (3020,2185)(3039,2161)(3062,2133)
    (3087,2099)(3108,2069)(3128,2039)
    (3145,2010)(3161,1983)(3174,1957)
    (3184,1935)(3193,1914)(3200,1895)
    (3205,1878)(3209,1861)(3213,1846)
    (3217,1830)(3221,1815)(3227,1799)
    (3233,1782)(3241,1764)(3250,1745)
    (3262,1724)(3275,1704)(3291,1684)
    (3308,1665)(3325,1649)(3350,1634)
    (3371,1630)(3386,1636)(3395,1649)
    (3400,1668)(3403,1689)(3404,1713)
    (3405,1736)(3407,1758)(3411,1778)
    (3418,1794)(3429,1804)(3444,1806)
    (3462,1799)(3476,1787)(3488,1772)
    (3498,1755)(3505,1736)(3511,1717)
    (3514,1699)(3517,1681)(3518,1663)
    (3519,1645)(3519,1627)(3520,1609)
    (3523,1591)(3526,1572)(3532,1553)
    (3539,1534)(3549,1516)(3561,1500)
    (3575,1487)(3595,1477)(3613,1478)
    (3626,1485)(3636,1498)(3643,1513)
    (3648,1530)(3653,1547)(3659,1564)
    (3668,1580)(3680,1594)(3695,1605)
    (3712,1612)(3727,1614)(3740,1611)
    (3749,1604)(3756,1596)(3762,1586)
    (3766,1574)(3770,1563)(3772,1553)
    (3774,1545)(3775,1540)(3775,1537)
\end{picture}
}
\end{center}
\caption{The curve $\gamma$ and a sequence of black squares it
generates}
\end{figure}
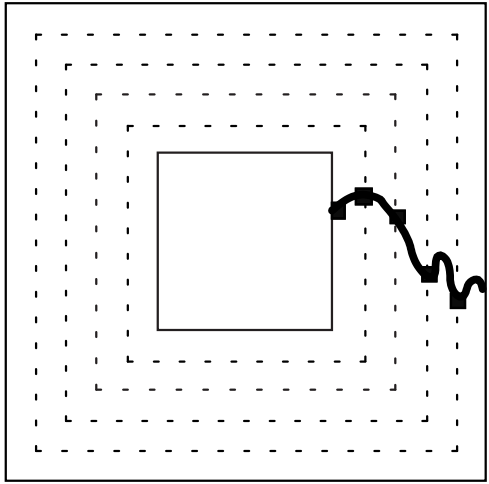

\medskip
Assume that the square $Q(w, r)$ is bad; i.e., there exists
a curve $\gamma$ joining $\partial Q(w, r)$ with $\partial Q(w,
2r)$ such that
\[
\int_{\gamma} |\nabla U (z)|\, |dz| < rR^{2\delta-1}\,.
\]
We fix an integer $ n= (r/B)^2 $ with $ B\gg 1 $. For any $ t\in
[0, 1] $, we put
\[
r_i(t) = r + \frac{(i-1)+ t}{n}\,r\,, \qquad 1\le i\le n\,.
\]
For each $t$, the squares $\partial Q(w, r_i(t))$ form a ``chain
of $n$ fences'', and the curve $\gamma$ crosses this chain at
least $n$ times. It may happen that, for some value $t$, the
gradient $\nabla U$ is not small at most of the crossing points,
or even at all of them. However, as we shall see, for a large
subset of $t\in [0, 1]$, the gradient $\nabla U$ is sufficiently
small at  $n$ crossing points $\gamma\cap
\partial Q(w, r_i(t))$, $1\le i \le n$, to guarantee that the
corresponding $\frac1{R}$-squares containing these points are
black. For each $t\in [0, 1]$, we denote by $ B(t)\subset \{1, 2,
\, ...\, n\}$ the subset of those $i$'s that at least one point
from the set $\gamma \cap\partial Q(w, r_i(t) )$ is covered by a
black square. By $m_1$ we denote the one-dimensional Lebesgue
measure.

\begin{lemma}\label{claim6.11}
Suppose that the square $Q(w, r)$ is bad. Then
\begin{equation}\label{eq6.b}
m_1 \big\{t\in [0, 1]\colon \#B(t)=n \big\} \ge \frac12\,.
\end{equation}
\end{lemma}

\par\noindent{\em Proof: } Let $L$ be the measure of the set of $\rho\in
[r, 2r]$ such that the intersection $\gamma \cap \partial Q(w,
\rho)$ is contained in white squares. Then
\[
\int_{\gamma} |\nabla U(z)|\, |dz| \ge R^{3\delta-1}L\,.
\]
Since the square $Q(w, r)$ is bad, the LHS does not exceed
$rR^{2\delta-1}$, and we see that $L \le rR^{-\delta}$. On the
other hand,
\[
L = \frac{r}{n} \int_0^1 (n - \# B(t))\, dt \ge \frac{r}{n} \cdot
m_1 \big\{ t\in [0, 1]\colon \#B(t) \le n-1\big\}\,,
\]
whence $\displaystyle m_1 \big\{  t\in [0, 1]\colon \#B(t) \le n-1
\big\} \le \frac{n}{R^\delta} \le \frac{r^2}{B^2 R^\delta} \ll 1 $
if $R\gg 1$. Hence the lemma. \nopagebreak\hfill $\Box$

\medskip For each $t\in [0, 1]$, consider  the collection $\mathfrak Z(t)$ of
``configurations'' $\mathfrak z = \{z_1, \, ...\, z_n \}$ of $n$
points such that each point $z_i$ is a center of a standard square
$Q\left(z_i, \frac 1{R} \right)$ from our partition that has a
non-void intersection with $\partial Q(w, r_i(t))$. Let us
introduce the events
\[
\Upsilon_j (\mathfrak z) = \left\{ \max_{1\le i \le n }
|f_j'(z_i-w)-(\overline{z_i-w}) f_j(z_i-w)|\le R^{6\delta -1}
\right\} \ {\rm and} \  \Omega_j(t) = \bigcup_{\mathfrak z \in
\mathfrak Z(t)} \Upsilon_j (\mathfrak z),
\]
and estimate their probabilities. Our estimate is based on the
lower bound for the determinant of the covariance matrix of
complex Gaussian random variables $ \left\{ f_j'(z_i) - \bar z_i
f_j(z_i) \right\}$ given in Lemma~\ref{lemma5.2}.

\begin{lemma}\label{lemma6.a} Given $t\in[0, 1]$,
$\Pr { \Omega_j(t)}\le e^{-cr^4} $.
\end{lemma}

\par\noindent{\em Proof:} First, we estimate the probability of
the event $\Upsilon_j (\mathfrak z)$.
\begin{claim}\label{claim6.18} For any configuration  $\mathfrak z
\in \mathfrak Z(t)$,  $\Pr {\Upsilon_j (\mathfrak z) } \le
R^{-2(1-6\delta)n}$.
\end{claim}
\noindent{\em Proof:} This is a straightforward combination of
Lemmas~\ref{lemma2.2} and~\ref{lemma5.2}. \hfill$\Box$

\medskip Next, we estimate the cardinality of the collection $\mathfrak Z(t)$
(recall that $t\in [0, 1]$ is fixed).

\begin{claim}\label{claim6.19} $\# \mathfrak Z(t) \le C^n (Rr)^n$.
\end{claim}
\noindent{\em Proof:} For each $i$, there are at most $CRr$
standard $\frac1{R}$-squares that intersect $\partial Q(w, r_i(t))
$. Therefore, there are at most $CRr$ choices for the centers
$z_i$ of these squares and the number of the corresponding
configurations $\mathfrak z$ cannot exceed $(CrR)^n$. Hence the
claim. \nopagebreak\hfill $\Box$

\medskip Using Claims~\ref{claim6.19}
and~\ref{claim6.18}, we get
\[
\Pr {\Omega_j(t)} \le C^n (Rr)^n \cdot R^{-2(1-6\delta)n} \le C^n
R^{-(1-13\delta)n}\,,
\]
provided that $R$ is sufficiently big. Recalling that
$\displaystyle 13\delta = \frac{13}{20} < 1$, $\displaystyle n =
\frac{r^2}{B^2}$, and $\displaystyle  r^2 = \frac1{16} \delta\log
R$, we see that $ \Pr {\Omega_j(t)} \le e^{-cr^4}$. This proves
Lemma~\ref{lemma6.a}. \nopagebreak\hfill$\Box$

\medskip Define the events
\[
\Omega_j' = \Big\{\omega\in\Omega\colon  m_1 \{t\in [0, 1]\colon
\omega\in\Omega_j(t) \} \ge \frac12 \Big\}\,,
\]
and
\[
\Omega_j = \Omega_j' \cup \Big\{ \max_{4r\D} |f_j| \ge R^\delta
\Big\} \,.
\]
Note that the event $\Omega_j$ depends only on $f_j$. By
Lemmas~\ref{sublemma6.12} and~\ref{claim6.11},
\[
\big\{ Q(w, r) \ {\rm is\ bad\ }  \big\}\subset \Omega_j \cup
\big\{ \max_{4r\D} |h_j| \ge R^{-4} \big\} \,.
\]

\begin{lemma} The probability of the event
$\displaystyle \big\{ \max_{4r\D} |f_j| \ge R^\delta \big\}$ does
not exceed $e^{-cr^4}$.
\end{lemma}

\par\noindent{\em Proof:} If $|f_j|\ge R^\delta $ somewhere in
the disk $4r\D$, then (at the same point) the corresponding
potential $U$ is not less than $\delta \log R - \frac12 (4r)^2 =
\frac12 \delta \log R$ (due to the choice of $r$). By
Lemma~\ref{lemma4.1}, the probability of this event does not
exceed $C r^2 e^{-cR^\delta} \ll e^{-cr^4}$ if $R\gg 1$.
\nopagebreak\hfill$\Box$

\medskip
Hence to complete the proof of Proposition~\ref{subthm6.4}, we
need to estimate the probability of the event $ \Omega_j' $.

\begin{lemma}\label{sublemma6.16}
$\Pr {\Omega_j'} \le e^{-cr^4}$.
\end{lemma}

\par\noindent{\em Proof: }
Define the random set $ A = \left\{ t\in [0,1]\colon
\Omega_j(t)\text{ occurs}\right\}$, and let $ X = m_1(A) $. Then,
by Chebyshev's inequality,
\[
\mathbb{P}(\Omega_j') = \mathbb{P}(X\ge 1/2) \le 2 \mathbb{E}(X) =
2 \int_0^1 \mathbb{P}( \Omega_j(t))dt \le 2 \max_{t\in [0,1]}
\mathbb{P}(\Omega_j(t)).
\]
By Lemma~\ref{lemma6.a}, the maximum on the right-hand side does
not exceed $e^{-cr^4}$. \hfill $\Box$

\medskip This completes the proof of Proposition~\ref{subthm6.4}.

\section{Proof of the partition theorem}

Set $\Delta(z) = U_{xx}(z)U_{yy}(z) - U_{xy}^2(z)$.

\begin{lemma}
For $z\in\C\setminus\mathcal Z_f$, $\displaystyle \Delta(z) = 1 -
\left| \left( \frac{f'}{f}\right)'(z)\right| $.
\end{lemma}

\par\noindent{\em Proof: } is a straightforward computation.
Since $\partial_x = \partial_z + \partial_{\bar z}$, $\partial_y =
i (\partial_z - \partial_{\bar z})$, we have
\[
\partial_{xx} = \partial_{zz} + 2\partial_{z\bar z} + \partial_{\bar z \bar
z}, \qquad \partial_{yy} = -\partial_{zz} + 2\partial_{z\bar z} -
\partial_{\bar z \bar z}, \qquad \partial_{xy} = i(\partial_{zz} -\partial_{\bar z
\bar z} )\,.
\]
Whence
\[
\Delta (z) = \left( U_{zz} + 2U_{z\bar z} + U_{\bar z \bar
z}\right) \left( - U_{zz} + 2U_{z\bar z} -U_{\bar z \bar z}
\right) + \left( U_{zz} - U_{\bar z \bar z}\right)^2 = 4 \left(
U_{z\bar z}^2 - U_{zz} \cdot U_{\bar z \bar z} \right)\,.
\]
Taking into account that $\displaystyle U_{zz} = \frac12 \left(
\frac{f'}{f}\right)'$, $\displaystyle U_{\bar z \bar z} = \frac12
\overline{\left( \frac{f'}{f} \right)'}$, and  $\displaystyle U_{z
\bar z} = -\frac12$, we get the result. \hfill $\Box$

\medskip Denote by $\operatorname{Crit}U$ the set $\{z\in\C\colon \nabla
U(z) = 0 \}$ of critical points of the potential $U$.

\begin{lemma}\label{lemma.16} Almost surely, the following hold:

\smallskip\noindent (i) each critical point of $U$ is
non-degenerate; i.e., $ \Delta(w) \ne 0$ for $w\in
\operatorname{Crit}(U)$;

\smallskip\par\noindent (ii) the critical set $\operatorname{Crit}U$
has no finite accumulation points.
\end{lemma}

\par\noindent{\em Proof:}
Note that the probability that $w=0$ is a critical point of $U$ is
$0$. At the critical points $w\ne 0$ of $U$, we have $f(w)=
\overline w ^{-1} f'(w)$. Hence
$$
\Delta (w) = 1 - \left|\overline w \frac{f''}{f'}(w) - \overline w
^2 \right|^2\,, \qquad w\in\crit\,.
$$

Now, let us set $f(z) = \xi_0 + \xi_1z +h(z)$, where $h(z)$ is a
random entire function determined by $\xi_2,\xi_3,\dots$. Then, on
$\crit$, the determinant $\Delta(w)$ coincides with
$$
\Delta_1(w):=1 - \left|\overline w \frac{h''(w)}{h'(w)+\xi_1} -
\overline w ^2 \right|^2\,.
$$
Observe that, for each $w$ with $|w|\ne 1$ and each
$\xi_2,\xi_3,\dots$, the set of $\xi_1$ where the last expression
is $0$ has zero measure. Thus, using the Fubini theorem, we
conclude that for almost all $\xi_1,\xi_2,\dots$, the set
$\{w\in\C\,:\Delta_1(w)=0\}$ has zero measure. Let now $g(z)=\xi_1
z+h(z)$. If $g$ is fixed (i.e., $\xi_1$, $\xi_2$, ..., are fixed)
and $w\ne 0$ is a critical point of $U$, then $\xi_0$ is
determined by equation
$$
\xi_0 = \frac{g'(w)}{\overline w} - g(w)\,.
$$
The right hand side defines a real-analytic mapping of the
punctured plane $\mathbb C\setminus \{0\}$ and, therefore, it maps
sets of zero area in the $w$-plane to sets of zero area in the
$\xi_0$-plane. Hence, for almost every choice of the independent
coefficients $\xi_1$, $\xi_2$, ..., the set of $\xi_0$ for which
there exist degenerate critical points of $U$ has measure zero.
Using Fubini's theorem once more, we get the conclusion of
statement (i) of the lemma.

Statement (ii) follows from (i). The planar map given by $\nabla
U$ is real-analytic outside the set where $U$ equals $-\infty$.
Note that, unless $f$ identically equals $0$ (which is an event of
zero probability), the gradient $\nabla U(z)=\overline{\left(\frac
{f'}f (z)\right)}-z$ tends to $\infty$ at every zero of $f$ and,
therefore, no point of $U^{-1}(-\infty)=\mathcal Z_f$ can be an
accumulation point of $\crit$. Thus, if the set
$\operatorname{Crit}U$ has a finite accumulation point, then this
point itself belongs to $\crit$ and, by the inverse function
theorem, the map given by $\nabla U$ is degenerate at this point.
\nopagebreak\hfill $\square$
\medskip

It is worth mentioning that there is another way to prove
statement (ii) of Lemma~\ref{lemma.16} elaborating on the fact
that, if $g$ is an analytic function and the solutions of the
equation $g(z)=\bar z$ have a finite accumulation point, then $g$
must be a M\"obius transformation.

\medskip

\begin{lemma}\label{lemma.17} Almost surely, the
following hold:

\smallskip\noindent (i) each oriented curve $\Gamma$ has a starting
point $s(\Gamma)\in \operatorname{Crit} (U)$ and a terminating
point $t(\Gamma)\in U^{-1}\{-\infty\} \cup
\operatorname{Crit}(U)$;

\smallskip\noindent (ii) at any limiting point, the oriented
gradient curve $\Gamma $ is tangent to a straight line passing
through that point.
\end{lemma}

\par\noindent{\em Proof:} We refer the reader to \cite[Chapter~4]{Hur} for
the facts from the standard ODE theory we use.

\smallskip\noindent (i) It follows from the long gradient
curve theorem that, almost surely, gradient curves cannot escape
to or come from infinity. Now it remains to observe that the
limiting set $\mathcal L$ of any gradient curve $\Gamma$ is
contained in the set of singular points of the gradient flow; that
is, in the set $\crit\cup U^{-1}\{-\infty\}$. Hence, by
Lemma~\ref{lemma.16}, $\mathcal L$ consists of isolated points.

\smallskip\par\noindent (ii) The critical points of $U$ are either local
maxima or saddle points. By Lemma~\ref{lemma.16}, almost surely,
all of them are non-degenerate. The rest follows from the standard
ODE theory: the behaviour of the integral curves in a
neighbourhood of these points is the same as the behaviour of the
integral curves for the linear ODE obtained by discarding the
non-linear terms in the Taylor expansion of $\nabla U$.
\nopagebreak\hfill $\square$
\medskip

\begin{lemma}\label{lemma.18} Each gradient curve is real
analytic everywhere except at the limiting points.
\end{lemma}

\par\noindent{\em Proof:} $\nabla U$ is real analytic everywhere
except on the set where $U=-\infty$. Hence, by the Cauchy theorem,
the gradient curves are real analytic at all points where $\nabla
U\ne 0$. \nopagebreak\hfill $\square$

\medskip

Now we are ready to prove the partition
theorem~\ref{th.partition}. By the long gradient curve
theorem~\ref{th.longcurve}, almost surely, all the basins are
bounded. We call a gradient curve $\Gamma$ {\it singular} if
$t(\Gamma)\in \crit$. Note that, almost surely, every point that
is not in one of the basins must lie on a singular curve.
Moreover, with probability $1$, for every compact $K$ on the
complex plane, there exists another compact $\ti K$ such that all
gradient curves intersecting $K$ are contained in $\ti K$.
(Otherwise, there exists an $N\in\mathbb N$ such that, for any
integer $M
> N$, there is a gradient curve connecting $\partial Q(0, N)$ and
$\partial Q(0, M)$. The probability of this event is $0$.) Also, a
gradient curve cannot terminate at a local maximum of $U$ and each
saddle point of $U$ serves as a terminating point for $2$ singular
curves. This allows us to conclude that, almost surely, we may
have only finitely many singular curves intersecting any compact
subset of $\C$. In particular, almost surely, each basin $B(a)$ is
bounded by finitely many singular curves and their limiting
points, which is enough to justify the area computation in the
introduction. \nopagebreak\hfill $\square$

\section{The upper bounds  in theorem~\ref{th.distance}}

First, we prove a useful ``length and area estimate'' of
deterministic nature valid for Liouville  vector fields; that is,
the fields with constant divergence.  Then we derive the upper
bounds for the probability that a given point $z$ is far from its
sink $a_z$.

\subsection{The length and area estimate}

Consider the disk $D = \{|z-a|<\epsilon\}$. Since $\displaystyle
\nabla U(z) = \frac{z-a}{|z-a|^2} + O(1)$, as $z\to a$, we can fix
a sufficiently small $\epsilon>0$ such that each gradient curve
hits the boundary circumference $\{|z-a|=\epsilon\}$ only once.
This gives us a one-to-one correspondence between the points of
the circumference $\big\{ a + \epsilon e^{i\theta} \big\}$ and the
gradient curves in $B(a)$; i.e., the gradient curves are
parameterized  by the ``angular coordinate'' $\theta$.

By $D(t)$ we denote the pre-image of $D$ under the gradient flow
of $\nabla U$ for time $t$; i.e., if $\displaystyle \frac{dZ}{dt}
= -\nabla U\big( Z(t) \big)$, then $D(t) = \{z = Z(0)\colon
Z(t)\in D\}$. By $A(t)$ we denote the area of $B(a)\setminus
D(t)$. Since $\operatorname{div} (\nabla U) = -2$ on
$B(a)\setminus \{a\}$, the evolution of the area is very simple:
$\displaystyle \frac{dA}{dt} = -2A $. This is Liouville's theorem
(which follows from the divergence theorem), see, for instance,
\cite[\S 16]{Arn}.

We will need an ``infinitesimal version'' of this equation. The
boundary $\partial B(a)$ contains finitely many saddle points of
$U$. By $\alpha_1<\,...\,<\alpha_s<\alpha_{s+1}=\alpha_1+2\pi$ we
denote the angular coordinates of the gradient curves that connect
the saddle points on $\partial B(a)$ with the sink $a$. Take any
$\theta$ different from $\alpha_1, \,...\,,\alpha_s$, say
$\alpha_l < \theta < \alpha_{l+1}$, and choose $\theta_1$ and
$\theta_2$ such that $\alpha_l < \theta_1 < \theta < \theta_2 <
\alpha_{l+1}$. The gradient curves $\Gamma(\theta_1)$,
$\Gamma(\theta_2)$ must terminate at the same local maximum. They
bound a ``diangle'' $Y(\theta_1, \theta_2)$ with the vertices at
$a$ and at a local maximum. Consider the ``triangle'' $T(t;
\theta_1, \theta_2) = Y(\theta_1, \theta_2)\setminus D(t)$ and its
area $A(t; \theta_1, \theta_2) = m_2 T(t; \theta_1, \theta_2) $.
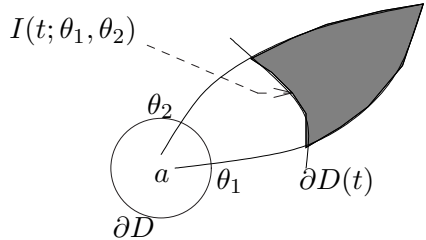
\begin{figure}[h]
\begin{center}
\setlength{\unitlength}{0.00078740in}
\begingroup\makeatletter\ifx\SetFigFont\undefined%
\gdef\SetFigFont#1#2#3#4#5{%
  \reset@font\fontsize{#1}{#2pt}%
  \fontfamily{#3}\fontseries{#4}\fontshape{#5}%
  \selectfont}%
\fi\endgroup%
{\renewcommand{\dashlinestretch}{30}
\begin{picture}(2858,1626)(0,-10)
\put(1005,519){\ellipse{656}{656}}
\dashline{60.000}(915,1329)(1635,1014)(1860,1014)
\path(1740.000,984.000)(1860.000,1014.000)(1740.000,1044.000)
\texture{44555555 55aaaaaa aa555555 55aaaaaa aa555555 55aaaaaa aa555555 55aaaaaa 
	aa555555 55aaaaaa aa555555 55aaaaaa aa555555 55aaaaaa aa555555 55aaaaaa 
	aa555555 55aaaaaa aa555555 55aaaaaa aa555555 55aaaaaa aa555555 55aaaaaa 
	aa555555 55aaaaaa aa555555 55aaaaaa aa555555 55aaaaaa aa555555 55aaaaaa }
\shade\path(1950,654)(2175,789)(2355,924)
	(2580,1194)(2670,1464)(2715,1599)
	(2445,1554)(2085,1464)(1860,1374)
	(1590,1239)(1725,1149)(1860,1014)
	(1950,879)(1950,789)(1950,654)
\path(1950,654)(2175,789)(2355,924)
	(2580,1194)(2670,1464)(2715,1599)
	(2445,1554)(2085,1464)(1860,1374)
	(1590,1239)(1725,1149)(1860,1014)
	(1950,879)(1950,789)(1950,654)
\path(1095,519)(1098,519)(1104,520)
	(1114,521)(1131,523)(1154,526)
	(1183,529)(1218,533)(1258,538)
	(1303,544)(1349,550)(1398,556)
	(1447,563)(1495,570)(1542,577)
	(1588,584)(1631,591)(1671,598)
	(1710,605)(1746,612)(1780,619)
	(1811,626)(1841,634)(1870,642)
	(1897,650)(1923,658)(1948,667)
	(1973,676)(1998,687)(2024,698)
	(2049,710)(2074,723)(2099,736)
	(2124,750)(2149,765)(2173,780)
	(2198,796)(2222,813)(2245,830)
	(2269,848)(2291,866)(2313,884)
	(2335,903)(2355,922)(2375,941)
	(2394,960)(2411,978)(2428,997)
	(2445,1015)(2460,1033)(2474,1051)
	(2488,1069)(2500,1087)(2513,1104)
	(2526,1125)(2539,1145)(2552,1167)
	(2563,1189)(2575,1212)(2587,1236)
	(2598,1262)(2609,1290)(2621,1319)
	(2633,1351)(2645,1383)(2656,1417)
	(2668,1451)(2679,1483)(2688,1514)
	(2697,1540)(2704,1562)(2709,1579)
	(2712,1590)(2714,1596)(2715,1599)
\path(1005,609)(1006,611)(1009,616)
	(1014,625)(1022,638)(1033,655)
	(1046,677)(1061,703)(1078,731)
	(1097,761)(1116,791)(1136,821)
	(1155,850)(1174,878)(1192,904)
	(1210,929)(1228,952)(1245,974)
	(1262,994)(1279,1013)(1296,1031)
	(1314,1048)(1331,1065)(1350,1081)
	(1368,1096)(1386,1111)(1405,1126)
	(1426,1141)(1447,1156)(1469,1171)
	(1491,1186)(1515,1201)(1540,1216)
	(1565,1231)(1591,1246)(1618,1261)
	(1645,1276)(1673,1290)(1701,1304)
	(1729,1317)(1756,1331)(1784,1343)
	(1812,1355)(1839,1367)(1866,1378)
	(1892,1389)(1918,1399)(1944,1408)
	(1970,1418)(1995,1426)(2020,1435)
	(2046,1443)(2072,1451)(2098,1459)
	(2126,1467)(2154,1475)(2184,1483)
	(2216,1491)(2249,1499)(2285,1507)
	(2322,1516)(2361,1525)(2401,1534)
	(2443,1543)(2484,1552)(2525,1560)
	(2564,1568)(2600,1576)(2632,1582)
	(2659,1588)(2681,1592)(2697,1595)
	(2707,1597)(2712,1598)(2715,1599)
\path(1455,1374)(1458,1372)(1463,1367)
	(1474,1358)(1488,1344)(1508,1327)
	(1530,1307)(1555,1285)(1581,1262)
	(1607,1239)(1631,1217)(1655,1196)
	(1676,1177)(1696,1159)(1714,1142)
	(1730,1127)(1745,1113)(1759,1099)
	(1772,1086)(1785,1074)(1800,1059)
	(1815,1043)(1829,1028)(1844,1012)
	(1857,996)(1871,980)(1884,963)
	(1896,946)(1907,930)(1917,913)
	(1927,897)(1935,881)(1942,865)
	(1948,850)(1953,834)(1958,819)
	(1961,803)(1964,787)(1965,770)
	(1966,751)(1967,730)(1966,707)
	(1965,681)(1963,654)(1961,625)
	(1958,597)(1956,570)(1953,549)
	(1952,533)(1951,523)(1950,520)(1950,519)
\put(1365,384){\makebox(0,0)[lb]{\smash{{\SetFigFont{11}{13.2}{\rmdefault}{\mddefault}{\updefault}$\theta_1$}}}}
\put(960,429){\makebox(0,0)[lb]{\smash{{\SetFigFont{11}{13.2}{\rmdefault}{\mddefault}{\updefault}$a$}}}}
\put(915,879){\makebox(0,0)[lb]{\smash{{\SetFigFont{11}{13.2}{\rmdefault}{\mddefault}{\updefault}$\theta_2$}}}}
\put(15,1374){\makebox(0,0)[lb]{\smash{{\SetFigFont{11}{13.2}{\rmdefault}{\mddefault}{\updefault}$I(t; 
\theta_1, \theta_2)$}}}}
\put(1905,384){\makebox(0,0)[lb]{\smash{{\SetFigFont{11}{13.2}{\rmdefault}{\mddefault}{\updefault}$\partial 
D(t)$}}}}
\put(690,69){\makebox(0,0)[lb]{\smash{{\SetFigFont{11}{13.2}{\rmdefault}{\mddefault}{\updefault}$\partial 
D$}}}}
\end{picture}
}
\end{center}
\caption{``Triangle'' $T(t; \theta_1, \theta_2)$}
\end{figure}

By Green's theorem,
\begin{equation}\label{eq8.0}
-2 A(t; \theta_1, \theta_2) = \iint_{T(t; \theta_1, \theta_2)}
\Delta U \, dm_2 = \int_{\partial T(t; \theta_1, \theta_2)}
\langle \nabla U, n \rangle\, |dz|
\end{equation}
where $n$ is the unit normal directed outward the triangle. The
boundary $\partial T(t; \theta_1, \theta_2)$ consists of parts of
the gradient curves $\Gamma(\theta_1)$ and $\Gamma(\theta_2)$,
where $\displaystyle \frac{\partial U}{\partial n} = 0$, and of
the part $I=I(t; \theta_1, \theta_2)$ of the curve $\partial
D(t)$. If $\big\{ (x(t, \theta), y(t, \theta))\colon \theta_1\le
\theta \le \theta_2 \big\}$ is the equation of the arc $I$, then
at the point $\theta$ the unit normal $n$ is given by
$\displaystyle -\frac{(y_\theta, -x_\theta )}{\sqrt{x_\theta^2 +
y_\theta^2}}$. Hence $\displaystyle  A(t; \theta_1, \theta_2) =
\frac12 \int_{\theta_1}^{\theta_2} \big( U_x y_\theta - U_y
x_\theta \big) \, d\theta$, and we conclude that the area $A(t;
\theta_1, \theta_2)$ has a smooth angular density $S(t, \theta) =
\frac12 \left(  U_x y_\theta - U_y x_\theta \right) $. By
Liouville's theorem, $\displaystyle \frac{\partial A(t; \theta_1,
\theta_2)}{\partial t} = -2A(t; \theta_1, \theta_2)$. Therefore,
the density $S(t, \theta)$ satisfies the same differential
equation $\displaystyle \frac{\partial S(t, \theta)}{\partial t} =
-2S(t, \theta)$.

Now, we re-parameterize the gradient curve $\Gamma(\theta)$ by its
length $l$ starting at the sink $a$. We treat the restrictions of
the density $S$ and of the gradient $\nabla U$ to $\Gamma
(\theta)$ as functions of the length $l$; i.e., $S(l)=S(t(l),
\theta)$, and similarly for $\nabla U$. Note that $\displaystyle
\frac{dl}{dt} = |\nabla U|$. We arrive at the ordinary
differential equation for the density $S$:
\[
\frac{\partial S}{\partial l}\cdot |\nabla U| = -2S\,.
\]
Solving this equation, we get
\[
S(l) = S(l_0) \exp \left( -2 \int_{l_0}^l \frac{dl}{|\nabla U|
}\right)\,.
\]

Denote by $z_l$ the point on the gradient curve $\Gamma (\theta)$
that cuts the arc of length $l$ from that curve. By the
Cauchy-Schwartz inequality,
\[
\int_{l_0}^l \frac{dl}{|\nabla U|} \ge (l-l_0)^2
\left(\int_{l_0}^l |\nabla U| \, dl\right)^{-1} =
\frac{(l-l_0)^2}{U(z_l) - U(z_{l_0})}\,.
\]
We arrive at the crucial
\begin{proposition}\label{lemma8.a}
In the same notation as above,
\[
S(l) \le S(l_0) \exp \left( - \frac{2(l-l_0)^2}{U(z_l) -
U(z_{l_0})} \right)\,.
\]
\end{proposition}

\subsection{Distance to the sink (the upper bound)}

Fix $\delta\in (0, 2]$. We define the {\em tentacles} $T_R(a)$ of
the basin $B(a)$ as follows. Given $\theta$, we move along the
gradient curve $\Gamma (\theta)$ in the direction of growth of the
potential $U$, starting at the sink $a$, till we hit the point
where $U = -R^\delta$. After that, we keep on moving along $\Gamma
(\theta)$ the distance $R$ (measured along $\Gamma (\theta)$), and
then stop. The rest of the curve is called the $\theta$-tentacle.
The tentacles $T_R(a)$ are the union of all $\theta$-tentacles. Of
course, it may happen that the tentacles $T_R(a)$ are empty.

\medskip
Now, we are ready to estimate the probability that $|z - a_z|$ is
large. By translation invariance, this probability does not depend
on the choice of $z$, so we choose $z=0$. Suppose that $|a_0|>2R$.
We know that at least one of the following happens:
\begin{itemize}
\item[(i)] either the distance from $0$ to the curve $\Gamma_0 \cap
\{U<-R^\delta\}$ measured along $\Gamma_0$ is less than $R$;
\item[(ii)] or $0\in T_R(a_0)$.
\end{itemize}
(Recall that $\Gamma_0$ is the gradient curve that passes through
the origin.)

In the first case, the curve $\gamma = \Gamma_0 \cap
\{U<-R^\delta\}$ connects the circumferences $\big\{ |z|=R \big\}$
and $\big\{ |z|=2R \big\}$. By Theorem~\ref{thm4.3}, the
probability of this event does not exceed $Ce^{-cR^{1+\frac32
\delta}}$.

Now, we estimate the probability of the event (ii) . By
translation invariance,
\begin{multline}\label{eq8.d}
\pi \Pr { 0\in T_R(a_0) } = \iint_{\D} \Pr {w\in T_R(a_w)}\,
dm_2(w) \\
= \int_\Omega m_2 \{w\in \D\colon w\in T_R(a_w) \}\, d\mathbb P\,.
\end{multline}
Thus, we need to estimate the area of the random set $\{w\in
\D\colon w\in T_R(a_w) \}$; that is, the area of the union of all
possible tentacles within $\D$.

\medskip We throw away three exceptional events.
Let $\Omega_1$ be the event that there exists a gradient curve
connecting the circumferences $\{|z|=1\}$ and $\{|z|=R^2\}$. By
the long gradient curve theorem, $\Pr {\Omega_1} \le e^{-cR^2} $.
If $\Omega_1$ does not occur, then $|a|<R^2$, for any basin $B(a)$
that intersects the unit disk. Let $\Omega_2$ be the event that
there exists a gradient curve connecting the circumferences
$\{|z|=R^2\}$ and $\{|z|=2R^2\}$. Again, $\Pr {\Omega_2} \le
e^{-cR^2}$. If $\Omega_1$ and $\Omega_2$ do not occur, then any
basin $B$ that intersects the unit disk $\D$ is contained in the
disk $2R^2\D$. Recalling that each basin has area $\pi$ and
comparing the areas, we see that the number of such basins does
not exceed $4R^4$. At last, we exclude the event $\displaystyle
\Omega_3 = \{\max_{2R^2\D} U > R^\delta\}$. By
Lemma~\ref{lemma4.1}, $\Pr {\Omega_3} < CR^4 e^{-ce^{R^\delta}} <
e^{-cR^4}$ if $R$ is big enough.

Now, after throwing away these three events, we can estimate the
area of the random set $\big\{w\in \D\colon w\in T_R(a_w) \big\}$.
First, we bound the area of one tentacle $T_R(a)$. Since $U \le
R^\delta$ everywhere in $B(a)$, for each $\theta$-tentacle, we can
apply the length and area estimate from Proposition~\ref{lemma8.a}
with $l-l_0\ge R$ and $U(z_l)-U(z_{l_0}) \le 2R^\delta$.
Integrating over $\theta$, we get
\[
m_2 T_R(a)  \le m_2 B(a) e^{- R^{2-\delta}} = \pi e^{-
R^{2-\delta}}\,.
\]
The number of tentacles coming from  different basins and hitting
the unit disk $\D$ does not exceed $4R^4$. We conclude that if the
events $\Omega_i$, $1\le i \le 3$, do not occur, then
$m_2\big\{w\in \D\colon w\in T_R(a_w) \big\} \le 4\pi R^4 e^{-
R^{2-\delta}}$. In view of (\ref{eq8.d}), we see that the
probability of the event $\left\{ 0\in T_R(a_0) \right\} $ is
bounded by $ e^{-cR^{2-\delta}} $ if $R\gg 1$.

Thus,
\begin{multline*}
\Pr {|a_0| > 2R} < \Pr { \operatorname{diam} \big(\Gamma_0 \cap
\{U<-R^\delta\} \big) > R} + \Pr { 0\in T_R(a_0)}  \\
< C e^{-cR^{1+\frac32 \delta}} + C e^{-cR^{2-\delta}}\,.
\end{multline*}
Choosing $\delta=\frac25$, we complete the proof. \hfill $\Box$

\section{The lower bounds in theorems~\ref{th.diameter} and
\ref{th.distance}}

The proofs of the lower bounds for the diameter of the basin and
the distance to the sink are based on the same idea. The function
$\displaystyle \frac{z^n}{\sqrt{n!}}$ has a singular line $\big\{
|z|=\sqrt{n} \big\} $ where the gradient of its potential
vanishes. Then after any analytic perturbations small in the
annulus $\big\{ |z-\sqrt{n}|\le 2 \big\}$, this annulus still
contains plenty of long gradient curves and of points that are far
from their sinks.

\subsection{Diameter of the basin (the lower bound)}

We choose a big $R\gg 1$ such that $n=R^2$ is an integer and
consider the function $F(z)=\dfrac{z^n}{\sqrt{n!}}\left(1+\dfrac
z{10R}\right)$ in the domain
\[
\mathcal D = \left\{z\in\C\colon R-1<|z|<R+1, |\arg z-\frac\pi 2|
< \frac 1{10}\right\}\,.
\]
Note that, for the corresponding potential $U$, we have
$$
\nabla U(z)=\overline{\frac {F'}F(z)}-z=\frac n{\bar z}-z+\frac
1{10R}\frac 1{1+\frac {\bar z}{10R}}\,.
$$
Since the vector $\frac n{\bar z}-z$ is purely radial and the sine
of the angle between the vectors $1+\frac z{10R}$ and $z$ is at
least $\frac12$ for $z\in\mathcal D$, we see that the angular
component of $ -\nabla U $ is oriented counter-clockwise and its
size is at least $\frac 1{30R}$ in $\mathcal D$. Also, the
gradient field $ -\nabla U $ is directed outside the domain
$\mathcal D$ on the boundary arcs $\big\{ |z|=R\pm 1, \ |\arg
z-\frac\pi 2|<\frac 1{10} \big\}$, and the radial component of $
-\nabla U $ is at least $1$ on these arcs.

Thus, there is a gradient curve that starts at the right boundary
interval $\big\{ \arg z-\frac\pi 2 = -\frac 1{10}, \ R-1<|z|<R+1
\big\}$, and hits the point $iR$. Thereby, its diameter must be at
least $\frac{R}{20}$.
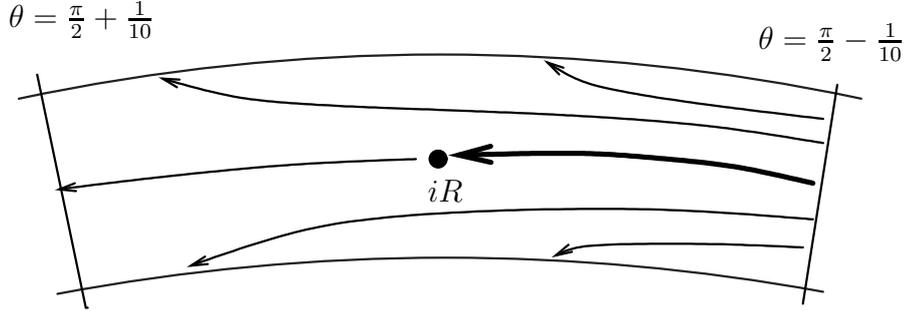
\begin{figure}[h]
\begin{center}
\setlength{\unitlength}{0.00083333in}
\begingroup\makeatletter\ifx\SetFigFont\undefined%
\gdef\SetFigFont#1#2#3#4#5{%
  \reset@font\fontsize{#1}{#2pt}%
  \fontfamily{#3}\fontseries{#4}\fontshape{#5}%
  \selectfont}%
\fi\endgroup%
{\renewcommand{\dashlinestretch}{30}
\begin{picture}(6895,1862)(0,-10)
\allinethickness{1.000pt}%
\put(2662.000,-10841.102){\arc{24856.218}{4.5016}{4.9231}}
\put(2697.324,-12551.892){\arc{25752.067}{4.5251}{4.8942}}
\put(2650,937){\blacken\ellipse{102}{102}}
\put(2650,937){\ellipse{102}{102}}
\path(175,1462)(475,12)(487,12)
\path(5112,1399)(4900,37)
\allinethickness{2.000pt}%
\path(4962,787)(4958,788)(4951,790)
    (4937,793)(4917,797)(4891,803)
    (4861,810)(4827,818)(4792,826)
    (4756,834)(4722,841)(4690,848)
    (4660,855)(4632,861)(4606,866)
    (4582,871)(4559,875)(4537,879)
    (4516,883)(4495,887)(4475,890)
    (4454,893)(4433,896)(4412,899)
    (4390,902)(4368,905)(4345,908)
    (4322,911)(4299,914)(4276,916)
    (4253,919)(4231,921)(4208,923)
    (4186,926)(4165,928)(4144,930)
    (4124,932)(4104,933)(4085,935)
    (4066,937)(4047,939)(4028,940)
    (4009,942)(3989,944)(3969,946)
    (3948,947)(3927,949)(3906,951)
    (3884,953)(3863,955)(3841,956)
    (3819,958)(3798,960)(3777,961)
    (3757,962)(3737,964)(3718,965)
    (3699,966)(3680,967)(3662,968)
    (3644,969)(3625,969)(3606,970)
    (3587,971)(3567,971)(3546,972)
    (3525,972)(3504,972)(3482,973)
    (3460,973)(3438,973)(3416,973)
    (3394,973)(3372,973)(3350,973)
    (3329,973)(3309,973)(3288,972)
    (3268,972)(3248,972)(3227,972)
    (3206,971)(3184,971)(3162,971)
    (3137,970)(3111,970)(3083,969)
    (3052,968)(3019,968)(2985,967)
    (2949,966)(2913,965)(2879,964)
    (2847,964)(2821,963)(2775,962)
\path(2973.866,1016.335)(2775.000,962.000)(2976.039,916.359)
\allinethickness{1.000pt}%
\path(5025,1037)(5023,1037)(5017,1038)
    (5008,1040)(4993,1042)(4973,1045)
    (4948,1050)(4918,1054)(4884,1060)
    (4847,1066)(4809,1072)(4770,1078)
    (4730,1084)(4692,1090)(4655,1096)
    (4619,1102)(4585,1107)(4553,1111)
    (4522,1116)(4492,1120)(4464,1124)
    (4436,1128)(4409,1131)(4382,1135)
    (4356,1138)(4329,1141)(4304,1144)
    (4278,1146)(4252,1149)(4226,1152)
    (4198,1155)(4171,1157)(4142,1160)
    (4112,1163)(4082,1166)(4052,1168)
    (4020,1171)(3988,1174)(3956,1176)
    (3923,1179)(3890,1181)(3857,1184)
    (3824,1186)(3791,1188)(3758,1191)
    (3726,1193)(3694,1195)(3662,1197)
    (3630,1199)(3599,1201)(3568,1203)
    (3538,1204)(3507,1206)(3477,1208)
    (3448,1209)(3420,1211)(3390,1212)
    (3361,1214)(3330,1215)(3299,1217)
    (3268,1218)(3236,1220)(3203,1221)
    (3169,1223)(3135,1224)(3101,1226)
    (3066,1227)(3031,1228)(2995,1230)
    (2959,1231)(2924,1233)(2888,1234)
    (2853,1236)(2818,1237)(2783,1238)
    (2749,1240)(2715,1241)(2681,1242)
    (2648,1243)(2615,1245)(2583,1246)
    (2551,1247)(2519,1248)(2488,1249)
    (2456,1250)(2424,1251)(2391,1253)
    (2359,1254)(2326,1255)(2292,1256)
    (2258,1258)(2223,1259)(2188,1261)
    (2153,1262)(2118,1264)(2082,1265)
    (2046,1267)(2010,1269)(1975,1271)
    (1939,1273)(1905,1274)(1870,1277)
    (1837,1279)(1804,1281)(1772,1283)
    (1741,1285)(1710,1288)(1681,1290)
    (1653,1292)(1626,1295)(1599,1297)
    (1574,1300)(1549,1303)(1525,1306)
    (1496,1309)(1468,1313)(1440,1317)
    (1412,1322)(1385,1327)(1357,1332)
    (1329,1338)(1300,1344)(1270,1351)
    (1239,1359)(1207,1366)(1174,1375)
    (1141,1384)(1108,1392)(1077,1401)
    (1047,1409)(1020,1417)(998,1423)
    (979,1428)(950,1437)
\path(1052.916,1431.237)(950.000,1437.000)(1038.096,1383.483)
\path(5025,1187)(5023,1187)(5018,1188)
    (5008,1189)(4994,1191)(4975,1193)
    (4951,1196)(4922,1199)(4889,1203)
    (4854,1207)(4816,1212)(4778,1217)
    (4741,1221)(4703,1226)(4667,1230)
    (4633,1234)(4600,1238)(4568,1242)
    (4538,1246)(4509,1250)(4482,1253)
    (4455,1257)(4428,1260)(4402,1263)
    (4376,1267)(4350,1270)(4323,1274)
    (4296,1278)(4268,1281)(4240,1285)
    (4212,1289)(4182,1293)(4152,1298)
    (4122,1302)(4091,1307)(4060,1311)
    (4028,1316)(3997,1321)(3966,1326)
    (3935,1331)(3905,1336)(3876,1342)
    (3847,1347)(3819,1352)(3792,1357)
    (3766,1362)(3742,1367)(3718,1372)
    (3696,1377)(3675,1381)(3655,1386)
    (3635,1391)(3609,1398)(3585,1405)
    (3562,1413)(3540,1421)(3518,1430)
    (3496,1439)(3475,1449)(3453,1460)
    (3431,1472)(3410,1485)(3389,1497)
    (3370,1508)(3354,1519)(3325,1537)
\path(3423.148,1505.505)(3325.000,1537.000)(3396.780,1463.023)
\path(2512,937)(2510,937)(2505,937)
    (2496,936)(2482,936)(2463,935)
    (2439,934)(2409,933)(2376,932)
    (2339,931)(2299,929)(2258,927)
    (2216,926)(2175,924)(2134,922)
    (2094,921)(2056,919)(2019,917)
    (1985,915)(1951,914)(1920,912)
    (1889,910)(1859,909)(1831,907)
    (1802,905)(1774,903)(1747,902)
    (1719,899)(1691,897)(1662,895)
    (1633,893)(1604,891)(1573,888)
    (1543,885)(1511,883)(1480,880)
    (1447,877)(1415,874)(1381,871)
    (1348,868)(1315,864)(1281,861)
    (1248,858)(1215,855)(1183,851)
    (1151,848)(1120,845)(1089,841)
    (1060,838)(1031,835)(1003,832)
    (977,829)(951,826)(925,824)
    (901,821)(877,818)(847,815)
    (819,811)(790,808)(762,805)
    (733,801)(703,798)(673,794)
    (641,790)(608,786)(574,782)
    (538,778)(502,773)(466,769)
    (432,764)(400,760)(372,757)
    (350,754)(312,749)
\path(407.884,786.832)(312.000,749.000)(414.407,737.259)
\path(4962,562)(4960,562)(4956,563)
    (4948,563)(4935,564)(4918,566)
    (4895,567)(4867,570)(4834,572)
    (4796,575)(4754,579)(4710,582)
    (4663,585)(4614,589)(4565,593)
    (4516,596)(4467,599)(4419,602)
    (4372,605)(4326,608)(4281,611)
    (4238,613)(4195,615)(4153,617)
    (4112,619)(4072,620)(4032,622)
    (3991,623)(3951,624)(3909,625)
    (3867,626)(3825,626)(3790,627)
    (3754,627)(3717,627)(3680,627)
    (3642,627)(3603,627)(3563,627)
    (3522,627)(3480,626)(3438,626)
    (3395,625)(3350,625)(3305,624)
    (3260,623)(3214,622)(3167,621)
    (3120,619)(3073,618)(3025,616)
    (2978,614)(2930,613)(2883,611)
    (2836,608)(2789,606)(2743,604)
    (2697,601)(2652,598)(2608,596)
    (2564,593)(2522,590)(2480,587)
    (2439,583)(2400,580)(2361,577)
    (2323,573)(2286,569)(2250,565)
    (2215,561)(2181,557)(2148,553)
    (2107,548)(2067,542)(2028,536)
    (1990,529)(1952,522)(1915,515)
    (1877,507)(1840,499)(1802,490)
    (1763,480)(1724,470)(1684,459)
    (1643,447)(1601,434)(1558,421)
    (1515,408)(1472,394)(1428,380)
    (1386,365)(1345,352)(1306,338)
    (1269,326)(1236,314)(1207,304)
    (1183,295)(1163,288)(1148,282)(1125,274)
\path(1211.237,330.464)(1125.000,274.000)(1227.663,283.240)
\path(4900,387)(4897,387)(4891,387)
    (4880,388)(4863,389)(4841,390)
    (4814,391)(4783,393)(4750,394)
    (4715,396)(4679,397)(4645,399)
    (4611,400)(4578,401)(4547,402)
    (4517,404)(4488,404)(4460,405)
    (4432,406)(4404,407)(4375,407)
    (4346,408)(4320,408)(4294,409)
    (4267,409)(4239,409)(4211,410)
    (4181,410)(4151,410)(4120,410)
    (4088,410)(4056,410)(4023,410)
    (3991,410)(3959,410)(3927,409)
    (3895,409)(3864,409)(3834,408)
    (3805,407)(3777,407)(3750,406)
    (3725,405)(3700,404)(3678,403)
    (3656,402)(3636,401)(3617,399)
    (3584,397)(3555,394)(3530,390)
    (3506,386)(3484,380)(3463,374)
    (3443,368)(3424,360)(3407,353)(3375,337)
\path(3453.262,404.082)(3375.000,337.000)(3475.623,359.361)
\put(0,1687){\makebox(0,0)[lb]{{\SetFigFont{12}{14.4}{\rmdefault}{\mddefault}{\updefault}$\theta
= \frac{\pi}2 + \frac1{10}$}}}
\put(2587,674){\makebox(0,0)[lb]{{\SetFigFont{12}{14.4}{\rmdefault}{\mddefault}{\updefault}$iR$}}}
\put(4562,1549){\makebox(0,0)[lb]{{\SetFigFont{12}{14.4}{\rmdefault}{\mddefault}{\updefault}
$\theta = \frac{\pi}2 - \frac1{10}$}}}
\end{picture}
}
\end{center}
\caption{The field $-\nabla U$ in $\mathcal D$}
\end{figure}

This conclusion will be preserved if, instead of the function $F$,
we consider its analytic perturbation $F+H$ with $H$ satisfying
$\left|\frac HF\right|\le R^{-2}$ in the annulus $\big\{ R-2\le
|z|\le R+2 \big\}$. Indeed, the absolute value of the perturbation
of $\nabla U$ the function $H$ creates in $\mathcal D$ is only
$\left|\frac{(F+H)'}{F+H} - \frac{F'}{F} \right| =
\left|\frac{(H/F)'}{1+(H/F)}\right|\le
\frac{R^{-2}}{1-R^{-2}}<\frac 1{60 R}$ for $R\gg 1$, which is too
small to change anything in the above picture.

Now it remains to estimate from below the probability of the event
that a G.E.F. $f$ is such a perturbation of $F$.
\begin{lemma}\label{lemma8.1}
If $R\gg 1$, then $\displaystyle \Pr {\max_{R-2 \le |z| \le R+2}
\left| \frac{f}{F}(z) -1 \right| \le \frac1{R^2}} \ge e^{-CR(\log
R)^{3/2}}$.
\end{lemma}
\noindent{\em Proof:} We write $f(z)=F(z)+H(z)$ where
$$
H(z)=\sum_{k\colon k\ne n,n+1}\xi_k\frac
{z^k}{\sqrt{k!}}+(\xi_n-1)\frac {z^n}{\sqrt{n!}}+
\left(\xi_{n+1}-\frac{\sqrt{n+1}}{10 R}\right)\frac
{z^{n+1}}{\sqrt{(n+1)!}}\,.
$$
Since $\displaystyle |F(z)|\ge \frac1{2}\frac{|z|^n}{\sqrt{n!}}$
in the annulus $R-2 \le |z| \le R+2$, it is enough to estimate
from below the probability of the event that
$$
\max_{R-2\le |z|\le R+2}\Bigl[\sum_{k\colon k\ne
n,n+1}|\xi_k|\frac {\sqrt{n!}}{\sqrt{k!}}|z|^{k-n}+|\xi_n-1|+
\left|\xi_{n+1}-\frac{\sqrt{n+1}}{10 R}\right|\frac
{|z|}{\sqrt{n+1}}\Bigr]<R^{-3}\,,
$$
say. Now, let us handle $\xi_n$ and $\xi_{n+1}$ first. We just
demand that the corresponding terms be both less than $R^{-4}$. It
is not hard to see that the probability of this event is about
$R^{-16}$. We may neglect it since the factor $R^{-16}$ does not
affect the lower bound $\displaystyle e^{-CR(\log R)^{3/2}}$ we
are trying to get. The remaining sum can be estimated as
$$
\sum_{k\colon k<n}|\xi_k|\frac
{\sqrt{n!}}{\sqrt{k!}}(R-2)^{k-n}+\sum_{k\colon k>n+1}|\xi_k|\frac
{\sqrt{n!}}{\sqrt{k!}}(R+2)^{k-n}\,.
$$
We shall show how to estimate from below the probability that the
second sum is less than $R^{-4}$. The estimate for the first sum
is very similar and we omit it (note that the corresponding events
depend on different $\xi_k$ and, therefore, are independent, so
the probability that both sums are small is just the product of
the probabilities that each of them is small). Let $k=n+m$,
$m=2,3,\dots$. We choose some big constant $A\gg 1$ and split the
sum into two: $\displaystyle S_1=\sum_{2\le m\le AR\sqrt{\log R}}$
and $\displaystyle S_2=\sum_{m > AR\sqrt{\log R}}$. We shall show
that the probability that $|S_2|<R^{-5}$ is very close to $1$ and
the probability that $|S_1|<R^{-5}\sqrt{\log R}$ is at least
$e^{-C R(\log R)^{3/2}}$.

To estimate $S_2$, we would like to use Lemma~\ref{lemma2.1} . To
this end, we need to estimate the sum
$$
\sum_{m > AR\sqrt{\log
R}}\frac{(R+2)^m}{\sqrt{(n+1)(n+2)\dots(n+m)}}\,.
$$
Note that, starting with $m=n$, the terms in this sum decay like a
geometric progression, more precisely, the ratio of each term to
the previous one is  $\frac
{R+2}{\sqrt{n+m+1}}<\frac{R+2}{\sqrt{2n}}=\frac{R+2}{\sqrt
2\,R}<\frac34$ if $R$ is large enough. Thus, it is enough to
estimate the sum over $m$ such that $AR\sqrt{\log R}<m\le n$. Now,
for $k=1,2,\dots,n$, we have $\displaystyle n+k\ge n e^{\frac
k{2n}}$. Thus, the $m$-th term of our sum does not exceed
\begin{equation}\label{eq8.*}
\left(1+\frac 2R\right)^m\prod_{k=1}^m e^{-\frac k{4n}}\le
e^{\frac {2m}{R}}e^{-\frac{m^2}{8n}} = e^{\frac
{2m}{R}}e^{-\frac{m^2}{8R^2}}\le C e^{-\frac{m^2}{16 R^2}}\,,
\qquad 1\le m \le n\,,
\end{equation}
and the whole sum does not exceed
\begin{multline*}
C\sum_{m>AR\sqrt{\log R}}e^{-\frac{m^2}{16 R^2}}\le
C\int_{AR\sqrt{\log R}-1}^{\infty} e^{-\frac{t^2}{16R^2}}\,dt
\\
\le CR\int_{\frac12 A\sqrt{\log R}}^{\infty}
e^{-\frac{t^2}{16}}\,dt\le CR e^{-\frac {A^2}{64}\log R}=
CR^{1-\frac{A^2}{64}}<R^{-6}
\end{multline*}
if $R$ is large enough. Thus, according to Lemma~\ref{lemma2.1},
the probability that $|S_2|<R^{-5}$ is very close to $1$ and, at
least, greater than $\frac 12$.

As to $S_1$, we just demand that each term in $S_1$ be less than
$R^{-6}$ (then $|S_1|<AR^{-5}\sqrt{\log R}$). Since the
coefficients
$\frac{(R+2)^m}{\sqrt{(n+1)(n+2)\dots(n+m)}}\overset{\eqref{eq8.*}}\le
C$, it is enough to demand that $|\xi_k|<C^{-1}R^{-6}$ for $n+2\le
k\le n + AR\sqrt{\log R}$. But the probability of this event is at
least $(cR^{-12})^{AR\sqrt{\log R}}\ge e^{-CR(\log R)^{3/2}}$.
This proves the lemma. \hfill $\Box$

\bigskip Thus, with probability $e^{-CR(\log R)^{3/2}}$, the point
$z=iR$ belongs to a basin of diameter greater than $\frac R{10}$.
It remains to note that, due to shift invariance of $U$, the same
is true for any other point $z$ on the complex plane. This proves
the lower bound in the diameter of the basin theorem. \hfill
$\Box$

\subsection{Distance to the sink (the lower bound)}

We choose a big $R\gg 1$ such that $n=R^2$ is an integer. This
time we start with the function $F(z)=\dfrac{z^n}{\sqrt{n!}}
e^{zR^{\delta-1} - R^\delta}$ with $0<\delta < 1$ (later, we'll
choose $\delta=\frac25$). The gradient of the corresponding
potential $U$ equals
$$
\nabla U(z) = \overline{\frac {F'}F(z)}-z=\frac n{\bar z} - z +
R^{\delta-1} = \frac{R^2-r^2}{r}\,e^{i\theta} + R^{\delta-1}\,,
\qquad z=re^{i\theta}\,.
$$
Let $\mathcal A = \{R-1 < |z| < R+1\}$, and let $\mathcal
D=\{z\in\mathcal A, |\arg z-\frac\pi 2|<\frac 1{10}\}$ be the same
sector as above. Note the following properties of the gradient
field:
\begin{itemize}
\item[(i)] on the boundary circumferences $|z|=R\pm 1$, the radial
component of the field $-\nabla U$ is directed  outward $\mathcal
A$, and its size is not less than $1$; inside $\mathcal A$, the
size of the radial component does not exceed $3$;
\item[(ii)] the field $-\nabla U$ has the horizontal drift
$R^{\delta-1}$ oriented to the left; in particular, inside the
sector $\mathcal D$, the angular component of the field $-\nabla
U$ is oriented counter-clockwise and its size is within the range
$[\frac12 R^{\delta-1}, 2R^{\delta-1} ]$.
\end{itemize}

By $\mathcal G$ we denote the set of points that hit the segment
\[ \displaystyle J = [ (R-1)e^{i\left( \frac{\pi}2+\frac1{10}
\right)}, (R+1)e^{i\left( \frac{\pi}2+\frac1{10} \right)}] \] when
moving along their trajectories. Because of the ``left-oriented
horizontal drift'' of the field $-\nabla U$, the points $z$ with
$\pi \ge |\arg z| > \frac{\pi}{2} + \frac1{10}$ cannot appear
within $\mathcal G$. (In fact, it is easy to see that $\mathcal
G\subset \{z\colon 0<\arg z < \frac{\pi}{2}+\frac1{10}\}$ but we
will not need this). By $\mathcal G(\theta)$ we denote the subset
of $\mathcal G$ that is located clock-wise with respect to the
segment \[ \displaystyle J(\theta) = [ (R-1)e^{i\left(
\frac{\pi}2+\frac1{10}-\theta\right)}, (R+1)e^{i\left(
\frac{\pi}2+\frac1{10}-\theta\right)}]\,. \] Note that $\mathcal
G(\theta_2) \subset \mathcal G(\theta_1) $ for
$\theta_2>\theta_1$. By $A(\theta)$ we denote the area of
$\mathcal G(\theta)$. We denote by $h(\theta)$ the length of the
intersection of the domain $\mathcal G$ with the segment
$J(\theta)$; i.e., the length of the ``left boundary wall'' of the
domain $\mathcal G(\theta)$.

\begin{lemma}\label{lemma8.2}
If $R \gg 1$, then $h(\theta)\ge e^{-CR^{2-\delta}}$ for
$0\le\theta\le \frac15$, and \[ m_2\big(\mathcal G \cap
\big\{\frac{\pi}2-\frac1{10} < \arg z < \frac{\pi}2 \big\}\big)
\ge e^{-CR^{2-\delta}}\,.
\]
\end{lemma}
\begin{figure}
\begin{center}
\setlength{\unitlength}{0.00083333in}
\begingroup\makeatletter\ifx\SetFigFont\undefined%
\gdef\SetFigFont#1#2#3#4#5{%
  \reset@font\fontsize{#1}{#2pt}%
  \fontfamily{#3}\fontseries{#4}\fontshape{#5}%
  \selectfont}%
\fi\endgroup%
{\renewcommand{\dashlinestretch}{30}
\begin{picture}(5891,2050)(0,-10)
\allinethickness{1.000pt}%
\put(2503.592,-10985.269){\arc{25091.961}{4.5121}{4.9851}}
\put(2546.928,-12777.859){\arc{26675.998}{4.5389}{4.9420}}
\path(4720,1372)(4570,410) \path(495,1372)(645,447)
\path(2133,1560)(2170,572)
\put(4220,60){\makebox(0,0)[lb]{{\SetFigFont{12}{14.4}{\rmdefault}{\mddefault}{\updefault}
$\frac{\pi}2 - \frac1{10}$}}}
\put(470,85){\makebox(0,0)[lb]{{\SetFigFont{12}{14.4}{\rmdefault}{\mddefault}{\updefault}
$\frac{\pi}2 + \frac1{10}$}}}
\put(308,822){\makebox(0,0)[lb]{{\SetFigFont{12}{14.4}{\rmdefault}{\mddefault}{\updefault}$J$}}}
\put(945,947){\makebox(0,0)[lb]{{\SetFigFont{12}{14.4}{\rmdefault}{\mddefault}{\updefault}$\mathcal
G$}}}
\put(1758,685){\makebox(0,0)[lb]{{\SetFigFont{12}{14.4}{\rmdefault}{\mddefault}{\updefault}$J(\theta)$}}}
\put(1083,1835){\makebox(0,0)[lb]{{\SetFigFont{12}{14.4}{\rmdefault}{\mddefault}{\updefault}$\theta$}}}
\put(2528.128,-12224.508){\arc{27980.924}{4.5629}{4.6832}}
\path(2020.874,1731.734)(2120.000,1760.000)(2019.234,1781.707)
\path(540.311,1649.253)(445.000,1610.000)(547.570,1599.782)
\thinlines \texture{aaffffff ffaaaaaa aaffffff ffaaaaaa aaffffff
ffaaaaaa aaffffff ffaaaaaa
    aaffffff ffaaaaaa aaffffff ffaaaaaa aaffffff ffaaaaaa aaffffff ffaaaaaa
    aaffffff ffaaaaaa aaffffff ffaaaaaa aaffffff ffaaaaaa aaffffff ffaaaaaa
    aaffffff ffaaaaaa aaffffff ffaaaaaa aaffffff ffaaaaaa aaffffff ffaaaaaa }
\shade\path(2133,1210)(2158,1010)(2383,1010)
    (3145,1022)(3708,997)(4170,960)
    (4633,922)(4658,997)(4195,1072)
    (3745,1122)(3295,1172)(2783,1185)
    (2145,1222)(2133,1210)
\path(2133,1210)(2158,1010)(2383,1010)
    (3145,1022)(3708,997)(4170,960)
    (4633,922)(4658,997)(4195,1072)
    (3745,1122)(3295,1172)(2783,1185)
    (2145,1222)(2133,1210)
\allinethickness{2.000pt}%
\path(508,1385)(510,1384)(515,1381)
    (524,1377)(538,1370)(555,1362)
    (576,1352)(599,1341)(624,1329)
    (649,1318)(675,1308)(699,1298)
    (723,1289)(747,1281)(770,1273)
    (793,1267)(816,1261)(839,1256)
    (863,1251)(889,1247)(908,1244)
    (927,1242)(948,1239)(970,1237)
    (992,1235)(1016,1233)(1040,1231)
    (1066,1229)(1093,1227)(1120,1225)
    (1149,1224)(1179,1222)(1209,1221)
    (1240,1220)(1272,1219)(1304,1218)
    (1337,1217)(1369,1216)(1402,1215)
    (1435,1215)(1468,1214)(1501,1214)
    (1534,1213)(1567,1213)(1599,1213)
    (1632,1212)(1664,1212)(1697,1212)
    (1726,1212)(1756,1212)(1785,1212)
    (1816,1211)(1847,1211)(1879,1211)
    (1912,1211)(1946,1211)(1980,1210)
    (2015,1210)(2051,1210)(2087,1209)
    (2125,1209)(2163,1208)(2201,1208)
    (2239,1207)(2278,1207)(2318,1206)
    (2357,1205)(2396,1205)(2435,1204)
    (2474,1203)(2513,1202)(2551,1201)
    (2590,1199)(2627,1198)(2665,1197)
    (2702,1196)(2739,1194)(2776,1193)
    (2813,1191)(2850,1189)(2882,1188)
    (2916,1186)(2949,1184)(2983,1182)
    (3018,1180)(3053,1178)(3088,1176)
    (3124,1173)(3161,1171)(3198,1168)
    (3236,1165)(3275,1162)(3314,1159)
    (3353,1156)(3392,1153)(3432,1150)
    (3472,1146)(3512,1143)(3552,1139)
    (3591,1135)(3631,1131)(3670,1128)
    (3709,1124)(3747,1120)(3785,1116)
    (3822,1112)(3859,1108)(3895,1104)
    (3930,1100)(3965,1095)(3999,1091)
    (4033,1087)(4066,1083)(4098,1079)
    (4130,1075)(4162,1070)(4197,1065)
    (4232,1060)(4267,1055)(4303,1050)
    (4338,1045)(4373,1039)(4409,1034)
    (4445,1028)(4481,1022)(4517,1016)
    (4554,1009)(4590,1003)(4626,997)
    (4663,990)(4699,983)(4734,976)
    (4770,970)(4805,963)(4839,956)
    (4873,949)(4906,942)(4938,936)
    (4970,929)(5000,922)(5030,916)
    (5059,909)(5087,903)(5115,897)
    (5142,891)(5168,884)(5193,878)
    (5218,872)(5254,863)(5290,854)
    (5325,845)(5359,836)(5393,827)
    (5427,818)(5460,809)(5492,800)
    (5522,792)(5552,783)(5579,775)
    (5605,768)(5629,760)(5650,754)
    (5669,748)(5685,743)(5698,738)
    (5709,734)(5718,731)(5723,729)
    (5726,727)(5727,726)(5724,726)
    (5719,727)(5710,729)(5699,731)
    (5684,734)(5667,738)(5646,743)
    (5624,748)(5598,754)(5571,760)
    (5542,767)(5511,774)(5480,781)
    (5448,788)(5415,795)(5383,802)
    (5350,809)(5318,816)(5287,822)
    (5256,828)(5233,833)(5211,837)
    (5188,841)(5165,845)(5141,849)
    (5117,853)(5092,856)(5067,860)
    (5041,864)(5014,868)(4986,872)
    (4957,875)(4928,879)(4898,883)
    (4867,886)(4835,890)(4804,893)
    (4771,897)(4738,900)(4705,904)
    (4671,907)(4637,910)(4603,913)
    (4568,916)(4533,919)(4497,922)
    (4460,925)(4422,928)(4394,931)
    (4364,933)(4334,935)(4304,938)
    (4272,940)(4239,942)(4205,945)
    (4171,947)(4135,950)(4099,953)
    (4061,955)(4023,958)(3984,960)
    (3944,963)(3903,965)(3862,968)
    (3820,971)(3778,973)(3735,976)
    (3692,978)(3649,980)(3606,983)
    (3563,985)(3520,987)(3478,989)
    (3435,991)(3393,993)(3352,995)
    (3311,997)(3270,998)(3230,1000)
    (3191,1001)(3152,1003)(3113,1004)
    (3075,1005)(3038,1006)(3001,1007)
    (2964,1008)(2926,1009)(2887,1009)
    (2849,1010)(2810,1010)(2771,1010)
    (2732,1011)(2693,1011)(2654,1011)
    (2614,1010)(2575,1010)(2535,1010)
    (2494,1009)(2454,1008)(2414,1007)
    (2374,1006)(2334,1005)(2294,1004)
    (2254,1002)(2215,1001)(2176,999)
    (2138,997)(2101,995)(2064,993)
    (2028,991)(1992,989)(1958,987)
    (1924,984)(1891,982)(1859,979)
    (1828,976)(1797,973)(1768,970)
    (1739,967)(1711,964)(1683,961)
    (1656,958)(1622,953)(1588,949)
    (1554,944)(1521,939)(1488,933)
    (1455,928)(1423,922)(1390,915)
    (1358,909)(1326,902)(1295,895)
    (1264,887)(1234,880)(1205,872)
    (1176,864)(1148,856)(1121,847)
    (1095,839)(1071,830)(1047,821)
    (1025,813)(1004,804)(984,795)
    (965,787)(947,778)(930,769)
    (914,760)(899,751)(877,737)
    (857,722)(838,707)(820,690)
    (802,671)(786,651)(769,629)
    (752,604)(736,578)(720,551)
    (704,524)(690,499)(679,476)
    (669,458)(663,446)(660,438)(658,435)
\end{picture}
}
\end{center}
\caption{The sets $\mathcal G$ and $\mathcal G(\theta)$}
\end{figure}
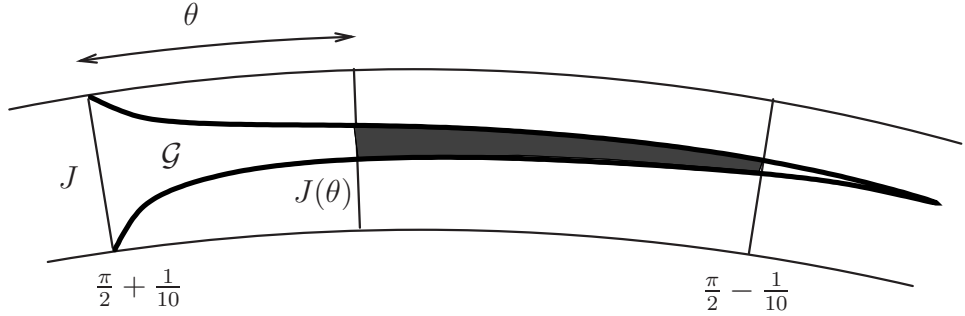

\par\noindent{\em Proof:} Note that the second estimate follows
from the first one by integration over $\theta$. We have
\[
A(\theta) = - \frac12 \iint_{\mathcal G(\theta)} \Delta U \, dm_2
= - \frac12 \int_{\partial\mathcal G(\theta) } \frac{\partial
U}{\partial n}\, |dz| = - \frac12 \int_{\mathcal G \cap J(\theta)
} \frac{\partial U}{\partial n}\, |dz|
\]
(since the rest of the boundary of $\mathcal G(\theta)$ consists
of gradient curves). In view of (ii),
\[
\frac12 R^{\delta-1} h(\theta) \le  - \int_{\mathcal G \cap
J(\theta) } \frac{\partial U}{\partial n}\, |dz| \le 2
R^{\delta-1} h(\theta)\,,
\]
whence
\begin{equation}\label{eq8.b}
\frac14 R^{\delta-1} h(\theta) \le A(\theta) \le  R^{\delta-1}
h(\theta)\,.
\end{equation}
We notice that $|A'(\theta)| \le (R+1)h(\theta) < 2Rh(\theta)$.
Combining this with the lower bound in \eqref{eq8.b}, we get the
differential inequality $ A'(\theta) \ge - 8
R^{2-\delta}A(\theta)$, whence $A(\theta) \ge A(0)
e^{-8R^{2-\delta}\theta}$.

To estimate $A(0)$ from below, recall that it equals $ \frac12 $
the flow of the field $-\nabla U$ through the interval $J$. Since
the length of $J$ is $2$, $A(0)$ cannot be less than the minimum
of the angular component of $-\nabla U$; i.e., $ A(0)\ge \frac12
R^{\delta-1}$. Thus, $A(\theta) \ge \frac12 R^{\delta-1}
e^{-8\theta R^{2-\delta}}$.

Now, using the upper bound in \eqref{eq8.b}, we get
\[
h(\theta) \ge R^{1-\delta}A(\theta) \ge e^{-8\theta R^{2-\delta}}
> e^{-2R^{2-\delta}}\,, \qquad \text{for}\quad  0 \le \theta \le \frac15\,.
\]
Hence the lemma. \hfill $\Box$

\medskip We can replace the function $F$ by its analytic
perturbation $F + H$ with $H$ satisfying $\left|\frac HF\right|\le
R^{-2}$ in the annulus $R-2\le |z|\le R+2$. After this
perturbation, the gradient field still satisfies the conditions
(i) and (ii), and the previous lemma applies to the new gradient
flow. The next lemma gives the lower bound for the probability of
the event that a G.E.F. $f$ is such a perturbation.
\begin{lemma}\label{lemma8.3}
If $0 < \delta < 1$ and $R\gg 1$, then $\displaystyle \Pr {
\max_{R-2 \le |z| \le R+2} \left| \frac{f}{F}(z) - 1 \right| \le
\frac1{R^2} } \ge e^{-C R^{1+\frac32 \delta}}$.
\end{lemma}
{\em Proof:} The proof we give is very similar to that of
Lemma~\ref{lemma8.1}. Actually, we estimate from below the
probability of the smaller event that $|f-F|\le e^{-R^\delta}|F|$
everywhere in the annulus $R-2\le |z|\le R+2$. Note that, in this
annulus, $\displaystyle |F(z)|\ge \frac12 e^{-2R^\delta}
\frac{|z|^n}{\sqrt{n!}}$.

First, we replace the exponent $e^{zR^{\delta-1}}$ by its Taylor
polynomial of degree $M = [20 R^\delta] $ in the disk $|z|\le 2R$.
It is easy to check that for $m\ge M$ and $|z|\le 2R$, the $m$-th
term in the Taylor expansion of the function $e^{zR^{\delta-1}}$
is bigger than twice the $m+1$-st term.  Hence the absolute value
of the tail that starts with the $M+1$-st term does not exceed the
absolute value of the $M$-th term. In particular, the relative
error we've made discarding the tail is at most
\[
e^{2R^\delta} \frac{(2R^\delta)^M}{M!} \le e^{2R^\delta} \left(
\frac{2e R^\delta}{M}\right)^M < e^{-10R^\delta}\,.
\]
Hence, for $|z|\le 2R$,
\[
|F(z) - P_M(z)| < e^{-10R^\delta} |F(z)|\,,
\]
where
\[
P_M(z) = \frac{z^n}{\sqrt{n!}} e^{-R^\delta} \sum_{m=0}^M
\frac{R^{\delta m}}{m!} \left( \frac{z}{R} \right)^m =
\sum_{m=0}^M a_m \frac{z^{n+m}}{\sqrt{(n+m)!}}
\]
is the Taylor polynomial of $F$. Note that
\[
a_m = \left( R^{-m} \sqrt{\frac{(n+m)!}{n!}}\right) \cdot \left(
e^{-R^\delta} \frac{R^{\delta m}}{m!}\right)\,.
\]
The second factor on the RHS is less than $1$. If $R\gg 1$, then
the first factor does not exceed $2$:
\[
\sqrt{\frac{(n+1)(n+2)\,...\, (n+m)}{n^m}} < \left( 1 +
\frac{M}{R^2}\right)^{M/2} < e^{\frac12 (M/R)^2} < 2\,.
\]
Thus, $0<a_m<2$.

Note that
\[
\left| \frac{f(z)}{F(z)} - 1 \right| \le \left|
\frac{f(z)-P_M(z)}{F(z)} \right| + \left| \frac{F(z)-P_M(z)}{F(z)}
\right|\,,
\]
and that we've already estimated the second term on the right-hand
side. We write
\begin{multline*}
\left| \frac{f(z)-P_M(z)}{F(z)} \right|
\\
\le 2e^{2R^\delta} \max_{R-2\le |z|\le R+2} \left[ \sum_{n \le k
\le n+M} |\xi_k-a_{k-n}| \sqrt{\frac{n!}{k!}} |z|^{k-n}+
\sum_{k\ne n, n+1, ..., n+M} |\xi_k| \sqrt{\frac{n!}{k!}}
|z|^{k-n} \right]
\end{multline*}
and show that with probability at least $e^{-C R^{1+\frac32 \delta
}}$ the maximum of the brackets on the right-hand side does not
exceed $C e^{-4R^\delta}$.

\medskip
We start with the first sum and demand that
\[
|\xi_{n+m}-a_m| < e^{-40R^\delta}, \quad m=0, 1,\, ...\,, M\,.
\]
The probability of this event is not less than $ \left( c
e^{-80R^\delta}\right)^{M+1} > e^{-C R^{2\delta}} >
e^{-CR^{1+\frac32 \delta}}$. For $|z|\le R+2$, $R\gg 1$, and $k\ge
n=R^2$, we have $\displaystyle \frac{|z|^{k+1}}{\sqrt{(k+1)!}} \le
2 \frac{|z|^k}{\sqrt{k!}}$. Hence $\displaystyle
\sqrt{\frac{n!}{k!}} |z|^{k-n} \le 2^{k-n} $, and the sum we are
estimating does not exceed $ 2^{M+1} e^{-40R^\delta} <
e^{-20R^\delta}$.

\medskip The second sum in the brackets does not exceed
\begin{equation}\label{eq9.*}
\sum_{0\le k < n} |\xi_k| \sqrt{\frac{n!}{k!}} (R-2)^{k-n}+
\sum_{k > n+M} |\xi_k| \sqrt{\frac{n!}{k!}} (R+2)^{k-n}\,.
\end{equation}
We estimate from below the probability that the first sum in
\eqref{eq9.*} is less than $3e^{-4R^\delta}$. The estimate for the
second sum is in the same spirit (cf. proof of
Lemma~\ref{lemma8.1}) and we omit it. We choose a large constant
$A\gg 1$ and split the first sum in \eqref{eq9.*} into two:
$\displaystyle S_1 = \sum_{0\le k < n-A\sqrt{Mn}}$ and
$\displaystyle S_2 = \sum_{n-A\sqrt{Mn} \le k < n}$.

\medskip
As in the proof of Lemma~\ref{lemma8.1}, we apply
Lemma~\ref{lemma2.1} to estimate the sum
$S_1$. For this, we need to estimate the sum
\begin{multline*}
\sum_{0\le k < n-A\sqrt{Mn}} \sqrt{n(n-1)\,...\, (n-(n-k-1))}
(R-2)^{k-n}
\\
=  \sum_{A\sqrt{Mn}< m \le n} \sqrt{n(n-1)\,...\,(n-(m-1))}
(R-2)^{-m}\,.
\end{multline*}
The $m$-th term of the sum on the right-hand side equals
\[
\sqrt{\left( 1-\frac1{n}\right) \left( 1-\frac2{n}\right)\, ...\,
\left( 1-\frac{m-1}{n}\right) } \left( 1 -
\frac2{R}\right)^{-m}\,.
\]
Using inequalities $1-\xi \le e^{-\xi}$, $0\le\xi\le 1$, and
$(1-\xi)^{-1} \le e^{2\xi}$, $0\le\xi \le \frac12$, we bound the
last expression by
\[
e^{\frac{4m}{R}} \prod_{j=0}^{m-1} e^{- \frac12 \frac{j}{n}} =
e^{\frac{4m}{R} - \frac{(m-1)m}{4n}} \le Ce^{-\frac{m^2}{8R^2}}\,,
\qquad 1\le m \le n\,.
\]
Then the sum we are estimating does not exceed
\[
C \int_{cAR^{1+ \frac{\delta}2}}^\infty e^{-\frac{t^2}{8R^2}}\, dt
\le CR \int_{cAR^{\delta/2}} e^{-t^2}\, dt \le CR e^{-cA^2
R^\delta} \le e^{-10R^\delta}\,,
\]
if $A$ is big enough.

Then, according to Lemma~\ref{lemma2.1}, the probability that
\[
\sum_{0\le k < n-A\sqrt{Mn}} |\xi_k| \sqrt{\frac{n!}{k!}}
(R-2)^{k-n} > e^{-4R^\delta}
\]
has a double exponential decay. We conclude modestly that $S_1 \le
e^{-4R^\delta}$ with probability at least $\frac12$.

\medskip
Now, we look at the sum $S_2$. In this case, we demand that
\[
|\xi_k|<e^{-10R^\delta}\,, \quad n-A\sqrt{Mn} \le k < n\,
\]
The probability of this event is not less than
\[
\left( \frac12 e^{-10R^\delta}\right)^{2(A\sqrt{Mn}+1)} \ge e^{-C
R^\delta \cdot R^{1+\frac{\delta}2}} = e^{-CR^{1+\frac32
\delta}}\,.
\]
Then the sum $S_2$ does not exceed
\begin{multline*}
 e^{-10R^\delta} \sum_{n -
A\sqrt{Mn}\le k < n } \sqrt{\frac{n!}{k!}} (R-2)^{k-n}
\\
= e^{-10R^\delta} \sum_{1\le m \le A\sqrt{Mn}} \sqrt{n(n-1)\,...\,
(n-(m-1))} (R-2)^{-m}\,.
\end{multline*}
We know from the discussion above that each term of the latter sum
is bounded by a constant. Hence
\[
S_2 \le CA\sqrt{Mn} e^{-10R^\delta} \le e^{-9R^\delta}\,,
\]
if $R$ is big enough. This completes the estimate of expression
\eqref{eq9.*} and proves the lemma. \hfill $\Box$

\medskip Now, let us fix the variables $\xi_k$ such that the function
$f(z)$ is a small perturbation of $F(z)$. For this function $f$,
we consider the corresponding ``tail'' $\mathcal G$. If $z$
belongs to the set $\mathcal G\cap \{z\colon
\frac{\pi}2-\frac1{10} < \arg z < \frac{\pi}2\}$ (the area of this
set was estimated in Lemma~\ref{lemma8.2}), then the trajectory
$\Gamma_z$ must traverse the whole set $\mathcal G \cap \{z\colon
\frac{\pi}2 < \arg z < \frac{\pi}2 + \frac1{10}\}$ before it hits
the radial interval $J$. Hence we expect that for such $z$'s the
distance from $z$ to its sink $a_z$ is comparable with $R$. We use
this idea to prove the following lemma.

\begin{lemma}\label{lemma8.4}
Suppose $R\gg 1$. With probability at least $e^{-C R^{1+\frac32
\delta}}$,
\[
m_2 \big\{z\in D(iR, R)\colon |z-a_z|\ge \frac{R}{100} \big\} \ge
e^{-C R^{2-\delta}}\,.
\]
\end{lemma}
\begin{figure}[h]
\begin{center}
\input{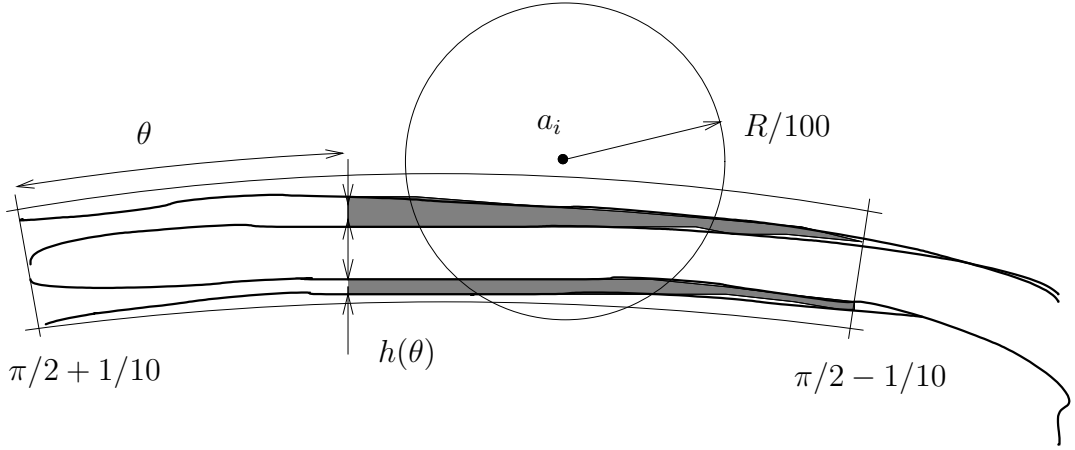}
\end{center}
\caption{The ``tail'' $\mathcal G_i$. The grey area equals
$A_i(\theta)$}
\end{figure}

\par\noindent{\em Proof:} After the trajectories from the tail
$\mathcal G$ leave the sector $\mathcal D$, they are attracted by
some of the zeroes of the function $f$. Let $a_1,\, ...\,, a_N$ be
the zeroes of $f$ that lie in the disk $2R\D$ and attract these
trajectories, and let $\mathcal G_i$ be the corresponding tails.
We discard the event $N\ge 100R^2$ since, by Theorem~2 in
\cite[part~III]{ST}, its probability is bounded by $e^{-CR^4}$
which is much less than $e^{-CR^{1+\frac32 \delta}}$. Hence we
assume that $N\le 100R^2$.

Let $A_i(\theta)$ be the area of the tail $\mathcal G_i\cap
\{\frac{\pi}2 - \frac1{10} < \arg z < \frac{\pi}2 + \frac1{10} -
\theta \}$, and let $h_i(\theta)$ be the length of the radial
section of $\mathcal G_i$ by the ray $\big\{ \arg z = \frac{\pi}2
+ \frac1{10} - \theta \big\}$; let $A_0(\theta)$, $h_0(\theta)$ be
the similar quantities that correspond to the trajectories
attracted by zeroes of $f$ lying outside the disk $2R\D$. By
Lemma~\ref{lemma8.2},
\[
\sum_{i=0}^N A_i(\tfrac1{10}) = m_2 \big(\mathcal G\cap \{z\colon
\tfrac{\pi}2-\tfrac1{10} < \arg z < \tfrac{\pi}2\} \big) \ge
e^{-CR^{2-\delta}}\,,
\]
thereby, $A_i(\frac1{10}) \ge e^{-C R^{2-\delta}}$ for some $i$.

If $i=0$, we are done: the points from the domain corresponding to
$A_0$ are far from their sinks. If $i\ne 0$, then, as in the proof
of Lemma~\ref{lemma8.2}, $h_i(\theta) \ge R^{1-\delta} A_i(\theta)
\ge e^{-C R^{2-\delta}}$ for $0 \le \theta \le \frac1{10}$. Hence,
after deleting the disk $D(a_i, \frac{R}{100})$, we still have a
set of points within $\mathcal D \subset D(iR, R)$ of area at
least $e^{-C R^{2-\delta}}$ that are attracted to $a_i$. This
proves the lemma. \hfill $\Box$

\medskip Now, we apply the same ``averaging trick'' that we've
already used in the proof of the upper bound for the distance to
the sink. Consider the (random) set $\mathcal C = \left\{z\colon
|z-a_z|\ge \frac{R}{100} \right\}$ and the event $\Omega^* =
\big\{ m_2 ( \mathcal C\cap D(z, R) ) \ge e^{-C R^{2-\delta}}
\big\}$. The probability of this event was estimated in the
previous lemma (for convenience, we took there $z=iR$ but, due to
the translation invariance, the probability of $\Omega^*$ does not
depend on the choice of $z$).

We aim at estimating from below the probability $\Pr {z\in\mathcal
C}$. We have
\begin{multline*}
\pi R^2 \Pr { z \in\mathcal C} = \iint_{R\D} \Pr {z + w \in
\mathcal C }\, dm_2(w) \\
= \int_\Omega m_2 \left( \mathcal C \cap D(z, R) \right) \,
d\mathbb P \ge \int_{\Omega^*} m_2 \left( \mathcal C \cap D(z, R)
\right) \, d\mathbb P \\
= \Pr{\Omega^*}\, e^{-C R^{2-\delta}} \ge e^{-C R^{1+\frac32
\delta} -C R^{2-\delta}}\,.
\end{multline*}
It remains to put $\delta=\frac25$ to balance the exponents. We
are done. \hfill $\Box$

\section{Diameter of the core}

Given $z\in \C$, we show that the probability of the event $
\big\{ m_2 \big( B_z \setminus D(a_z, R) \big) > \epsilon \big\} $
behaves as $e^{-c R^4}$ when $R$ is sufficiently large.

\subsection{The upper bound}

Given $z\in \C$, we show that the probability of the event $
\big\{ m_2 \big( B _z \setminus D(a_z, R) \big) > \epsilon \big\}
$ cannot be bigger than  $e^{-c R^4}$ when $R \gg 1$.

We take a small positive $\eta$ depending on $\epsilon$ only and
assume that $U \le \eta R^2$ everywhere in the basin $B_z$. It is
not difficult to see that the probability of the opposite event
does not exceed $C e^{-c_\eta R^4}$. Indeed, the event
$\displaystyle \big\{ \max_{B_z} U > \eta R^2\big\}$ is contained
in the union of the events $\big\{ \operatorname{diam} (B_z) > R^4
\big\}$ and $\displaystyle \big\{ \max_{D(z, R^4)} U > \eta R^2
\big\}$. By the long gradient curve theorem, the probability of
the first event does not exceed $Ce^{-cR^4}$. By
Lemma~\ref{lemma4.1}, the probability of the second event does not
exceed $CR^8 e^{-ce^{2\eta R^2}}$.

Similarly, we also assume that $ U \ge -\eta R^2 $ everywhere in $
B_z \setminus D(a_z, \frac{R}2)$. The opposite event is contained
in the union of the events $\big\{ \operatorname{diam} (B_z) > R^4
\big\}$ and
\[
\Big\{ {\rm there\ exists\ a\ curve\ } \gamma\subset R^4\D {\rm\
with\ } \operatorname{diam}(\gamma)\ge \frac12 R {\rm\ such\ that\
} \max_\gamma U < -\eta R^2 \Big\}\,,
\]
and by Theorem~\ref{thm4.3}, the probability of the second event
is bounded by $e^{-c_\eta R^4}$.

Thus, discarding events of probability less than $e^{-c_\eta
R^4}$, we may assume that $\displaystyle \max_{B_z\setminus D(a_z,
R)} |U| \le \eta R^2$. Then, by our length and area estimate
(Proposition~\ref{lemma8.a}),
\[
m_2 \big( B_z \setminus D(a_z, R) \big) \le \pi \exp \left( -
\frac{(R/2)^2}{\eta R^2} \right) = \pi \exp \left(
-\frac1{4\eta}\right) < \epsilon
\]
if $\eta$ is sufficiently small. This proves the upper bound.
\hfill $\Box$

\subsection{The lower bound}

We fix a positive $ \kappa < \pi $ and consider the random set
$\mathcal C = \left\{z\colon m_2 \big( B_z \setminus D(a_z, R)
\big) \ge \kappa \right\} $. We need to estimate from below the
probability $\Pr {z \in \mathcal C} $, which does not depend on
the choice of $z$. We apply the averaging again, but this time we
average over the disk of radius $R^5$. We get
\[
\pi R^{10} \Pr {0\in\mathcal C} = \iint_{R^5\D} \Pr {w\in \mathcal
C}\, dm_2(w) = \int_\Omega m_2(\mathcal C\cap R^5\D )\, d\mathbb
P\,.
\]

Introduce the event $\Omega^*$ that the following two conditions
hold:
\begin{itemize}
\item[(i)] $\displaystyle \# \big( \mathcal Z_f \cap R\D \big)
\ge \frac{4\pi}{\pi-\kappa}\,  R^2 $;
\item[(ii)] there is no gradient curve connecting the circumferences $\{|z|=R\}$
and $\{|z|=R^5\}$.
\end{itemize}

The probability of the first event is not less than $e^{-CR^4}$.
This estimate can be derived using the same techniques as in
\cite[Part~III]{ST} and in \cite{Krishnapur}, though it was not
explicitly proved in these papers. To get this estimate, denote by
$m$ the least integer that is not less than $\displaystyle
\frac{4\pi}{\pi-\kappa}\, R^2$, and estimate from below the
probability that
\[
\left| \xi_m \frac{z^m}{\sqrt{m!}} \right| > \left| f(z) - \xi_m
\frac{z^m}{\sqrt{m!}} \right|
\]
everywhere on the circumference $\{ |z| = R^2 \}$. We skip the
estimate since it repeats the one used in the proof of Theorem~3
in \cite{Krishnapur}.

Next, by the long gradient curve theorem, the probability that the
second event does not hold is less than $e^{-cR^5}$. Hence $\Pr
{\Omega^*} \ge e^{-CR^4}$.

Now, assuming that $\Omega^*$ happens, we can easily give a lower
bound for the area of the set $\mathcal C \cap R^5\D$. Actually,
we need to find only {\em one basin} $B(a)$ with $|a|\le R$ and
$m_2\big( B(a)\setminus 2R\D \big) \ge \kappa$. Then, by
assumption (ii), this basin lies within the disk $R^5\D$. Thereby,
$m_2 \big( \mathcal C \cap R^5\D \big) \ge \pi $, and we are done:
\[
\Pr {0\in \mathcal C} \ge  R^{-10} \Pr {\Omega^*} \ge
ce^{-CR^4}\,.
\]

To find a basin $B(a)$ with $|a|\le R$ and $m_2\big( B(a)\setminus
2R\D \big) \ge \kappa$, we do a simple counting. Consider the
basins $B(a)$ with $|a|\le R$ but $m_2 \big( B(a)\cap 2R\D \big)
> \pi - \kappa $. Let $N$ be the number of such basins.
Comparing the areas, we get
\[
4\pi R^2 = m_2(2R\D) \ge \sum_a m_2(B(a)\cap 2R\D) > (\pi-\kappa)
N\,;
\]
that is, $N<4\pi (\pi-\kappa)^{-1} R^2$. Hence, by assumption (i),
there is at least one basin $B(a)$, with $|a|\le R$ and
$m_2(B(a)\setminus 2R\D) \ge \kappa$. This finishes the proof.
\hfill $\Box$

\section{Modified basins}

In this section, we prove the remaining Theorem~\ref{th.cutoff}.
First, we describe a deterministic algorithm that ``improves''
partitions of the plane into domains of equal areas by cutting off
the tentacles of the basins and re-allocating them closer to the
sinks. Then we'll prove the probabilistic estimates for the sizes
of the modified basins of our random partition.

\subsection{Cutting off the tentacles}

Suppose we are given a partition of the plane $\displaystyle
\mathbb C = \cup_i E_i$ into bounded open domains of equal area,
say $\pi$, with marked points $c_i$, the ``centers'' of $E_i$. Let
\[
R_i = \inf\{R\colon E_i\subset D(c_i, R)\}\,.
\]
Clearly, $R_i \ge 1$.

Given $\epsilon\in (0, 1)$, we choose the least $r_i$ satisfying the
condition
\[
m_2\left( E_i\setminus D(c_i, r_i)\right) \le \frac1{A R_i^3}
\]
with $A=10^4 \epsilon^{-1}$ and define the ``{\em kernel}''
$\displaystyle K_i = E_i\cap D(c_i, r_i) $ and the ``{\em
tentacle}'' $T_i = E_i \setminus D(c_i, r_i)$ of the domain $E_i$.
Note that $m_2 T_i < 10^{-4}\epsilon$. It is worth mentioning that
this definition of the tentacle differs from the one we used in
Section~8.2. Later on, the factor $R_i^{-3}$ will help us to avoid
large tangles of different tentacles.

\begin{proposition}\label{prop11}
Given $\epsilon>0$, there exist open pairwise disjoint sets $E_i'$
with the following properties:

\smallskip\par\noindent (i) $m_2 E_i' = \pi$;

\smallskip\par\noindent (ii) $\displaystyle \C = \bigcup_{i} E_i'$
(up to a set of measure $0$);

\smallskip\par\noindent (iii) $m_2 \big(E_i\cap E_i'\big) \ge \pi - \epsilon$;

\smallskip\par\noindent (iv) $E_i' \subset D(c_i, r_i + \sqrt{5})$.
\end{proposition}

\medskip This proposition is useful when some of the domains $E_i$
have long tentacles; that is, $r_i \ll R_i$. The sets $E_i$ may be
assumed only measurable. Then the resulting sets $E_i'$ will be
measurable too.

\medskip\par\noindent{\em Proof of Proposition~\ref{prop11}:}
Split the plane $\C$ into standard unit squares. Suppose that $Q$
is one of them. First, we check that the union of the tentacles
$T_i$ can cover only a small portion of the square $Q$:
\begin{lemma}\label{lemma11.1}
\[
m_2 \left( Q\cap \left( \cup_i T_i \right) \right) \le \frac1{10}
\epsilon \,.
\]
\end{lemma}

\par\noindent{\em Proof: } If the domain $E_i$ with $R_i\le R$
intersects the square $Q$, then $E_i$ is contained in the square
with side length $4R+1$ homothetic to $Q$. Hence, comparing the
areas, we note that
\[
N_Q(R) \stackrel{\rm def} = \# \left\{ i\colon E_i \cap Q \ne
\varnothing, R_i \le R \right\} \le \frac1{\pi} (4R+1)^2\,.
\]
Thus
\begin{multline*}
\sum_{i\colon E_i\cap Q\ne \varnothing} m_2 T_i = \frac1{A}
\sum_{i\colon E_i\cap Q\ne \varnothing}
\frac1{R_i^3} \\
= \frac3{A} \int_1^\infty \frac{N_Q(R)}{R^4}\, dR \le \frac3{\pi
A} \int_1^\infty \frac{(4R+1)^2}{R^4}\, dR \le \frac1{10}
\epsilon\,.
\\ \hfill \Box
\end{multline*}

\medskip Now, let $\widehat{E}_i$ be a minimal square that is a
union of several standard unit squares and that contains the set
$E_i$.

\begin{lemma}\label{lemma11.2}
\[
\sum_{ i\colon Q\subset\widehat{E}_i } m_2 T_i \le \frac1{10}
\epsilon\,.
\]
\end{lemma}

\par\noindent{\em Proof: } Comparing the areas, we see that
\[
\# \left\{i\colon Q\subset \widehat{E}_i, R_i\le R  \right\} \le
\frac1{\pi} (4R+3)^2\,.
\]
The rest is the same as in the previous lemma. \hfill $\Box$

\medskip Let $\widehat{Q}_i$ be the square that contains the
center $c_i$ of $E_i$ (if $c_i$ lies on the grid, it does not
matter which one of several squares containing $c_i$ to choose).
For each pair $(i, Q)$ with $Q\subset \widehat{E}_i\setminus
\widehat{Q}_i$, we choose a ``{\em storage}'' $S_i(Q) \subset Q
\cap \left( \cup_j K_j \right) $ according to the following rules:
\begin{itemize}
\item[(a)] $ m_2 S_i(Q)  = m_2 T_i $;
\item[(b)] for different $i$'s, the storages $S_i(Q)$ are mutually
disjoint;
\item[(c)] for each pair $(i, Q)$, the area of the storage
$S_i(Q)$ is distributed between the kernels $K_j\cap Q$
proportionally to their areas; i.e.,
\[
m_2\left( S_i(Q)\cap K_j\right) = m_2\big( S_i(Q) \big)
\,\frac{m_2(Q\cap K_j)}{\sum_l m_2(Q\cap K_l)}\,.
\]
\end{itemize}
By Lemma~\ref{lemma11.2}, the total area within $Q$ that we need
to allocate to all the storages does not exceed $
\frac{\epsilon}{10} $, while by Lemma~\ref{lemma11.1}, the area of
$Q \cap \left( \cup_j K_j \right)$ is not less than $ 1 -
\frac{\epsilon}{10}$. Hence, we can meet the requirements (a) and
(b). The requirement (c) does not impose any additional
restriction.

\medskip Now we describe the cut-off algorithm. It consists of
countably many parallel and independent of each other processes.
During the $i$-th process, for each square
$Q\subset\widehat{E_i}$, the piece of the tentacle $T_i\cap Q$ is
re-allocated to some centers $c_l$ such that $K_l\cap Q \ne
\varnothing$. At the same time, some subsets of $K_l\cap S_i(Q)$
are re-allocated to some centers $c_m$ whose kernels $K_m$
intersect one of the squares neighbouring $Q$.
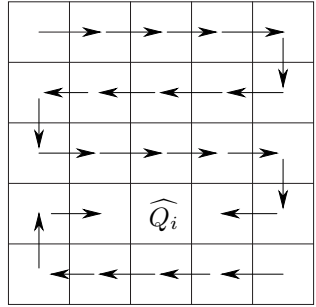
\begin{figure}[h]
\begin{center}
\setlength{\unitlength}{0.00069991in}
\begingroup\makeatletter\ifx\SetFigFont\undefined%
\gdef\SetFigFont#1#2#3#4#5{%
  \reset@font\fontsize{#1}{#2pt}%
  \fontfamily{#3}\fontseries{#4}\fontshape{#5}%
  \selectfont}%
\fi\endgroup%
{\renewcommand{\dashlinestretch}{30}
\begin{picture}(2274,2289)(0,-10)
\path(12,2262)(2262,2262)(2262,12)
	(12,12)(12,2262)
\path(12,1812)(2262,1812)
\path(12,1362)(2262,1362)
\path(12,912)(2262,912)
\path(12,462)(2262,462)
\path(462,2262)(462,12)
\path(912,2262)(912,12)
\path(1362,2262)(1362,12)
\path(1812,2262)(1812,12)
\path(237,2037)(687,2037)
\blacken\path(567.000,2007.000)(687.000,2037.000)(567.000,2067.000)(603.000,2037.000)(567.000,2007.000)
\path(732,2037)(1137,2037)
\blacken\path(1017.000,2007.000)(1137.000,2037.000)(1017.000,2067.000)(1053.000,2037.000)(1017.000,2007.000)
\path(1182,2037)(1542,2037)
\blacken\path(1422.000,2007.000)(1542.000,2037.000)(1422.000,2067.000)(1458.000,2037.000)(1422.000,2007.000)
\path(1587,2037)(2037,2037)
\blacken\path(1917.000,2007.000)(2037.000,2037.000)(1917.000,2067.000)(1953.000,2037.000)(1917.000,2007.000)
\path(2037,1992)(2037,1632)
\blacken\path(2007.000,1752.000)(2037.000,1632.000)(2067.000,1752.000)(2037.000,1716.000)(2007.000,1752.000)
\path(2037,1587)(1632,1587)
\blacken\path(1752.000,1617.000)(1632.000,1587.000)(1752.000,1557.000)(1716.000,1587.000)(1752.000,1617.000)
\path(1587,1587)(1137,1587)
\blacken\path(1257.000,1617.000)(1137.000,1587.000)(1257.000,1557.000)(1221.000,1587.000)(1257.000,1617.000)
\path(1092,1587)(687,1587)
\blacken\path(807.000,1617.000)(687.000,1587.000)(807.000,1557.000)(771.000,1587.000)(807.000,1617.000)
\path(597,1587)(282,1587)
\blacken\path(402.000,1617.000)(282.000,1587.000)(402.000,1557.000)(366.000,1587.000)(402.000,1617.000)
\path(237,1542)(237,1182)
\blacken\path(207.000,1302.000)(237.000,1182.000)(267.000,1302.000)(237.000,1266.000)(207.000,1302.000)
\path(237,1137)(642,1137)
\blacken\path(522.000,1107.000)(642.000,1137.000)(522.000,1167.000)(558.000,1137.000)(522.000,1107.000)
\path(687,1137)(1137,1137)
\blacken\path(1017.000,1107.000)(1137.000,1137.000)(1017.000,1167.000)(1053.000,1137.000)(1017.000,1107.000)
\path(1182,1137)(1542,1137)
\blacken\path(1422.000,1107.000)(1542.000,1137.000)(1422.000,1167.000)(1458.000,1137.000)(1422.000,1107.000)
\path(1632,1137)(1992,1137)
\blacken\path(1872.000,1107.000)(1992.000,1137.000)(1872.000,1167.000)(1908.000,1137.000)(1872.000,1107.000)
\path(2037,1092)(2037,732)
\blacken\path(2007.000,852.000)(2037.000,732.000)(2067.000,852.000)(2037.000,816.000)(2007.000,852.000)
\path(1992,687)(1587,687)
\blacken\path(1707.000,717.000)(1587.000,687.000)(1707.000,657.000)(1671.000,687.000)(1707.000,717.000)
\path(2037,237)(1587,237)
\blacken\path(1707.000,267.000)(1587.000,237.000)(1707.000,207.000)(1671.000,237.000)(1707.000,267.000)
\path(1497,237)(1137,237)
\blacken\path(1257.000,267.000)(1137.000,237.000)(1257.000,207.000)(1221.000,237.000)(1257.000,267.000)
\path(1092,237)(687,237)
\blacken\path(807.000,267.000)(687.000,237.000)(807.000,207.000)(771.000,237.000)(807.000,267.000)
\path(642,237)(327,237)
\blacken\path(447.000,267.000)(327.000,237.000)(447.000,207.000)(411.000,237.000)(447.000,267.000)
\path(237,282)(237,687)
\blacken\path(267.000,567.000)(237.000,687.000)(207.000,567.000)(237.000,603.000)(267.000,567.000)
\path(327,687)(687,687)
\blacken\path(567.000,657.000)(687.000,687.000)(567.000,717.000)(603.000,687.000)(567.000,657.000)
\put(1047,597){\makebox(0,0)[lb]{\smash{{\SetFigFont{10}{12.0}{\rmdefault}{\mddefault}{\updefault}$\widehat{Q_i}$}}}}
\end{picture}
}
\end{center}
\caption{The square $\hat{E_i}$  and two sequences of unit
squares}
\end{figure}

We split the unit squares from $ \widehat{E}_i\setminus
\widehat{Q}_i $ into two disjoint sequences $\{ Q_1,  Q_2,\,
...\,, Q_{m_1} \}$ and  $\{ Q_{m_1+1},  Q_{m_1+2},\, ...\,,
Q_{m_2} \}$ such that in each sequence any two consecutive squares
$Q_l$ and $Q_{l+1}$ have a common boundary side, and the last
squares $Q_{m_1}$, $Q_{m_2}$ of each sequence have a common
boundary side with the square $\widehat{Q}_i$ (see Figure~9).
\begin{figure}[h]
\begin{center}
\input{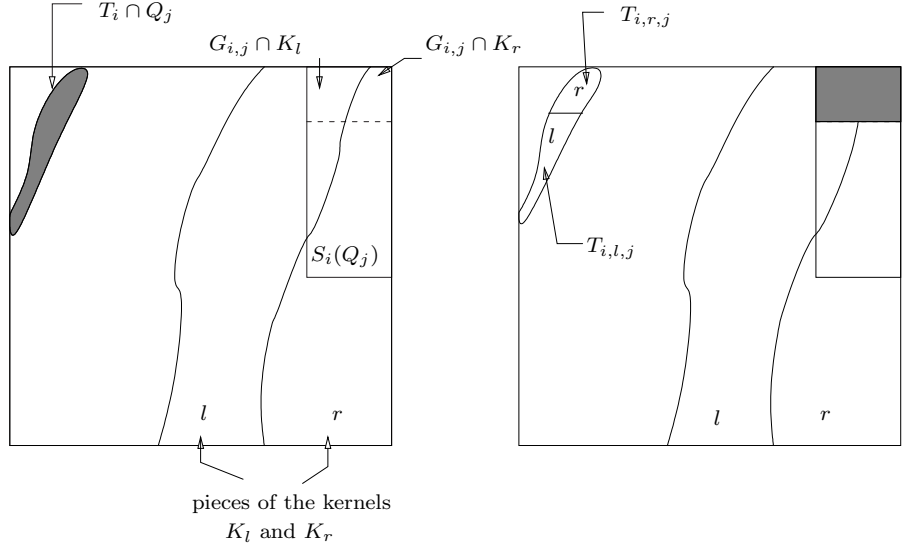}
\end{center}
\caption{Putting the grey area $T_i\cap Q_j$ into the storage
$S_i(Q_j)$}
\end{figure}

Let us call $T_i\setminus \widehat{Q}_i$ the {\em grey area}.
First, for each $j$, $1\le j \le m_1$, we swap the set $T_i\cap
Q_j$ with a part of the storage $S_i(Q_j)$. More precisely, we
\begin{itemize}
\item[i.] choose parts of the storages $G_{i, j} \subset S_i(Q_j)$ with
$m_2 G_{i, j} = m_2 (T_i\cap Q_j)$;
\item[ii.] decompose the tentacle $T_i\cap Q_j$ into disjoint union of
subsets $T_{i, l, j}$, $l\ne i$, with $m_2 T_{i, l, j} = m_2(G_{i,
j} \cap K_l)$;
\item[iii.] for $1\le j \le m_1$, re-allocate the grey area from
$T_i\cap Q_j$ to $G_{i, j}$;
\item[iv.] for each $l\ne i$, remove the set $\displaystyle K_l
\cap \bigcup_{1\le j \le m_1} G_{i, j}$ from $E_l$, and
re-allocate the set $\displaystyle \bigcup_{1\le j \le m_1} T_{i,
l, j}$ of equal measure to $E_l$.
\end{itemize}
Now, the grey area occupies some parts of the storages $S_i(Q_j)$.

\begin{figure}[h]
\begin{center}
\input{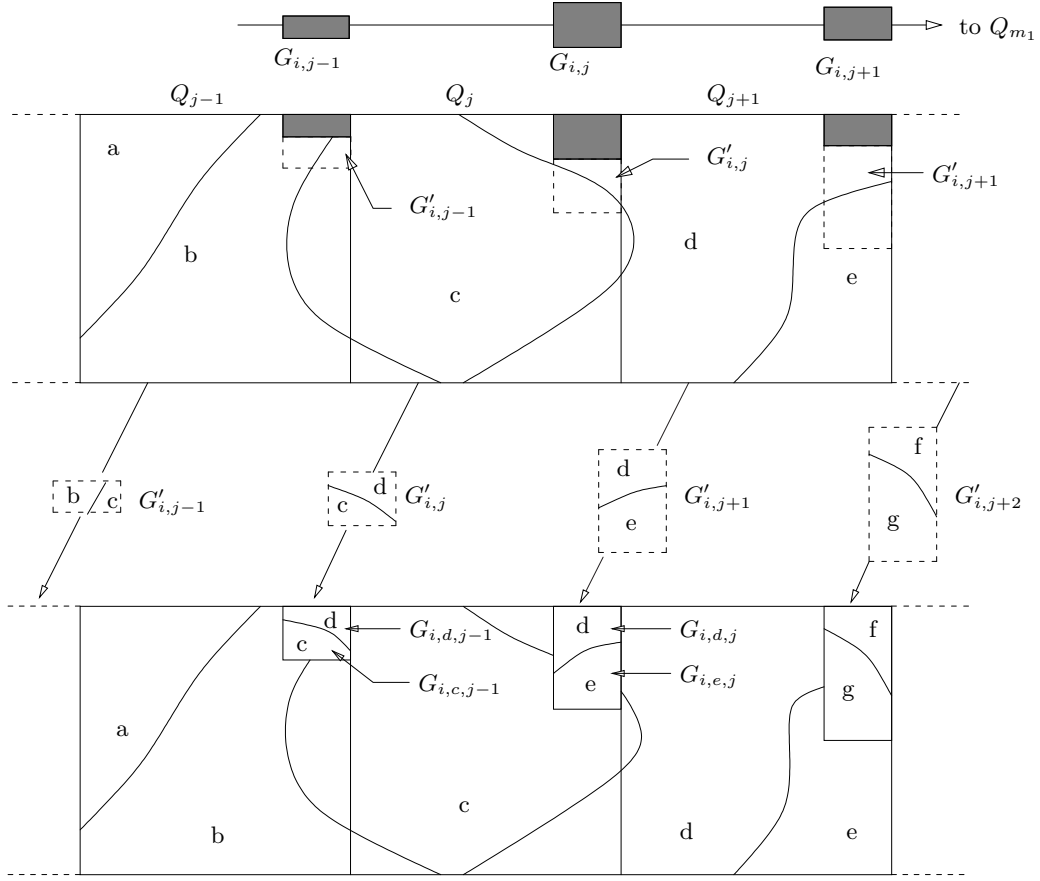}
\end{center}
\caption{Moving the grey area from the storages $S_i(Q_j)$ 
to the storage $S_i(Q_{m_1})$}
\end{figure}

\medskip

At the next step, starting with the square $Q_1$, square after
square, we move the grey area from $S_i(Q_j)$ to $S_i(Q_{j+1})$,
until the whole grey area appears in the last storage
$S_i(Q_{m_1})$ of the sequence of squares we are traversing. After
that, we allocate the grey area to the center $c_i$. More
formally,
\begin{itemize}
\item[i.] for $2\le j \le m_1$, we choose subsets $G_{i, j}' \subset
S_i(Q_j) \setminus G_{i, j}$ such that $\displaystyle m_2 G_{i,
j}' = \sum_{k=1}^{j-1} m_2 G_{i, k}$, and set $G_{i, j}'' = G_{i,
j}' \cup G_{i, j}$, $G_{i, 1}'' = G_{i, 1}$;
\item[ii.] for $1\le j \le m_1-1$, we decompose the sets $G_{i, j}''$ into
disjoint union of subsets $G_{i, l, j} $ with $m_2 G_{i, l, j} =
m_2 \left( G_{i, j+1}' \cap K_l \right)$;
\item[iii.] within each $E_l$, $l\ne i$, we replace the set
$\displaystyle K_l \cap \bigcup_{1\le j\le m_1-1} G'_{i, j+1 } $
by the set $\displaystyle \bigcup_{1\le j \le m_1-1} G_{i, l, j}$
of equal measure.
\item[iv.] In the end, the tentacles $T_i\cap Q_j$, $1\le j \le m_1$,
are cut off from $E_i$, and the set $G_{i, m_1 }''$ with $
\displaystyle m_2 G_{i, m_1 }'' = \sum_{1\le j \le m_1} m_2 \left(
T_i\cap Q_j\right) $ is added to $E_i$.
\end{itemize}
Then we apply the same process to the second sequence of squares
$\{Q_{m_1+1}, \,...\,, Q_{m_2} \}$.

\medskip

Note that all points re-allocated during the $i$-th process will
appear either in $T_i$ or in one of the storages $S_i$. Hence, due
to the choice of the storages, these points are not displaced
during the other steps. We see that for different $i$'s the
processes are independent of each other.

We conclude that the new sets $E_i'$ are located in the
$\sqrt{5}$-neighbourhoods of the kernels $K_i$ and have the same
area as $E_i$. ($\sqrt{5}$ is the length of the diagonal of the
rectangle comprised of two adjacent standard squares.) By
construction,
\[
E_i\setminus E_i' \subset T_i \cup \bigcup_{(j, Q)} \left(
S_j(Q)\cap K_i \right)\,.
\]
Due to the choice of the storages and Lemmas~\ref{lemma11.1} and~
\ref{lemma11.2},
\[
\sum_{(j, Q)} m_2(S_j(Q)\cap K_i) \le 2 \sum_Q \left[
\sum_{j\colon Q\subset\widehat{E}_j} m_2 ( T_j ) \right] m_2(Q\cap
K_i) \le \frac{\pi\epsilon}5\,.
\]
Recall that $m_2 T_i \le 10^{-4} \epsilon$. Hence, for each $i$,
$m_2(E_i\setminus E_i')<\epsilon$. This proves
Proposition~\ref{prop11}. \hfill $\Box$

\subsection{Probabilistic estimate}

We fix $\epsilon >0$ and apply the cut-off algorithm  to the
basins $B(a)$. The sink $a$ is the ``center'' of $B(a)$, $R(a)$ is
the least number $R$ such that $B(a)\subset D(a, R)$. As above, we
set $A=10^4 \epsilon^{-1}$, and
\[
r(a) = \min \left\{r\colon m_2\left( B(a) \setminus D(a, r) \right)
\le \frac1{AR^3(a)} \right\}\,.
\]
Proposition~\ref{prop11} gives us the modified basins
$B'(a)\subset D(a, r(a)+\sqrt{5})$ satisfying conditions
(i)--(iii) of Theorem~\ref{th.cutoff}. Let $B'(\alpha)$ be the
modified basin with center at $\alpha $ that contains the origin.
Since $\operatorname{diam} B'(\alpha) \le 2\big( r(\alpha) +
\sqrt{5} \big)$, the proof of condition (iv) in
Theorem~\ref{th.cutoff} boils down to the estimate
\begin{equation}\label{eq11.a}
\Pr {r(\alpha) > R } \le e^{-cR^4/(\log R)^{3/2}}
\end{equation}
for $R\gg A$.

\begin{claim}\label{claim11.1}
$\Pr {|\alpha| \ge \frac12 R^4} \le e^{-cR^4}$.
\end{claim}

\par\noindent{\em Proof:}
Assume that $|\alpha|  \ge \frac12 R^4$. Since the origin lies at
the distance at most $\sqrt{5}$ from the basin $B(\alpha)$, we
know that there is a gradient curve that connects the
circumferences $\{|z|=\sqrt{5}\}$ and  $\{|z|=\frac12 R^4\}$. This
gradient curve connects the boundaries of the squares $Q(0,
\frac1{4\sqrt{2}} R^4)$ and $Q(0, \frac1{2\sqrt{2}} R^4)$. By the
long gradient curve theorem, the probability of this event is less
than $e^{-cR^4}$. \hfill$\Box$

\medskip Now, we prove \eqref{eq11.a}.
First, we suppose that $R(\alpha) > R^4$. In view of the claim, we
also assume that $|\alpha|< \frac12 R^4$. We cover  the disk $D(0,
\frac12 R^4)$ by a bounded number of standard squares $Q(w,
\frac1{2\sqrt{2}} R^4)$, and consider the square that contains the
point $\alpha$. We know that there is a gradient curve of diameter
$R(\alpha)$ that terminates at the sink $\alpha$. This gradient
curve must connect $\partial Q(w, \frac1{2\sqrt{2}} R^4)$ with
$\partial Q(w, \frac1{\sqrt{2}} R^4)$. By the long gradient curve
theorem, the probability of this event is less than $e^{-cR^4}$.
Hence, $\Pr {R(\alpha) > R^4 } < e^{-cR^4}$.

\medskip
Now, we suppose that $ R^4\ge R(\alpha) $. Set $M =
\frac1{52}\,\frac{R^2}{\log R}$. By Lemma~\ref{lemma4.1}, throwing
away an event of probability much less than $e^{-cR^4}$, we may
assume that $ U \le M $ everywhere in $D(0, 3R^4)$, in particular,
everywhere in $D(\alpha, R(\alpha) )$. By Claim~\ref{claim11.1},
we may assume that $|\alpha|< \frac12 R^4$. Hence, if
$\displaystyle \min_{B(\alpha)\setminus D(\alpha, \frac12 R)} U <
-M$, then the disk $R^4\D$ contains a curve of diameter at least
$\frac 12 R$ where $U<-M$. By Theorem~\ref{thm4.3}, the
probability of this event does not exceed $ e^{-cRM^{3/2}}$. Thus,
discarding the event of probability at most $e^{-cR^4/(\log
R)^{3/2}}$, we may assume that $\displaystyle |U|\le M$ in
$B(\alpha)\setminus D\big(\alpha, \frac12 R\big)$. Then by the
length and area estimate (Proposition~\ref{lemma8.a}), the area of
the set $B(\alpha)\setminus D(\alpha, R)$ cannot exceed
\[
\pi e^{-2 \frac{(R/2)^2}{2M}} = \pi R^{-13} < \frac1{AR^{12}} \le
\frac1{AR^3(\alpha)}\,,
\]
provided that $R\ge \pi A$. Hence, after the events described
above have been thrown away, we get $r(\alpha) \le R$.  Therefore,
the probability of the event $\{r(\alpha)>R\}$ does not exceed the
sum of probabilities of the events thrown away, and we are done.
\hfill $\Box$

\section{Discussion and  questions}

\subsection{Optimal transportation to the zero set of G.E.F.}
\begin{question}
Does there exists a transportation $T$ of the Lebesgue measure
$\frac1{\pi} m_2$ to the random zero set $\mathcal Z_f$ such that
the tails $\displaystyle \sup_{z\in\C} \Pr {|T(z)-z|>R}$ decay as
$e^{-cR^4}$ as $R\to\infty$?
\end{question}
Recall that the estimate $e^{-cR^4 (\log R)^{-1}}$ can be achieved
by modification of the proof in \cite[Part~II]{ST}. Note that, in
view of the lower bound for the ``hole probability'' $\Pr {
\mathcal Z_f \cap R\D = \varnothing} \ge ce^{-CR^4}$ proved in
\cite[Part~III]{ST}, one cannot get a better estimate than
$ce^{-CR^4}$.

\subsection{Length of the gradient curve and the travel time}

Given $z$, consider the gradient curve $\Gamma_z$ that passes
through the point $z$. Let $\ell_z$ be the length of the part of
the curve $\Gamma_z$ that starts at $z$ and terminates at $a_z$.

\begin{question}
Find the order of decay of the tails $\Pr {\ell_z >R}$ as
$R\to\infty$.
\end{question}

An interesting characteristic of the ``random landscape'' of the
potential $U$ is the time $\tau_z$ needed for the point $z$ to
roll down to the sink $a_z$ along the gradient curve $\Gamma_z$.
By analogy with some models from astrophysics, Michael Douglas
asked us about the {\em order of decay of the tails $\Pr {\tau_z >
t}$ as $t\to\infty$}. Since $\operatorname{div} (\nabla U) = -2$
everywhere on $\C\setminus \mathcal Z_f$, one can show using
Liouville's theorem that this probability {\em equals} $e^{-2t}$
(cf. Section~8). The length $l$ measured along the gradient curve
and the travel time $\tau$ are connected by relation
$\displaystyle \frac{dl}{d\tau} = |\nabla U |$. Since we know the
distribution of the gradient field (recall that $\displaystyle
\nabla U = \overline{\frac{f'(z)}{f(z)}} - z$), it looks tempting
to use this information  to simplify the proofs of our main
results and to achieve a better understanding of the properties of
the random partition.

\subsection{Statistics of the basins}

There are several interesting questions related to the statistics
of our random partition of the plane. We say that two basins are
neighbours if they have a common gradient curve on the boundary.
By $N_z$ we denote the number of basins $B$ neighbouring the basin
$B_z$. Clearly, $N_z$ equals the number of saddle points of the
potential $U$ connected with the sink $a_z$ by gradient curves.
Heuristically, since almost surely each saddle point is connected
with two sinks,
\[
\mathbb E N_z = 2\, \frac{{\rm mean\ number\ of\ saddle\ points\
per\ unit\ area}}{{\rm mean\ number\ of\ zeroes\ per\ unit\
area}}\,.
\]
Douglas, Shiffman and Zelditch proved in \cite{DSZ} that the mean
number of saddle points of $U$ per unit area is $\frac4{3\pi}$.
(They proved this for another closely related ``elliptic model''
of Gaussian polynomials. It seems that their proof also works for
G.E.F.) Hence the question:
\begin{question}
Prove that $\mathbb E N_z = \frac83$.
\end{question}

We are also interested in the behaviour of the tails of the random
variable $N_z$:
\begin{question}
Find the order of decay of $\Pr { N_z > N}$ as $N\to \infty$.
\end{question}

Another characteristic of the random partition is the number $M$
of basins that meet at the same local maximum. Taking into account
the result from \cite{DSZ}, we expect that its average equals $8$.
It is also interesting to look at the decay of the tails of $M$.
Probably, some lower bound can be extracted from the analysis of
perturbations of the polynomial $z^n-1$ similar to the one we did
in Section~9.

\subsection{The skeleton topology}

By the {\em skeleton} of the gradient flow we mean the connected
planar graph with vertices at local maxima of $U$ and edges
corresponding to the boundary curves of the basins. The graph may
have multiple edges and loops. Our question is
\begin{question}
Are there any non-trivial topological restrictions on finite parts
of the skeleton that hold almost surely?
\end{question}

There is an interesting  finite counterpart of this question.
Choose $N$ independent points $ a_1, \,...\,, a_N $ uniformly
distributed on the Riemann sphere $\widehat{\mathbb C} $ and
consider the gradient flow of the random spherical potential
\[
V(z) = \sum_{i} \log|z-a_i| - \frac{N}2 \log(1+|z|^2)\,.
\]
\begin{question}
Describe all possible skeletons of the gradient flow on $\widehat
{\mathbb C}$ of the potential $V$ that are realized with positive
probability.
\end{question}

\bigskip
\filbreak { \small
\begin{sc}
\parindent=0pt\baselineskip=12pt
\parbox{2.4in}{
F.N. and A.V.:\\
Department of Mathematics\\
Michigan State University\\
East Lansing, MI 48824, USA
\smallskip
\emailwww{fedja@math.msu.edu} {volberg@math.msu.edu} } \hfill
\parbox{2.1in}{
M.S.:\\
School of Mathematics\\
Tel Aviv University\\
Tel Aviv 69978, Israel
\smallskip
\emailwww{sodin@post.tau.ac.il} {} } \hfill
\end{sc}
}

\end{document}